\documentclass{article}

\usepackage{hyperref}
\usepackage{amsmath}
\usepackage{unicode-math}
\usepackage{amsthm}
\usepackage{tikz}
\usepackage{tikz-cd}

\usepackage[backend=biber, style=alphabetic]{biblatex}

\theoremstyle{definition}
\newtheorem{defin}{Definition}[section]

\theoremstyle{remark}
\newtheorem{rmk}[defin]{Remark}%[section]

\theoremstyle{remark}
\newtheorem{prop}[defin]{Proposition}%[section]

\theoremstyle{definition}
\newtheorem{thm}[defin]{Theorem}%[section]

\theoremstyle{remark}
\newtheorem{lm}[defin]{Lemma}%[section]

\theoremstyle{remark}
\newtheorem{ex}[defin]{Example}%[section]

\theoremstyle{definition}
\newcommand{\comment}[1]{}
\newcommand{\thmheader}[1]{\emph{#1}\label{#1}}
\newcommand{\colorformula}[2]{\color{#1}#2\color{black}}

\newcommand{\todo}[1]{}

\tikzset{
bullet/.style={circle, draw=black, fill,minimum size=4pt, inner sep=0pt, outer sep=0pt}%,
}

\usetikzlibrary{calc}

\title{Parametric Distributive Laws: \\ uniform monad composition}
\author{Lorenzo Perticone}
\date{\today}

\begin{document}

\maketitle

\begin{abstract}
Monads play an important role in both the syntax and semantics of modern functional programming languages (\cite{notions-of-computation}). The problem of combining them has been of profound interest at least since the 90s, and different approaches have been employed to tackle it: currently, the most prevalent notion is that of \textit{monad transformers} (\cite{syntactic-modular-denotational-semantics, monad-transformers-modular-interpreters}). We provide a novel abstract framework to describe such ``compositions'' which we call \textit{parametric distributive laws}.

Our description of this framework hinges on the theory of 2-categories and relies strongly on the construction of a left (1-) adjoint to the construction of monads internal to a 2-category: since distributive laws are monads in monads, we obtain a semi-strict ``walking distributive law'' by Gray-tensoring the walking monad with itself. By Gray-tensoring again, we can then construct a parametric variant thereof.

We demonstrate the applicability of such a framework by providing two concrete examples (involving the \texttt{Writer} and \texttt{Either} monads), explicitly describing morphisms of such parametric distributive laws, and showing how (an iterated version of) our construction can be employed to describe (parametric) iterated distributive laws (in the sense of \cite{iterated-distributive-laws}), motivating the appearance of the Yang-Baxter equations as coherences between the distributive laws involved.
\end{abstract}

\thanks{This research was funded by SSF, the Swedish Foundation for Strategic Research, grant number FUS21-0063}

\section*{Introduction}

The first structured uses of monads in computer science are due to Moggi \cite{notions-of-computation}, with the purpose of modeling different features (e.g. non-termination, side effects, nondeterminism) in the semantics of a programming language. The problem of composing them has been explored long before that, in contexts such as categorical algebra. Different authors have proposed different ways to describe such compositions; the main approaches seem to be \textit{distributive laws} (mainly employed in algebraic contexts) and \textit{monad transformers} (mainly used in computer science, e.g. the Haskell standard library).

\begin{center}
A \textit{monad transformer} for a monad $\mathcal{T}$ over the category $\mathcal{C}$ is a pointed endofunctor $(\mathbf{Tr}_{\mathcal{T}}, u)$ over the category of monads over $\mathcal{C}$, such that $\mathbf{Tr}_{\mathcal{T}}(\mathbf{Id}) = \mathcal{T}$. That is, an endofunctor together with a natural transformation
\[u : \mathbf{Id} \to \mathbf{Tr}_{\mathcal{T}} : \mathbf{Mnd}(\mathcal{C}) \to \mathbf{Mnd}(\mathcal{C})\]
\end{center}

In this framework, one considers $\mathbf{Tr}_{\mathcal{T}}(\mathcal{S})$ the ``composite monad'' of $\mathcal{T}$ and $\mathcal{S}$. 

\begin{center}
A \textit{distributive law} of a monad $\mathcal{T} = (T, \eta^T, \mu^T)$ over the monad $\mathcal{S} = (S, \eta^S, \mu^S)$ is a natural transformation $\gamma : S \circ T \to T \circ S$, subject to appropriate axioms.
\end{center}

In this framework, one considers $\mathcal{T} \circ \mathcal{S}$ the composite monad (the axioms we glossed over above guarantee that there is a well-defined multiplication for it).

What we propose in this paper is a generalization of the latter notion, allowing for the treatment of a \textit{class} of monads instead of just one such. Intuitively, one should think of that notion as a ``coherent family'' of distributive laws, in a sense that will be made precise in the next chapter. It also has the pleasant property of being easily generalized to the context of composing $n$ monads (rather than just two), as we'll see later.

For a fixed monad $\mathcal{T}$, our notion will be equivalent to fixing a coherent family of \textit{lifts} $\hat{\mathcal{S}}$ for monads $\mathcal{S}$ on $\mathcal{C}$, obtaining monads on $\mathcal{T}$'s Eilenberg-Moore category $\mathcal{C}^{\mathcal{T}}$, which correspond 1-1 with compositions (by composing the corresponding monadic adjunctions).

%% ## SECTION 1
\section{Monads in 2-category theory}

The notion of a \textit{monad over a category} $\mathcal{T} : \mathcal{C} \to \mathcal{C}$ is well-known to mathematicians (dating back to \cite{topologie-algebrique} and \cite{homotopy-theory-general-categories}; originally called \textit{standard constructions} or \textit{triples}) and computer scientists (the first reference we're aware of in this context is \cite{notions-of-computation}). We'll discuss a direct generalization, \textit{monads in a 2-category} (first introduced in \cite{introduction-to-bicategories}, where the name ``monad'' is first used). We'll start by introducing basic notions, notations (notably, string diagram notation) and results about 2-categories; our treatment will closely follow B\'enabou's \cite{introduction-to-bicategories} and Leinster's \cite{basic-bicategories}.

\begin{defin}\thmheader{2-category}

A 2-category $\mathcal{B}$ is given by
\begin{enumerate}
\item A class of \textit{0-cells} or \textit{objects}, denoted $X : \mathcal{B}$,
\item For any objects $X, Y : \mathcal{B}$, a (possibly large) \textit{hom category} $\mathcal{B}[X, Y]$ whose
\begin{itemize}
\item objects are called \textit{1-cells} or \textit{morphisms}
\item morphisms are called \textit{2-cells} or \textit{transformations}
\item composition $\circ$ is called \textit{vertical} composition
\end{itemize}
\item For any object $X : \mathcal{B}$, a 1-cell $\mathbf{Id}_{X} : \mathcal{B}[X, X]$ (sometimes denoted by $X$)
\item For any objects $X, Y, Z : \mathcal{B}$ a \textit{horizontal composition} functor 
\[\cdot : \mathcal{B}[Y, Z] \times \mathcal{B}[X, Y] \to \mathcal{B}[X, Z]\]
\item For any objects $X, Y : \mathcal{B}$, \textit{left} and \textit{right unitors} natural isomorphisms
\[\lambda : (Y \cdot -) \to \mathbf{Id}_{\mathcal{B}[X, Y]} : \mathcal{B}[X, Y] \to \mathcal{B}[X, Y]\]
\[\rho : (- \cdot X) \to \mathbf{Id}_{\mathcal{B}[X, Y]} : \mathcal{B}[X, Y] \to \mathcal{B}[X, Y]\]
\item For any objects $W, X, Y, Z : \mathcal{B}$, an \textit{associator} natural isomorphism
\[\alpha : ((- \cdot -) \cdot -) \to (- \cdot (- \cdot -)) : \mathcal{B}[Y, Z] \times \mathcal{B}[X, Y] \times \mathcal{B}[W, X] \to \mathcal{B}[W, Z]\]
\end{enumerate}
satisfying the triangle and pentagon identities (whenever well-typed)
\begin{center}\begin{tikzcd}[scale=0.5]
(g \cdot Y) \cdot f \arrow[dd, "\alpha"'] \arrow[dr, "\rho \cdot f"] &
& ((l \cdot h) \cdot g) \cdot f \arrow[r, "\alpha \cdot f"] \arrow[dd, "\alpha"] & (l \cdot (h \cdot g)) \cdot f \arrow[r, "\alpha"] & l \cdot ((h \cdot g) \cdot f) \arrow[dd, "l \cdot \alpha"] \\
& g \cdot f \\
g \cdot (Y \cdot f) \arrow[ur, "g \cdot \lambda"'] &
& (l \cdot h) \cdot (g \cdot f) \arrow[rr, "\alpha"] && l \cdot (h \cdot (g \cdot f))
\end{tikzcd}\end{center}
These axioms are to be read as ordinary commutative diagrams in the appropriate hom categories.
\end{defin}

We say that a 2-category is \textit{strict} if its unitors and associator are identities (so that horizontal composition is unital and associative \textit{on the nose}). It is also worth stressing that what we call 2-categories are sometimes called \textit{bicategories}, in which case the name \textit{2-categories} is reserved for the strict ones (or, in select cases, to some semi-strict notion).

\begin{rmk}\thmheader{Diagrammatic and string diagram notation}

There are two widespread ways of depicting cells and their compositions in a 2-category $\mathcal{B}$: \textit{diagrammatic notation} and \textit{string diagram notation}. We first give the example of a 2-cell $\colorformula{brown}{\phi} : \colorformula{green}{f} \to \colorformula{purple}{g} : \colorformula{orange}{X} \to \colorformula{cyan}{Y} : \mathcal{C}$: its depiction would then be (diagrammatic notation is on the left, string diagram notation on the right)

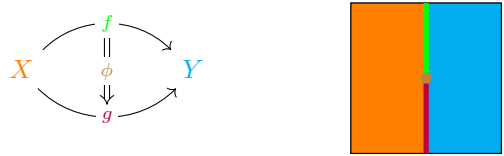
\begin{center}\begin{tikzcd}
\colorformula{orange}{X} \arrow[rr, bend left=45, "\colorformula{green}{f}"{name=f, description}] \arrow[rr, bend right=45, "\colorformula{purple}{g}"{name=g, description}] \arrow[Rightarrow, from=f, to=g, "\colorformula{brown}{\phi}" description] && \colorformula{cyan}{Y} &&
\end{tikzcd}
\begin{tikzpicture}[baseline=(current bounding box.center)] %[scale=0.5]
    \filldraw[draw=black, fill=white] (0, -1) rectangle (2, 1);
    \draw[fill=orange] (0, 1) -- (1, 1) -- (1, -1) -- (0, -1) -- cycle;
    \draw[fill=cyan] (1, 1) -- (2, 1) -- (2, -1) -- (1, -1) -- cycle;

    \node (c1) at (1, 0) [bullet, brown] {\colorformula{black}{\phi}};

    \path[draw, line width=2pt, line width=2pt, green] (1, 1) -- (c1);
    \path[draw, line width=2pt, line width=2pt, purple] (c1) -- (1, -1);
\end{tikzpicture}\end{center}

In diagrammatic notation 0-cells are points, 1-cells are arrows between the points, and 2-cells are (directed) areas enclosed by arrows. String diagram notation is the dual (in the sense of Poincar\`e duality) of diagrammatic notation: 0-cells are areas, 1-cells are strings, and 2-cells are points. Moreover, identity 1-cells, unitors, and associators are (almost) never depicted in string diagram notation\footnote{When they are, they are depicted with dotted strings. }. 
\end{rmk}

There are a lot of 2-categories which people are (more or less) used to dealing with:
\begin{ex}\thmheader{Examples of 2-categories}
\begin{enumerate}
\item Most notions of a ``structured category'' yield a 2-category: $\mathbf{Cat}$ is the 2-category of categories, functors and natural transformations; $\mathbf{MonCat}$ is the 2-category of monoidal categories, monoidal functors and monoidal natural transformations; $\mathbf{BrMonCat}$ for braided such; $\mathbf{SymMonCat}$ for symmetric such and $\mathbf{CartMonCat}$ for cartesian monoidal categories\footnote{Each of these comes in four distinct incarnation, using strict, weak, lax and colax monoidal functors as morphisms.}. We also have the 2-categories of (co-) complete categories, (co-) continuous functors and natural transformation, and so on. All these are examples of \textit{strict} 2-categories.
\item Every category $\mathcal{C}$ can be seen as a strict 2-category $\mathbf{LDisc}(\mathcal{C})$ (the \textit{locally discrete} 2-category generated by $\mathcal{C}$) with the same objects as $\mathcal{C}$, and whose hom-categories are the discrete categories on the hom-sets of $\mathcal{C}$; meaning $\mathbf{LDisc}(\mathcal{C})[X, Y] := \mathbf{Disc}(\mathcal{C}[X, Y])$. More concretely, this 2-category has the same objects and morphisms of the original category, and only identity 2-cells. From now on, we will abuse notation and omit any mention of $\mathbf{LDisc}$, treating every category as the corresponding locally discrete 2-category.
\item Given a monoidal category $\mathbb{C} = (\mathcal{C}, \otimes, I, \dots)$ with (an explicit choice for) all binary coequalizers, and such that the tensor product preserves (the choices of) such coequalizers in each argument separately, we can define a 2-category $\mathbf{BiMod}(\mathbb{C})$ whose 0-cells are monoid objects in $\mathcal{C}$, 1-cells are bimodules, 2-cells are bimodule homomorphisms, identities are the monoids themselves seen as bimodules, and composition of 1-cells is given by tensoring bimodules along the monoid in the middle. The associators can be defined using the universal property of coequalizers, while unitors witness the statement that tensoring a bimodule by a monoid (seen as a bimodule over itself) is an equivalence. Under appropriate assumptions, the same can be done using comonoids and bi-comodules (assuming the existence of all binary equalizers and that the tensor product preserves them).
\item Given a cosmos $\mathbb{V}$ (in the sense of B\'enabou \cite{elementary-cosmoi}), there is a strict 2-category $\mathbb{V}\mathbf{Cat}$ with objects $\mathbb{V}$-enriched categories, hom-categories the (ordinary, i.e. \textit{not} $\mathbf{V}$-enriched) categories of $\mathbf{V}$-enriched functors and $\mathbf{V}$-enriched natural transformations (with the obvious choices for identity 1-cells, unitors and associators). In the case $\mathbb{V} = \mathbf{Set}$ this recovers the familiar 2-category $\mathbf{Cat} = \mathbf{SetCat}$ of categories, functors and natural transformations.
\item Given a monoidal category $(\mathcal{C}, I, \otimes, \lambda, \rho, \alpha)$, we can construct a 2-category $\mathbb{B}\mathcal{C}$ with just one object $* : \mathbb{B}\mathcal{C}$, such that $\mathbb{B}\mathcal{C}[*, *] := \mathcal{C}$: 1-cells are objects in $\mathcal{C}$, 2-cells are morphisms in $\mathcal{C}$, the identity is $I : \mathcal{C}$ and horizontal composition is the tensor product $\otimes$. We call it the \textit{delooping} of the monoidal category $\mathcal{C}$, and it is strict precisely when the monoidal structure on $\mathcal{C}$ is. In particular, the terminal category $\mathbb{1}$ admits a strict monoidal structure. Its delooping, which we still denote with $\mathbb{1} := \mathbb{B}\mathbb{1}$, is the terminal (strict) 2-category.
\item Given a category $\mathcal{C}$ with binary pullbacks, we can define the 2-category $\mathbf{Span}(\mathcal{C})$ whose 0-cells are $\mathcal{C}$'s objects, whose 1-cells $X \to Y$ are spans $X \leftarrow Z \to Y$, 2-cells are morphisms between the vertices of the related spans, identities given by identity spans $X \xleftarrow{\mathbf{Id}} X \xrightarrow{\mathbf{Id}} X$ and horizontal composition given by pulling back the spans.
\item Given a 2-category $\mathcal{B}$, we can define 2-categories $\mathcal{C}^{co}, \mathcal{C}^{op}, \mathcal{C}^{coop}$ by (resp.) reversing 2-cells, 1-cells or both. These are strict precisely when $\mathcal{B}$ is.
\end{enumerate}
\end{ex}

\begin{rmk}\thmheader{Loop space of a pointed 2-category}

Given an object $X$ in a 2-category $\mathcal{B}$, the hom category $\mathcal{B}[X, X]$ is monoidal.
\end{rmk}

\begin{defin}\thmheader{Lax 2-functor}

Given 2-categories $\mathcal{B}, \mathcal{C}$, a lax 2-functor $\mathcal{F} : \mathcal{B} \to \mathcal{C}$ is given by
\begin{enumerate}
\item For any object $X : \mathcal{B}$, an object $\mathcal{F}(X) : \mathcal{C}$,
\item For any objects $X, Y : \mathcal{B}$, a functor $\mathcal{F} : \mathcal{B}[X, Y] \to \mathcal{C}[\mathcal{F}(X), \mathcal{F}(Y)]$,
\item For any object $X : \mathcal{B}$, a 2-cell
\[\mathcal{F}_u : \mathbf{Id}_{\mathcal{F}(X)} \to \mathcal{F}(\mathbf{Id}_{X})\]
\item For any $X \xrightarrow{f} Y \xrightarrow{g} Z : \mathcal{B}$, a 2-cell (natural in $f, g$)
\[\mathcal{F}_c: \mathcal{F}(g) \cdot \mathcal{F}(f) \to \mathcal{F}(g \cdot f)\]
\end{enumerate}
satisfying the following axioms:
\begin{enumerate}
\item Unitality: for $\colorformula{green}{f} : \colorformula{orange}{X} \to \colorformula{cyan}{Y} : \mathcal{B}$
\begin{center}\begin{tikzcd}
\mathbf{Id}_{\mathcal{F}(Y)} \cdot \mathcal{F}(\colorformula{green}{f}) \arrow[r, "\lambda"] \arrow[d, "\mathcal{F}_u \cdot \mathbf{Id}_{\mathcal{F}(\colorformula{green}{f})}"'] &
\mathcal{F}(\colorformula{green}{f}) \arrow[dd, equals] &
\mathcal{F}(\colorformula{green}{f}) \cdot \mathbf{Id}_{\mathcal{F}(X)} \arrow[l, "\rho"'] \arrow[d, "\mathbf{Id}_{\mathcal{F}(\colorformula{green}{f})} \cdot \mathcal{F}_u"] \\
\mathcal{F}(\mathbf{Id}_Y) \cdot \mathcal{F}(\colorformula{green}{f}) \arrow[d, "\mathcal{F}_c"'] &&
\mathcal{F}(\colorformula{green}{f}) \cdot \mathcal{F}(\mathbf{Id}_X) \arrow[d, "\mathcal{F}_c"] \\
\mathcal{F}(\mathbf{Id}_Y \cdot \colorformula{green}{f}) \arrow[r, "\mathcal{F}(\lambda)"'] &
\mathcal{F}(\colorformula{green}{f}) &
\mathcal{F}(\colorformula{green}{f} \cdot \mathbf{Id}_X) \arrow[l, "\mathcal{F}(\rho)"]
\end{tikzcd}\end{center}

Or, in string diagram notation (where the black line is $\mathcal{F}(\mathbf{Id})$ and we don't depict identity 1-cells or coherences),

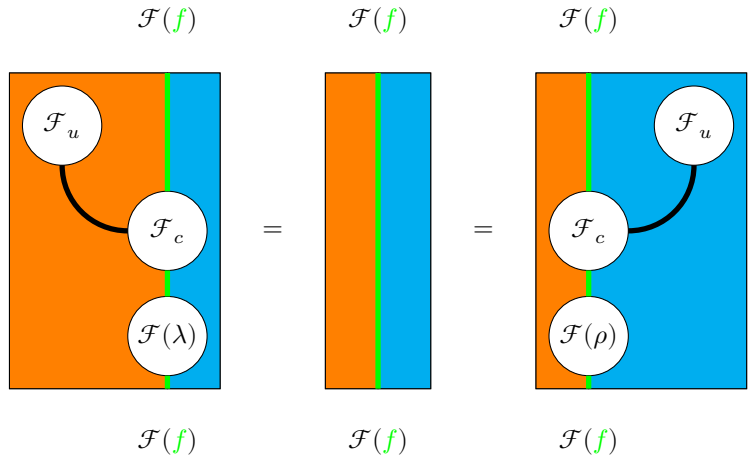
\begin{center}\begin{tikzpicture}[baseline=(current bounding box.center), scale=0.7]
    \filldraw[draw=black, fill=white] (-1, -3) rectangle (3, 3);
    \draw[fill=orange] (-1, -3) rectangle (2, 3);
    \draw[fill=cyan] (2, 3) rectangle (3, -3);
    
    \node (ul) at (0, 2) [circle, draw=black, fill=white, minimum size=30pt, inner sep=0pt, outer sep=0pt] {$\mathcal{F}_u$};
    \node (dr) at (2, 0) [circle, draw=black, fill=white, minimum size=30pt, inner sep=0pt, outer sep=0pt] {$\mathcal{F}_c$};
    \node (ddr) at (2, -2) [circle, draw=black, fill=white, minimum size=30pt, inner sep=0pt, outer sep=0pt] {$\mathcal{F}(\lambda)$};

    \node at (2, 4) [] {$\mathcal{F}(\colorformula{green}{f})$};
    \node at (2, -4) [] {$\mathcal{F}(\colorformula{green}{f})$};

    \path[draw, line width=2pt, line width=2pt, green] (2, 3) -- (dr);
    \path[draw, line width=2pt, line width=2pt, black] (ul) to [out=-90, in=180] (dr);
    \path[draw, line width=2pt, line width=2pt, green] (dr) -- (ddr);
    \path[draw, line width=2pt, line width=2pt, green] (ddr) -- (2, -3);

    \node at (4, 0) [] {=};

    \filldraw[draw=black, fill=white] (5, -3) rectangle (7, 3);
    \draw[fill=orange] (5, -3) rectangle (6, 3);
    \draw[fill=cyan] (6, -3) rectangle (7, 3);

    \node at (6, 4) [] {$\mathcal{F}(\colorformula{green}{f})$};
    \node at (6, -4) [] {$\mathcal{F}(\colorformula{green}{f})$};

    \path[draw, line width=2pt, line width=2pt, green] (6, 3) -- (6, -3);

    \node at (8, 0) [] {=};

    \filldraw[draw=black, fill=white] (9, -3) rectangle (13, 3);
    \draw[fill=orange] (9, -3) rectangle (10, 3);
    \draw[fill=cyan] (10, -3) rectangle (13, 3);

    \node (ur) at (12, 2) [circle, draw=black, fill=white, minimum size=30pt, inner sep=0pt, outer sep=0pt] {$\mathcal{F}_u$};
    \node (dl) at (10, 0) [circle, draw=black, fill=white, minimum size=30pt, inner sep=0pt, outer sep=0pt] {$\mathcal{F}_c$};
    \node (ddl) at (10, -2) [circle, draw=black, fill=white, minimum size=30pt, inner sep=0pt, outer sep=0pt] {$\mathcal{F}(\rho)$};

    \node at (10, 4) [] {$\mathcal{F}(\colorformula{green}{f})$};
    \node at (10, -4) [] {$\mathcal{F}(\colorformula{green}{f})$};

    \path[draw, line width=2pt, line width=2pt, green] (10, 3) -- (dl);
    \path[draw, line width=2pt, line width=2pt, green] (dl) -- (ddl);
    \path[draw, line width=2pt, line width=2pt, green] (ddl) -- (10, -3);
    \path[draw, line width=2pt, line width=2pt, black] (ur) to [out=-90, in=0] (dl);
\end{tikzpicture}\end{center}

\item Associativity: for $\colorformula{orange}{W} \xrightarrow{\colorformula{green}{f}} \colorformula{cyan}{X} \xrightarrow{\colorformula{purple}{g}} \colorformula{yellow}{Y} \xrightarrow{\colorformula{brown}{h}} \colorformula{lime}{Z} : \mathcal{B}$
\begin{center}\begin{tikzcd}
(\mathcal{F}(h) \cdot \mathcal{F}(g)) \cdot \mathcal{F}(f) \arrow[rr, "\alpha"] \arrow[d, "\mathcal{F}_c \cdot \mathbf{Id}_{\mathcal{F}(f)}"'] &&
\mathcal{F}(h) \cdot (\mathcal{F}(g) \cdot \mathcal{F}(f)) \arrow[d, "\mathbf{Id}_{\mathcal{F}(h)} \cdot \mathcal{F}_c"] \\
\mathcal{F}(h \cdot g) \cdot \mathcal{F}(f) \arrow[d, "\mathcal{F}_c"'] &&
\mathcal{F}(h) \cdot \mathcal{F}(g \cdot f) \arrow[d, "\mathcal{F}_c"] \\
\mathcal{F}((h \cdot g) \cdot f) \arrow[rr, "\mathcal{F}(\alpha)"'] &&
\mathcal{F}(h \cdot (g \cdot f))
\end{tikzcd}\end{center}

Which corresponds to the following string diagram equation

\begin{center}\begin{tikzpicture}

\end{tikzpicture}\end{center}
\end{enumerate}
We say that
\begin{itemize}
\item a \textit{strict 2-functor} $\mathcal{F} : \mathcal{B} \to \mathcal{C}$ is a lax 2-functor where $\mathcal{F}_u$ and $\mathcal{F}_c$ are all identities;
\item a \textit{(pseudo) 2-functor} $\mathcal{F} : \mathcal{B} \to \mathcal{C}$ is a lax 2-functor where $\mathcal{F}_u$ and $\mathcal{F}_c$ are all isomorphisms (we will omit the ``pseudo'' part of the name most of the times);
\item a \textit{colax 2-functor} $\mathcal{F} : \mathcal{B} \to \mathcal{C}$ is a lax 2-functor $\mathcal{F} : \mathcal{B} \to \mathcal{C}^{co}$ (i.e. if $\mathcal{F}_u$ and $\mathcal{F}_c$ have their sources and target flipped).
\end{itemize}
\end{defin}

\begin{ex}\thmheader{Examples of 2-functors}
\begin{enumerate}
\item Every ordinary functor $\mathcal{F} : \mathcal{C} \to \mathcal{D}$ can be seen as a strict 2-functor $\mathbf{LDisc}(\mathcal{F}) : \mathcal{C} \to \mathcal{D}$ in the obvious way: we map objects to objects as $\mathcal{F}$ did, and we do the same for 1-cells. There is no choice for 2-cells, and 2-functoriality follows immediately from $\mathcal{F}$'s functoriality. We will abuse notation and write $\mathcal{F}$ for $\mathbf{LDisc}(\mathcal{F})$.
\item Monoidal functors preserving binary coequalizers\footnote{It seems that if $\mathcal{F}$ is lax monoidal, there is no need for it to preserve coequalizer at all in order for it to induce a lax 2-functor.} $\mathbb{F} : \mathbb{C} \to \mathbb{D}$ between monoidal categories as in point 3 in \ref{Examples of 2-categories} induce 2-functors
\[\mathbf{BiMod}(\mathbb{F}) : \mathbf{BiMod}(\mathbb{C}) \to \mathbf{BiMod}(\mathbb{D})\]
which are as strict as $\mathbb{F}$ is.
\item Any lax monoidal functor $\mathcal{F} : \mathbb{V} \to \mathbb{W}$ between cosmoi induces a strict ``change-of-base'' 2-functor $\mathcal{F}_*$ (for the details, see prop. 6.4.3 in \cite{handbook-categorical-algebra}, vol. 2).
\item Every monoidal functor $\mathcal{F} : \mathcal{C} \to \mathcal{D}$ can be delooped to a 2-functor $\mathbb{B}\mathcal{F} : \mathbb{B}\mathcal{C} \to \mathbb{B}\mathcal{D}$ in the obvious way; it will be as strict as $\mathcal{F}$ is.
\end{enumerate}
\end{ex}

\begin{defin}\thmheader{Lax 2-natural transformation}

Given (strict, ordinary, lax or colax) 2-functors $\mathcal{F}, \mathcal{G} : \mathcal{B} \to \mathcal{C}$, a lax 2-natural transformation $\Phi : \mathcal{F} \to \mathcal{G}$ is given by
\begin{enumerate}
\item For any $X : \mathcal{B}$, a 1-cell $\Phi_X : \mathcal{F}(X) \to \mathcal{G}(X)$,
\item For any $f : X \to Y : \mathcal{B}$ an 2-cell $\Phi_f : \mathcal{G}(f) \cdot \Phi_X \to \Phi_Y \cdot \mathcal{\mathcal{F}(f)}$ naturally in $f$, meaning that for every $\phi : f \to g : X \to Y : \mathcal{B}$ the following diagram is to be commutative in $\mathcal{C}[\mathcal{F}(X), \mathcal{G}(Y)]$
\begin{center}\begin{tikzcd}
\mathcal{G}(f) \cdot \Phi_X \arrow[r, "\mathcal{G}(\phi) \cdot \Phi_X"] \arrow[d, "\Phi_f"'] &
\mathcal{G}(g) \cdot \Phi_X \arrow[d, "\Phi_g"] \\
\Phi_Y \cdot \mathcal{F}(f) \arrow[r, "\Phi_Y \cdot \mathcal{F}(\phi)"] &
\Phi_Y \cdot \mathcal{F}(g)
\end{tikzcd}\end{center}
\end{enumerate}
satisfying the following axioms:
\begin{enumerate}
\item Compatibility with unitors:
\begin{center}\begin{tikzcd}
\mathbf{Id}_{\mathcal{G}(X)} \cdot \Phi_X \arrow[d, "\lambda"] \arrow[r, "\mathcal{G}_u \cdot \mathbf{Id}_{\Phi_X}"] & \mathcal{G}(\mathbf{Id}_X) \cdot \Phi_X \arrow[dd, "\Phi_{\mathbf{Id}_X}"] \\
\Phi_X \arrow[d, "\rho^{-1}"] & \\
\Phi_X \cdot \mathbf{Id}_{\mathcal{F}(X)} \arrow[r, "\mathbf{Id}_{\Phi_X} \cdot \mathcal{F}_u"] & \Phi_X \cdot \mathcal{F}(\mathbf{Id}_X)
\end{tikzcd}\end{center}
\item Compatibility with associators:
\begin{center}\begin{tikzcd}
(\mathcal{G}(g) \cdot \mathcal{G}(f)) \cdot \Phi_X \arrow[d, "\alpha"] \arrow[r, "\mathcal{G}_c \cdot \mathbf{Id}_{\Phi_X}"] & \mathcal{G}(g \cdot f) \cdot \Phi_X \arrow[ddddd, "\Phi_{g \cdot f}"] \\
\mathcal{G}(g) \cdot (\mathcal{G}(f) \cdot \Phi_X) \arrow[d, "\mathbf{Id}_{\mathcal{G}(g)} \cdot \Phi_f"] & \\
\mathcal{G}(g) \cdot (\Phi_Y \cdot \mathcal{F}(f)) \arrow[d, "\alpha^{-1}"] & \\
(\mathcal{G}(g) \cdot \Phi_Y) \cdot \mathcal{F}(f) \arrow[d, "\Phi_g \cdot \mathbf{Id}_{\mathcal{F}(f)}"] & \\
(\Phi_Z \cdot \mathcal{F}(g)) \cdot \mathcal{F}(f) \arrow[d, "\alpha"] & \\
\Phi_Z \cdot (\mathcal{F}(g) \cdot \mathcal{F}(f)) \arrow[r, "\mathbf{Id}_{\Phi_Z} \cdot \mathcal{F}_c"] & \Phi_Z \cdot \mathcal{F}(g \cdot f)
\end{tikzcd}\end{center}
\end{enumerate}
We say that
\begin{itemize}
\item a \textit{strict 2-natural transformation} $\Phi : \mathcal{F} \to \mathcal{G}$ is a lax 2-natural transformation where all the $\Phi_f$ are identities;
\item a \textit{2-natural transformation} $\Phi : \mathcal{F} \to \mathcal{G}$ is a lax 2-natural transformation where all the $\Phi_f$ are isomorphisms;
\item a \textit{colax 2-natural transformation} $\Phi : \mathcal{F} \to \mathcal{G}$ is a lax 2-natural transformation $\Phi : \mathcal{F} \to \mathcal{G} : \mathcal{B} \to \mathcal{C}^{co}$ (i.e. if $\Phi_f$ has its source and target flipped).
\end{itemize}
\end{defin}

\begin{ex}\thmheader{Examples of 2-natural transformations}
\begin{itemize}
\item Every ordinary natural transformation $\phi : \mathcal{F} \to \mathcal{G} : \mathcal{C} \to \mathcal{D}$ can be seen as a strict 2-natural transformation $\mathbf{LDisc}(\mathcal{\phi}) : \mathbf{LDisc}(\mathcal{F}) \to \mathbf{LDisc}(\mathcal{G})$ in the obvious way.
\item Given monoidal categories $\mathbb{C}, \mathbb{D}$ with all binary coequalizers and functors $\mathbb{F}, \mathbb{G} : \mathbb{C} \to \mathbb{D}$ preserving them, a monoidal natural transformation $\phi : \mathbb{F} \to \mathbb{G}$ extends to a strict 2-transformation
\[\mathbf{BiMod}(\phi) : \mathbf{BiMod}(\mathbb{F}) \to \mathbf{BiMod}(\mathbb{G})\]
\item Any monoidal natural transformation $\phi : \mathcal{F} \to \mathcal{G} : \mathcal{C} \to \mathcal{D}$ between (lax) monoidal functors induces a strict 2-natural transformation between the induced change-of-base strict 2-functors (see 4.3.1 in \cite{change-of-base} for the details).
\item Any monoidal natural transformation between (strict, weak, co/lax) monoidal functors gives rise to a 2-natural transformation between the corresponding deloopings.
\end{itemize}
\end{ex}

\begin{defin}\thmheader{Modification}

Given 2-natural transformations $\Phi, \Psi : \mathcal{F} \to \mathcal{G} : \mathcal{B} \to \mathcal{C}$, a modification
\[\Gamma : \Phi \to \Psi : \mathcal{F} \to \mathcal{G} : \mathcal{B} \to \mathcal{C}\]
consists of the datum of a 2-cell $\Gamma_X$ for any $X : \mathcal{B}$, such that for any $f : X \to Y : \mathcal{B}$, the following diagram commutes
\begin{center}\begin{tikzcd}
\mathcal{G}(f) \cdot \Phi_X \arrow[rr, "\Phi_f"] \arrow[dd, "\mathbf{Id}_{\mathcal{G}(f)} \cdot \Gamma_X"'] &&
\Phi_Y \cdot \mathcal{F}(f) \arrow[dd, "\Gamma_Y \cdot \mathbf{Id}_{\mathcal{F}(f)}"] \\
\\
\mathcal{G}(f) \cdot \Psi_X \arrow[rr, "\Psi_f"'] &&
\Psi_Y \cdot \mathcal{F}(f)
\end{tikzcd}\end{center}
\end{defin}

As one might have guessed, all these gadgets (2-functors, 2-natural transformations, 2-modifications) make up 2-categories and (sometimes) compose appropriately. We will not get into the details, just stating the following result without any proofs.

\begin{prop}\thmheader{2-functor 2-categories}

Given 2-categories $\mathcal{B}, \mathcal{C}$ there are 2-categories
\begin{enumerate}
\item 
Given 
$\mathbf{St}_{\mathbf{St}}[\mathcal{B}, \mathcal{C}]$ made of strict 2-functors, strict 2-natural transformations and modifications,
\item $\mathbf{St}_{\mathbf{Ps}}[\mathcal{B}, \mathcal{C}]$ made of strict 2-functors, 2-natural transformations and modifications,
\item $\mathbf{St}_{\mathbf{Lax}}[\mathcal{B}, \mathcal{C}]$ made of strict 2-functors, lax 2-natural transformations and modifications,
\item $\mathbf{St}_{\mathbf{CoLax}}[\mathcal{B}, \mathcal{C}]$ made of strict 2-functors, colax 2-natural transformations and modifications,
\item $\mathbf{Ps}_{\mathbf{Ps}}[\mathcal{B}, \mathcal{C}]$ made of 2-functors, 2-natural transformations and modifications,
\item $\mathbf{Ps}_{\mathbf{Lax}}[\mathcal{B}, \mathcal{C}]$ made of 2-functors, lax 2-natural transformations and modifications,
\item $\mathbf{Ps}_{\mathbf{CoLax}}[\mathcal{B}, \mathcal{C}]$ made of 2-functors, colax 2-natural transformations and modifications,
\item $\mathbf{Lax}[\mathcal{B}, \mathcal{C}]$ made of lax 2-functors, lax 2-natural transformations and modifications,
\item $\mathbf{CoLax}[\mathcal{B}, \mathcal{C}]$ made of colax 2-functors, colax 2-natural transformations and modifications.
\end{enumerate}
We also get composition 2-functors, but only in cases (1, 2, 5). Moreover, the above list is \textit{not} complete: we can consider, e.g., lax 2-functors and 2-natural transformations, or colax 2-functors and lax 2-natural transformations, etc. Finally, a note about notation: the use of $\mathbf{Ps}$ is motivated by the widespread terminology calling pseudofunctors what we just call 2-functors.

In what follows we'll mostly be interested in (1, 3, 8).
\end{prop}

\begin{defin}\thmheader{Monad in a 2-category}

Given a 2-category $\mathcal{B}$, a monad in it is given by the data of a 0-cell $X$ and a monoid object in $\mathcal{B}[X, X]$. More explicitly, this means a 4-tuple
\[(X, t : X \to X, \eta : \mathbf{Id}_X \to t, \mu : t \cdot t \to t)\]
such that the following diagrams commute in $\mathcal{B}[X, X]$, encoding unitality and associativity (resp.)
\begin{center}\begin{tikzcd}
\mathbf{Id}_X \cdot t \arrow[r, "\lambda"] \arrow[d, "\eta \cdot t"'] &
t \arrow[d, equal] &
t \cdot \mathbf{Id}_X \arrow[l, "\rho"'] \arrow[d, "t \cdot \eta"] &
(t \cdot t) \cdot t \arrow[rr, "\alpha"] \arrow[d, "\mu \cdot t"'] &&
t \cdot (t \cdot t) \arrow[d, "t \cdot \mu"] \\
t \cdot t \arrow[r, "\mu"'] &
t &
t \cdot t \arrow[l, "\mu"] &
t \cdot t \arrow[r, "\mu"'] &
t &
t \cdot t \arrow[l, "\mu"]
\end{tikzcd}\end{center}
\end{defin}

Analogously, we can define a comonad in a 2-category $\mathcal{B}$ as a monad in $\mathcal{B}^{co}$. It's easy to see that 2-functors interact nicely with monads and comonads

\begin{prop}\thmheader{Lax 2-functors preserve monads}

Let $\mathcal{F} : \mathcal{B} \to \mathcal{C}$ be a lax 2-functor, and $(X, t, \eta, \mu)$ be a monad in $\mathcal{B}$. Then we get a monad in $\mathcal{C}$
\[\mathcal{F}(X, t, \eta, \mu) = (\mathcal{F}(X), \mathcal{F}(t), \mathcal{F}(\eta) \circ \mathcal{F}_u, \mathcal{F}(\mu) \circ \mathcal{F}_c)\]

\begin{proof}
Start with unitality:
\begin{center}\begin{tikzcd}[scale=0.8]
\mathbf{Id}_{\mathcal{F}(X)} \cdot \mathcal{F}(t) \arrow[rr, "\lambda"] \arrow[d, "\mathcal{F}_u \cdot \mathbf{Id}_{\mathcal{F}(t)}"'] &&
\mathcal{F}(t) \arrow[d, equals] &&
\mathcal{F}(t) \cdot \mathbf{Id}_{\mathcal{F}(X)} \arrow[ll, "\rho"'] \arrow[d, "\mathbf{Id}_{\mathcal{F}(t)} \cdot \mathcal{F}_u"] \\
\mathcal{F}(\mathbf{Id}_X) \cdot \mathcal{F}(t) \arrow[d, "\mathcal{F}(\eta) \cdot \mathbf{Id}_{\mathcal{F}(t)}"'] \arrow[r, "\mathcal{F}_c"] &
\mathcal{F}(\mathbf{Id}_X \cdot t) \arrow[r, "\mathcal{F}(\lambda)"] \arrow[d, "\mathcal{F}(\eta \cdot \mathbf{Id}_t)"'] &
\mathcal{F}(t) \arrow[d, equals] &
\mathcal{F}(t \cdot \mathbf{Id}_X) \arrow[d, "\mathcal{F}(\mathbf{Id}_t \cdot \eta)"] \arrow[l, "\mathcal{F}(\rho)"'] &
\mathcal{F}(t) \cdot \mathcal{F}(\mathbf{Id}_X) \arrow[d, "\mathbf{Id}_{\mathcal{F}(t)} \cdot \mathcal{F}(\eta)"] \arrow[l, "\mathcal{F}_c"'] \\
\mathcal{F}(t) \cdot \mathcal{F}(t) \arrow[r, "\mathcal{F}_c"'] &
\mathcal{F}(t \cdot t) \arrow[r, "\mathcal{F}(\mu)"'] &
\mathcal{F}(t) &
\mathcal{F}(t \cdot t) \arrow[l, "\mathcal{F}(\mu)"] &
\mathcal{F}(t) \cdot \mathcal{F}(t) \arrow[l, "\mathcal{F}_c"]
\end{tikzcd}\end{center}

The pentagons on the top commute by the unitality assumption on the 2-functor $\mathcal{F}$, the squares at the extreme left and right by naturality of $\mathcal{F}_c$, while the remaining two squares are the unitality axioms on the monad $(X, t, \eta, \mu)$.
\begin{center}\begin{tikzcd}
(\mathcal{F}(t) \cdot \mathcal{F}(t)) \cdot \mathcal{F}(t) \arrow[rrrr, "\alpha"] \arrow[d, "\mathcal{F}_c \cdot \mathbf{Id}_{\mathcal{F}(t)}"'] &&&&
\mathcal{F}(t) \cdot (\mathcal{F}(t) \cdot \mathcal{F}(t)) \arrow[d, "\mathbf{Id}_{\mathcal{F}(t)} \cdot \mathcal{F}_c"] \\
\mathcal{F}(t \cdot t) \cdot \mathcal{F}(t) \arrow[d, "\mathcal{F}(\mu) \cdot \mathbf{Id}_{\mathcal{F}(t)}"'] \arrow[r, "\mathcal{F}_c"] &
\mathcal{F}((t \cdot t) \cdot t) \arrow[rr, "\mathcal{F}(\alpha)"] \arrow[d, "\mathcal{F}(\mu \cdot \mathbf{Id}_t)"'] &&
\mathcal{F}(t \cdot (t \cdot t)) \arrow[d, "\mathcal{F}(\mathbf{Id}_t \cdot \mu)"] &
\mathcal{F}(t) \cdot \mathcal{F}(t \cdot t) \arrow[d, "\mathbf{Id}_{\mathcal{F}(t)} \cdot \mathcal{F}(\mu)"] \arrow[l, "\mathcal{F}_c"] \\
\mathcal{F}(t) \cdot \mathcal{F}(t) \arrow[r, "\mathcal{F}_c"'] &
\mathcal{F}(t \cdot t) \arrow[r, "\mathcal{F}(\mu)"'] &
\mathcal{F}(t) &
\mathcal{F}(t \cdot t) \arrow[l, "\mathcal{F}(\mu)"] &
\mathcal{F}(t) \cdot \mathcal{F}(t) \arrow[l, "\mathcal{F}_c"]
\end{tikzcd}\end{center}

The hexagon on top is the associativity assumption on the 2-functor $\mathcal{F}$, the squares on the left and right commute by naturality of $\mathcal{F}_c$ and the remaining pentagon is the associativity axiom on the monad.
\end{proof}
\end{prop}

An analogous result holds for comonads and colax 2-functors

\begin{prop}\thmheader{Colax 2-functors preserve comonads}
Let $\mathcal{F} : \mathcal{B} \to \mathcal{C}$ a colax 2-functor, and $(X, t, \epsilon, \Delta)$ a comonad in $\mathcal{B}$. Then we get a comonad in $\mathcal{C}$
\[\mathcal{F}(X, t, \epsilon, \Delta) := (\mathcal{F}(X), \mathcal{F}(t), \mathcal{F}_u \circ \mathcal{F}(\epsilon), \mathcal{F}_c \circ \mathcal{F}(\Delta))\]
\begin{proof} Identical to that in \ref{Lax 2-functors preserve monads}.
\end{proof}
\end{prop}

More is true: every monad can be encoded as a lax 2-functor

\begin{prop}\thmheader{Monads as lax 2-functors}

Given a monad $(X, t, \eta, \mu)$ in a 2-category $\mathcal{B}$, there is a lax 2-functor that sends the unique monad $(*, \mathbf{Id}_*, \mathbf{Id}_{\mathbf{Id}_*}, \mathbf{Id}_{\mathbf{Id}_*})$ in $\mathbb{1}$ to $(X, t, \eta, \mu)$.
\begin{proof}

Clearly the lax 2-functor $\mathcal{F} : \mathbb{1} \to \mathcal{B}$ will send $*$ to $\mathcal{F}(*) = X$ and $\mathbf{Id}_*$ to $\mathcal{F}(\mathbf{Id}_*) = t$. In order for it to be a lax functor, we need 1-cells
\[\mathcal{F}_u : \mathbf{Id}_{\mathcal{F}(*)} := \mathbf{Id}_X \xrightarrow{\eta} t =: \mathcal{F}(\mathbf{Id}_*)\]
\[\mathcal{F}_c : \mathcal{F}(\mathbf{Id}_*) \cdot \mathcal{F}(\mathbf{Id}_*) := t \cdot t \xrightarrow{\mu} t =: \mathcal{F}(\mathbf{Id}_*) = \mathcal{F}(\mathbf{Id}_* \cdot \mathbf{Id}_*)\]
where the last equality follows from the fact that $\mathbb{1}$ is a strict 2-category. Since the only 1-cell in $\mathbb{1}$ is $\mathbf{Id}_*$, the lax 2-functor axioms become
\begin{center}\begin{tikzcd}
\mathbf{Id}_X \cdot t \arrow[r, "\lambda"] \arrow[d, "\eta \cdot \mathbf{Id}_t"'] &
t \arrow[dd, equals] &
t \cdot \mathbf{Id}_X \arrow[l, "\rho^{-1}"'] \arrow[d, "\mathbf{Id}_t \cdot \eta"] &
(t \cdot t) \cdot t \arrow[r, "\alpha"] \arrow[d, "\mu \cdot \mathbf{Id}_t"'] &
t \cdot (t \cdot t) \arrow[d, "\mathbf{Id}_t \cdot \mu"] \\
t \cdot t \arrow[d, "\mu"'] &&
t \cdot t \arrow[d, "\mu"] &
t \cdot t \arrow[d, "\mu"'] &
t \cdot t \arrow[d, "\mu"] \\
t \arrow[r, equals] &
t &
t \arrow[l, equals] &
t \arrow[r, equals] &
t
\end{tikzcd}\end{center}
which are exactly the monad axioms.
\end{proof}
\end{prop}

In total analogy, we can prove the following

\begin{prop}\thmheader{Comonads as colax 2-functors}

Given a comonad $(X, t, \epsilon, \Delta)$ in a 2-category $\mathcal{B}$ there's a colax 2-functor that sends the unique comonad $(*, \mathbf{Id}_*, \mathbf{Id}_{\mathbf{Id}_*}, \mathbf{Id}_{\mathbf{Id}_*})$ in $\mathbb{1}$ to $(X, s, \epsilon, \Delta)$.
\begin{proof}

Identical to that in \ref{Monads as lax 2-functors}.
\end{proof}
\end{prop}

These results suggest the following

\begin{defin}\thmheader{2-categories of monads and comonads}

Let $\mathcal{B}$ be a 2-category. We define the 2-categories
\[\mathbf{Mnd}(\mathcal{B}) := \mathbf{Lax}[\mathbb{1}, \mathcal{B}]\]
\[\mathbf{CoMnd}(\mathcal{B}) := \mathbf{CoLax}[\mathbb{1}, \mathcal{B}]\]
\end{defin}

This is worth unpacking, and we'll take the chance to introduce a notation we'll use throughout the rest of this section: string diagrams for 2-categories. In string diagram notation 0-cells are areas, 1-cells are strings and 2-cells are nodes; this comes with two big advantages: there's no way to depict identities, unitors and associators (so there will be no coherence issues), and it'll allow us to reason graphically about complex commutative diagrams.
\begin{enumerate}
\item A monad $\colorformula{blue}{\mathcal{T}} = (\colorformula{orange}{X}, \colorformula{blue}{t}, \colorformula{blue}{\eta}, \colorformula{blue}{\mu})$ would be depicted as
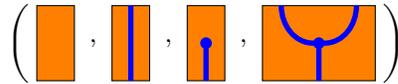
\begin{center}\begin{equation*}\begin{tikzpicture}[scale=0.5]
    \node at (-1, 0) [] {\Bigg{(}};

    \filldraw[draw=black, fill=orange] (-0.5, -1) rectangle (0.5, 1);
    
    \node at (1, 0) [] {,};

    \filldraw[draw=black, fill=orange] (1.5, -1) rectangle (2.5, 1);
    
    \node (u1) at (2, 1) [] {};
    \node (d1) at (2, -1) [] {};
    \path[draw, line width=2pt, line width=2pt, blue] (u1.center) -- (d1.center);

    \node at (3, 0) [] {,};

    \filldraw[draw=black, fill=orange] (3.5, -1) rectangle (4.5, 1);

    \node (c2) at (4, 0) [bullet, blue] {};
    \node (d2) at (4, -1) [] {};
    \path[draw, line width=2pt, line width=2pt, blue] (c2) -- (d2.center);

    \node at (5, 0) [] {,};

    \filldraw[draw=black, fill=orange] (5.5, -1) rectangle (8.5, 1);

    \node (ul3) at (6, 1) [] {};
    \node (ur3) at (8, 1) [] {};
    \node (c3) at (7, 0) [bullet, blue] {};
    \node (d3) at (7, -1) [] {};
    \path[draw, line width=2pt, line width=2pt, blue] (ul3.center) to [out=-90, in=180] (c3);
    \path[draw, line width=2pt, line width=2pt, blue] (ur3.center) to [out=-90, in=0] (c3);
    \path[draw, line width=2pt, line width=2pt, blue] (c3) -- (d3.center);

    \node at (9, 0) [] {\Bigg{)}};
\end{tikzpicture}\end{equation*}\end{center}
and the axioms would translate to
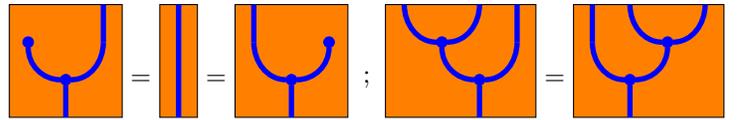
\begin{center}\begin{tikzpicture}[scale=0.5]
        \filldraw[draw=black, fill=orange] (-1.5, -1) rectangle (1.5, 2);
        \filldraw[draw=black, fill=orange] (2.5, -1) rectangle (3.5, 2);
        \filldraw[draw=black, fill=orange] (4.5, -1) rectangle (7.5, 2);
        \filldraw[draw=black, fill=orange] (8.5, -1) rectangle (12.5, 2);
        \filldraw[draw=black, fill=orange] (13.5, -1) rectangle (17.5, 2);

        \node (1-uur) at (1, 2) [] {};
        \node (1-ul) at (-1, 1) [bullet, blue] {};
        \node (1-ur) at (1, 1) [] {};
        \node (1-c) at (0, 0) [bullet, blue] {};
        \node (1-d) at (0, -1) [] {};

        \path[draw, line width=2pt, line width=2pt, blue] (1-uur.center) to (1-ur.center);
        \path[draw, line width=2pt, line width=2pt, blue] (1-ul) to [out=-90, in=180] (1-c);
        \path[draw, line width=2pt, line width=2pt, blue] (1-ur.center) to [out=-90, in=0] (1-c);
        \path[draw, line width=2pt, line width=2pt, blue] (1-c) to (1-d.center);

        \node at (2, 0) [] {=};

        \node (2-uu) at (3, 2) [] {};
        \node (2-d) at (3, -1) [] {};

        \path[draw, line width=2pt, line width=2pt, blue] (2-uu.center) to (2-d.center);

        \node at (4, 0) [] {=};

        \node (3-uul) at (5, 2) [] {};
        \node (3-ul) at (5, 1) [] {};
        \node (3-ur) at (7, 1) [bullet, blue] {};
        \node (3-c) at (6, 0) [bullet, blue] {};
        \node (3-d) at (6, -1) [] {};

        \path[draw, line width=2pt, line width=2pt, blue] (3-uul.center) to (3-ul.center);
        \path[draw, line width=2pt, line width=2pt, blue] (3-ul.center) to [out=-90, in=180] (3-c);
        \path[draw, line width=2pt, line width=2pt, blue] (3-ur) to [out=-90, in=0] (3-c);
        \path[draw, line width=2pt, line width=2pt, blue] (3-c) to (3-d.center);
        
        \node at (8, 0) [] {;};

        \node (4-uul) at (9, 2) [] {};
        \node (4-uuc) at (11, 2) [] {};
        \node (4-uur) at (12, 2) [] {};
        \node (4-uc) at (10, 1) [bullet, blue] {};
        \node (4-ur) at (12, 1) [] {};
        \node (4-c) at (11, 0) [bullet, blue] {};
        \node (4-d) at (11, -1) [] {};

        \path[draw, line width=2pt, line width=2pt, blue] (4-uul.center) to [out=-90, in=180] (4-uc);
        \path[draw, line width=2pt, line width=2pt, blue] (4-uuc.center) to [out=-90, in=0] (4-uc);
        \path[draw, line width=2pt, line width=2pt, blue] (4-uur.center) to (4-ur.center);
        \path[draw, line width=2pt, line width=2pt, blue] (4-uc) to [out=-90, in=180] (4-c);
        \path[draw, line width=2pt, line width=2pt, blue] (4-ur.center) to [out=-90, in=0] (4-c);
        \path[draw, line width=2pt, line width=2pt, blue] (4-c) to (4-d.center);

        \node at (13, 0) [] {=};

        \node (5-uul) at (14, 2) [] {};
        \node (5-uuc) at (15, 2) [] {};
        \node (5-uur) at (17, 2) [] {};
        \node (5-ul) at (14, 1) [] {};
        \node (5-uc) at (16, 1) [bullet, blue] {};
        \node (5-c) at (15, 0) [bullet, blue] {};
        \node (5-d) at (15, -1) [] {};

        \path[draw, line width=2pt, line width=2pt, blue] (5-uul.center) to (5-ul.center);
        \path[draw, line width=2pt, line width=2pt, blue] (5-uuc.center) to [out=-90, in=180] (5-uc);
        \path[draw, line width=2pt, line width=2pt, blue] (5-uur.center) to [out=-90, in=0] (5-uc);
        \path[draw, line width=2pt, line width=2pt, blue] (5-ul.center) to [out=-90, in=180] (5-c);
        \path[draw, line width=2pt, line width=2pt, blue] (5-uc) to [out=-90, in=0] (5-c);
        \path[draw, line width=2pt, line width=2pt, blue] (5-c) to (5-d.center);

        \node at (18, 0) [] {;};
\end{tikzpicture}\end{center}

\item A 1-cell $\colorformula{green}{\Phi} = (\colorformula{green}{f}, \colorformula{green}{\phi}) : (\colorformula{orange}{X}, \colorformula{blue}{t^X}, \colorformula{blue}{\eta^X}, \colorformula{blue}{\mu^X}) \to (\colorformula{cyan}{Y}, \colorformula{red}{t^Y}, \colorformula{red}{\eta^Y}, \colorformula{red}{\mu^Y})$ in $\mathbf{Mnd}(\mathcal{B})$ is given by a 1-cell $\colorformula{green}{f} : \colorformula{orange}{X} \to \colorformula{cyan}{Y}$ and a 2-cell
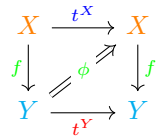
\begin{center}\begin{tikzcd}
\colorformula{orange}{X} \arrow[r, "\colorformula{blue}{t^X}"] \arrow[d, "\colorformula{green}{f}"'] & \colorformula{orange}{X} \arrow[d, "\colorformula{green}{f}"] \\
\colorformula{cyan}{Y} \arrow[r, "\colorformula{red}{t^Y}"'] \arrow[ur, Rightarrow, "\colorformula{green}{\phi}" description] & \colorformula{cyan}{Y}
\end{tikzcd}\end{center}
depicted as
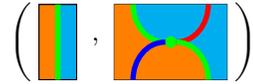
\begin{center}\begin{tikzpicture}[scale=0.5]
    \node at (0, 0) [] {\Bigg{(}};

    \filldraw[draw=black, fill=white] (0.5, -1) rectangle (1.5, 1);
    \draw[fill=orange] (0.5, 1) -- (1, 1) -- (1, -1) -- (0.5, -1) -- cycle;
    \draw[fill=cyan] (1, 1) -- (1.5, 1) -- (1.5, -1) -- (1, -1) -- cycle;

    \node (u1) at (1, 1) [] {};
    \node (d1) at (1, -1) [] {};

    \path[draw, line width=2pt, line width=2pt, green] (u1.center) -- (d1.center);

    \node at (2, 0) [] {,};

    \filldraw[draw=black, fill=white] (2.5, -1) rectangle (5.5, 1);
    \draw[draw=white, draw opacity=0, fill=orange] (2.5, -1) -- (2.5, 1) -- (3, 1) to [out=-90, in=180] (4, 0) to [out=0, in=90] (5, -1) -- cycle;
    \draw[draw=white, draw opacity=0, fill=cyan] (5.5, 1) -- (3, 1) to [out=-90, in=180] (4, 0) to [out=0, in=90] (5, -1) -- (5.5, -1) -- cycle;

    \node (ul2) at (3, 1) [] {};
    \node (ur2) at (5, 1) [] {};
    \node (c2) at (4, 0) [bullet, green] {};
    \node (dl2) at (3, -1) [] {};
    \node (dr2) at (5, -1) [] {};

    \path[draw, line width=2pt, line width=2pt, green] (ul2.center) to [out=-90, in=180] (c2);
    \path[draw, line width=2pt, line width=2pt, red] (ur2.center) to [out=-90, in=0] (c2);
    \path[draw, line width=2pt, line width=2pt, blue] (c2) to [out=180, in=90] (dl2.center);
    \path[draw, line width=2pt, line width=2pt, green] (c2) to [out=0, in=90] (dr2.center);

    \node at (6, 0) [] {\Bigg{)}};
\end{tikzpicture}\end{center}
and since there's no non-trivial 2-cells in $\mathbb{1}$, naturality just means $\color{green} \phi \color{black} = \color{green} \phi \color{black}$; the two axioms translate to the commutativity of the following diagram
\begin{center}\begin{tikzcd}
\colorformula{cyan}{Y} \cdot \colorformula{green}{f} \arrow[r, "\colorformula{red}{\eta^Y} \cdot \colorformula{green}{f}"] \arrow[ddd, "\lambda"'] &
\colorformula{red}{t^Y} \cdot \colorformula{green}{f} \arrow[dddddd, "\colorformula{green}{\phi}" description] &
(\colorformula{red}{t^Y} \cdot \colorformula{red}{t^Y}) \cdot \colorformula{green}{f} \arrow[l, "\colorformula{red}{\mu^Y} \cdot \colorformula{green}{f}"'] \arrow[d, "\alpha"] \\
&& \colorformula{red}{t^Y} \cdot (\colorformula{red}{t^Y} \cdot \colorformula{green}{f}) \arrow[d, "\colorformula{red}{t^Y} \cdot \colorformula{green}{\phi}"] \\
&& \colorformula{red}{t^Y} \cdot (\colorformula{green}{f} \cdot \colorformula{blue}{t^X}) \arrow[dd, "\alpha^{-1}"] \\
\colorformula{green}{f} \arrow[ddd, "\rho^{-1}"'] && \\
&& (\colorformula{red}{t^Y} \cdot \colorformula{green}{f}) \cdot \colorformula{blue}{t^X} \arrow[d, "\colorformula{green}{\phi} \cdot \colorformula{blue}{t^X}"] \\
&& (\colorformula{green}{f} \cdot \colorformula{blue}{t^X}) \cdot \colorformula{blue}{t^X} \arrow[d, "\alpha"] \\
\colorformula{green}{f} \cdot \colorformula{orange}{X} \arrow[r, "\colorformula{green}{f} \cdot \colorformula{blue}{\eta^X}"'] &
\colorformula{green}{f} \cdot \colorformula{blue}{t^X} &
\colorformula{green}{f} \cdot (\colorformula{blue}{t^X} \cdot \colorformula{blue}{t^X}) \arrow[l, "\colorformula{green}{f} \cdot \colorformula{blue}{\mu^X}"]
\end{tikzcd}\end{center}
that is
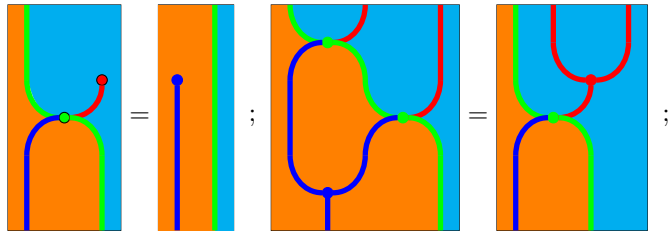
\begin{center}\begin{tikzpicture}[scale=0.5]
    \filldraw[draw=black, fill=white] (-0.5, -3) rectangle (2.5, 3);
    \draw[draw=white, draw opacity=0, fill=orange] (-0.5, 3) -- (0, 3) -- (0, 1) to [out=-90, in=180] (1, 0) to [out=0, in=90] (2, -1) -- (2, -3) -- (-0.5, -3) -- cycle;
    \draw[draw=white, draw opacity=0, fill=cyan] (2.5, 3) -- (0, 3) to [out=-90, in=180] (1, 0) to [out=0, in=90] (2, -1) -- (2, -3) -- (2.5, -3) -- cycle;
    
    \node (uuul1) at (0, 3) [] {};
    \node (ul1) at (0, 1) [] {};
    \node (ur1) at (2, 1) [bullet, fill=red] {};
    \node (c1) at (1, 0) [bullet, fill=green] {};
    \node (dl1) at (0, -1) [] {};
    \node (dr1) at (2, -1) [] {};
    \node (dddl1) at (0, -3) [] {};
    \node (dddr1) at (2, -3) [] {};

    \path[draw, line width=2pt, line width=2pt, green] (uuul1.center) -- (ul1.center);
    \path[draw, line width=2pt, line width=2pt, green] (ul1.center) to [out=-90, in=180] (c1);
    \path[draw, line width=2pt, line width=2pt, red] (ur1) to [out=-90, in=0] (c1);
    \path[draw, line width=2pt, line width=2pt, blue] (c1) to [out=180, in=90] (dl1.center);
    \path[draw, line width=2pt, line width=2pt, green] (c1) to [out=0, in=90] (dr1.center);
    \path[draw, line width=2pt, line width=2pt, blue] (dl1.center) -- (dddl1.center);
    \path[draw, line width=2pt, line width=2pt, green] (dr1.center) -- (dddr1.center);

    \node at (3, 0) [] {=};

    \filldraw[draw=black, fill=white] (3.5, -3) rectangle (5.5, 3);
    \draw[draw=white, draw opacity=0, fill=orange] (3.5, 3) -- (5, 3) -- (5, -3) -- (3.5, -3) -- cycle;
    \draw[draw=white, draw opacity=0, fill=cyan] (5, 3) -- (5.5, 3) -- (5.5, -3) -- (5, -3) -- cycle;

    \node (uuur2) at (5, 3) [] {};
    \node (ul2) at (4, 1) [bullet, blue] {};
    \node (dddl2) at (4, -3) [] {};
    \node (dddr2) at (5, -3) [] {};

    \path[draw, line width=2pt, line width=2pt, green] (uuur2.center) -- (dddr2.center);
    \path[draw, line width=2pt, line width=2pt, blue] (ul2) -- (dddl2.center);

    \node at (6, 0) [] {;};

    \filldraw[draw=black, fill=white] (6.5, -3) rectangle (11.5, 3);
    \draw[draw=white, draw opacity=0, fill=orange] (6.5, 3) -- (7, 3) to [out=-90, in=180] (8, 2) to [out=0, in=90] (9, 1) to [out=-90, in=180] (10, 0) to [out=0, in=90] (11, -1) -- (11, -3) -- (6.5, -3) -- cycle;
    \draw[draw=white, draw opacity=0, fill=cyan] (11.5, 3) -- (7, 3) to [out=-90, in=180] (8, 2) to [out=0, in=90] (9, 1) to [out=-90, in=180] (10, 0) to [out=0, in=90] (11, -1) -- (11, -3) -- (11.5, -3) -- cycle;

    \node (uuull3) at (7, 3) [] {};
    \node (uuuc3) at (9, 3) [] {};
    \node (uuurr3) at (11, 3) [] {};
    \node (uul3) at (8, 2) [bullet, green] {};
    \node (ull3) at (7, 1) [] {};
    \node (uc3) at (9, 1) [] {};
    \node (urr3) at (11, 1) [] {};
    \node (cr3) at (10, 0) [bullet, green] {};
    \node (dll3) at (7, -1) [] {};
    \node (dc3) at (9, -1) [] {};
    \node (drr3) at (11, -1) [] {};
    \node (ddl3) at (8, -2) [bullet, blue] {};
    \node (dddl3) at (8, -3) [] {};
    \node (dddrr3) at (11, -3) [] {};

    \path[draw, line width=2pt, line width=2pt, green] (uuull3.center) to [out=-90, in=180] (uul3);
    \path[draw, line width=2pt, line width=2pt, red] (uuuc3.center) to [out=-90, in=0] (uul3);
    \path[draw, line width=2pt, line width=2pt, red] (uuurr3.center) -- (urr3.center);
    \path[draw, line width=2pt, line width=2pt, blue] (uul3) to [out=180, in=90] (ull3.center);
    \path[draw, line width=2pt, line width=2pt, green] (uul3) to [out=0, in=90] (uc3.center);
    \path[draw, line width=2pt, line width=2pt, blue] (ull3.center) -- (dll3.center);
    \path[draw, line width=2pt, line width=2pt, green] (uc3.center) to [out=-90, in=180] (cr3);
    \path[draw, line width=2pt, line width=2pt, red] (urr3.center) to [out=-90, in=0] (cr3);
    \path[draw, line width=2pt, line width=2pt, blue] (cr3) to [out=180, in=90] (dc3.center);
    \path[draw, line width=2pt, line width=2pt, green] (cr3) to [out=0, in=90] (drr3.center);
    \path[draw, line width=2pt, line width=2pt, blue] (dll3.center) to [out=-90, in=180] (ddl3);
    \path[draw, line width=2pt, line width=2pt, blue] (dc3.center) to [out=-90, in=0] (ddl3);
    \path[draw, line width=2pt, line width=2pt, green] (drr3.center) -- (dddrr3.center);
    \path[draw, line width=2pt, line width=2pt, blue] (ddl3.center) -- (dddl3.center);

    \node at (12, 0) [] {=};

    \filldraw[draw=black, fill=white] (12.5, -3) rectangle (16.5, 3);
    \draw[draw=white, draw opacity=0, fill=orange] (12.5, 3) -- (13, 3) -- (13, 1) to [out=-90, in=180] (14, 0) to [out=0, in=90] (15, -1) -- (15, -3) -- (12.5, -3) -- cycle;
    \draw[draw=white, draw opacity=0, fill=cyan] (16.5, 3) -- (13, 3) -- (13, 1) to [out=-90, in=180] (14, 0) to [out=0, in=90] (15, -1) -- (15, -3) -- (16.5, -3) -- cycle;

    \node (uuull4) at (13, 3) [] {};
    \node (uuul4) at (14, 3) [] {};
    \node (uuurr4) at (16, 3) [] {};
    \node (uul4) at (14, 2) [] {};
    \node (uurr4) at (16, 2) [] {};
    \node (ull4) at (13, 1) [] {};
    \node (ur4) at (15, 1) [bullet, red] {};
    \node (cl4) at (14, 0) [bullet, green] {};
    \node (dll4) at (13, -1) [] {};
    \node (dr4) at (15, -1) [] {};
    \node (dddll4) at (13, -3) [] {};
    \node (dddr4) at (15, -3) [] {};

    \path[draw, line width=2pt, line width=2pt, green] (uuull4.center) -- (ull4.center);
    \path[draw, line width=2pt, line width=2pt, red] (uuul4.center) -- (uul4.center);
    \path[draw, line width=2pt, red] (uuurr4.center) -- (uurr4.center);
    \path[draw, line width=2pt, red] (uul4.center) to [out=-90, in=180] (ur4);
    \path[draw, line width=2pt, red] (uurr4.center) to [out=-90, in=0] (ur4);
    \path[draw, line width=2pt, green] (ull4.center) to [out=-90, in=180] (cl4);
    \path[draw, line width=2pt, red] (ur4) to [out=-90, in=0] (cl4);
    \path[draw, line width=2pt, blue] (cl4) to [out=180, in=90] (dll4.center);
    \path[draw, line width=2pt, green] (cl4) to [out=0, in=90] (dr4.center);
    \path[draw, line width=2pt, blue] (dll4.center) -- (dddll4.center);
    \path[draw, line width=2pt, green] (dr4.center) -- (dddr4.center);

    \node at (17, 0) [] {;};
\end{tikzpicture}\end{center}

\item A 2-cell $\gamma : (\colorformula{green}{f}, \colorformula{green}{\phi}) \to (\colorformula{violet}{g}, \colorformula{violet}{\psi}) : (\colorformula{orange}{X}, \colorformula{blue}{t^X}, \colorformula{blue}{\eta^X}, \colorformula{blue}{\mu^X}) \to (\colorformula{cyan}{Y}, \colorformula{red}{t^Y}, \colorformula{red}{\eta^Y}, \colorformula{red}{\mu^Y}) : \mathbf{Mnd}(\mathcal{B})$ consists of a 2-cell $\gamma : \colorformula{green}{f} \to \colorformula{violet}{g} : \colorformula{orange}{X} \to \colorformula{cyan}{Y} : \mathcal{B}$, depicted as
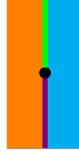
\begin{center}\begin{tikzpicture}
    \filldraw[draw=black, fill=white] (-0.5, -1) rectangle (+0.5, 1);
    \draw[draw=white, draw opacity=0, fill=orange] (-0.5, 1) -- (0, 1) -- (0, -1) -- (-0.5, -1) -- cycle;
    \draw[draw=white, draw opacity=0, fill=cyan] (0.5, 1) -- (0, 1) -- (0, -1) -- (0.5, -1)-- cycle;

    \node (u) at (0, 1) [] {};
    \node (c) at (0, 0) [bullet] {};
    \node (d) at (0, -1) [] {};

    \path[draw, line width=2pt, green] (u.center) -- (c);
    \path[draw, line width=2pt, violet] (c) -- (d.center);
\end{tikzpicture}\end{center}
such that
\begin{center}\begin{tikzcd}
\colorformula{red}{t^Y} \cdot \colorformula{green}{f} \arrow[r, "\colorformula{green}{\phi}"] \arrow[d, "\colorformula{red}{t^Y} \cdot \gamma"'] &
\colorformula{green}{f} \cdot \colorformula{blue}{t^X} \arrow[d, "\gamma \cdot \colorformula{blue}{t^X}"] \\
\colorformula{red}{t^Y} \cdot \colorformula{violet}{g} \arrow[r, "\colorformula{violet}{\psi}"'] &
\colorformula{violet}{g} \cdot \colorformula{blue}{t^X}
\end{tikzcd}\end{center}
that is
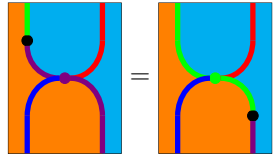
\begin{center}\begin{tikzpicture}[scale=0.5]
    \filldraw[draw=black, fill=white] (-0.5, -2) rectangle (2.5, 2);
    \draw[draw=white, draw opacity=0, fill=orange] (-0.5, 2) -- (0, 2) -- (0, 1) to [out=-90, in=180] (1, 0) to [out=0, in=90] (2, -1) -- (2, -2) -- (-0.5, -2) -- cycle;
    \draw[draw=white, draw opacity=0, fill=cyan] (2.5, 2) -- (0, 2) -- (0, 1) to [out=-90, in=180] (1, 0) to [out=0, in=90] (2, -1) -- (2, -2) -- (2.5, -2) -- cycle;

    \node (uul1) at (0, 2) [] {};
    \node (uur1) at (2, 2) [] {};
    \node (ul1) at (0, 1) [bullet] {};
    \node (ur1) at (2, 1) [] {};
    \node (cc1) at (1, 0) [bullet, violet] {};
    \node (dl1) at (0, -1) [] {};
    \node (dr1) at (2, -1) [] {};
    \node (ddl1) at (0, -2) [] {};
    \node (ddr1) at (2, -2) [] {};

    \path[draw, line width=2pt, green] (uul1.center) -- (ul1);
    \path[draw, line width=2pt, red] (uur1.center) -- (ur1.center);
    \path[draw, line width=2pt, violet] (ul1) to [out=-90, in=180] (cc1);
    \path[draw, line width=2pt, red] (ur1.center) to [out=-90, in=0] (cc1);
    \path[draw, line width=2pt, blue] (cc1) to [out=180, in=90] (dl1.center);
    \path[draw, line width=2pt, violet] (cc1) to [out=0, in=90] (dr1.center);
    \path[draw, line width=2pt, blue] (dl1.center) -- (ddl1.center);
    \path[draw, line width=2pt, violet] (dr1.center) -- (ddr1.center);

    \node at (3, 0) [] {=};

    \filldraw[draw=black, fill=white] (3.5, -2) rectangle (6.5, 2);
    \draw[draw=white, draw opacity=0, fill=orange] (3.5, 2) -- (4, 2) -- (4, 1) to [out=-90, in=180] (5, 0) to [out=0, in=90] (6, -1) -- (6, -2) -- (3.5, -2) -- cycle;
    \draw[draw=white, draw opacity=0, fill=cyan] (6.5, 2) -- (4, 2) -- (4, 1) to [out=-90, in=180] (5, 0) to [out=0, in=90] (6, -1) -- (6, -2) -- (6.5, -2) -- cycle;

    \node (uul2) at (4, 2) [] {};
    \node (uur2) at (6, 2) [] {};
    \node (ul2) at (4, 1) [] {};
    \node (ur2) at (6, 1) [] {};
    \node (cc2) at (5, 0) [bullet, green] {};
    \node (dl2) at (4, -1) [] {};
    \node (dr2) at (6, -1) [bullet] {};
    \node (ddl2) at (4, -2) [] {};
    \node (ddr2) at (6, -2) [] {};

    \path[draw, line width=2pt, green] (uul2.center) -- (ul2.center);
    \path[draw, line width=2pt, red] (uur2.center) -- (ur2.center);
    \path[draw, line width=2pt, green] (ul2.center) to [out=-90, in=180] (cc2);
    \path[draw, line width=2pt, red] (ur2.center) to [out=-90, in=0] (cc2);
    \path[draw, line width=2pt, blue] (cc2) to [out=180, in=90] (dl2.center);
    \path[draw, line width=2pt, green] (cc2) to [out=0, in=90] (dr2);
    \path[draw, line width=2pt, blue] (dl2.center) -- (ddl2.center);
    \path[draw, line width=2pt, violet] (dr2) -- (ddr2.center);
\end{tikzpicture}\end{center}
\end{enumerate}

Another central notion in the theory of monads is that of a distributive law. This too can be internalized in an arbitrary 2-category

\begin{defin}\thmheader{Distributive law}

Let $\mathcal{B}$ be a 2-category and $(X, \colorformula{blue}{t_1}, \colorformula{blue}{\eta_1}, \colorformula{blue}{\mu_1}), (X, \colorformula{red}{t_2}, \colorformula{red}{\eta_2}, \colorformula{red}{\mu_2})$ be two monads on the same 0-cell. A distributive law
\[\gamma : (X, \colorformula{blue}{t_1}, \colorformula{blue}{\eta_1}, \colorformula{blue}{\mu_1}) \rightsquigarrow (X, \colorformula{red}{t_2}, \colorformula{red}{\eta_2}, \colorformula{red}{\mu_2})\]
is given by a 2-cell $\gamma : \colorformula{red}{t_2} \cdot \colorformula{blue}{t_1} \to \colorformula{blue}{t_1} \cdot \colorformula{red}{t_2}$, depicted\footnote{Anytime there's only one 0-cell involved, we'll depict it as a white background.}
\begin{center}\begin{tikzpicture}[scale=0.5]
    \filldraw[draw=black, fill=white] (-0.5, -2) rectangle (2.5, 2);

    \node (uul) at (0, 2) [] {};
    \node (uur) at (2, 2) [] {};
    \node (ul) at (0, 1) [] {};
    \node (ur) at (2, 1) [] {};
    \node (c) at (1, 0) [bullet] {};
    \node (dl) at (0, -1) [] {};
    \node (dr) at (2, -1) [] {};
    \node (ddl) at (0, -2) [] {};
    \node (ddr) at (2, -2) [] {};

    \path[draw, line width=2pt, blue] (uul.center) -- (ul.center);
    \path[draw, line width=2pt, red] (uur.center) -- (ur.center);
    \path[draw, line width=2pt, blue] (ul.center) to [out=-90, in=180] (c);
    \path[draw, line width=2pt, red] (ur.center) to [out=-90, in=0] (c);
    \path[draw, line width=2pt, red] (c) to [out=180, in=90] (dl.center);
    \path[draw, line width=2pt, blue] (c) to [out=0, in=90] (dr.center);
    \path[draw, line width=2pt, red] (dl.center) -- (ddl.center);
    \path[draw, line width=2pt, blue] (dr.center) -- (ddr.center);
\end{tikzpicture}\end{center}
making the following diagrams commute
\begin{center}\begin{tikzcd}
X \cdot \colorformula{blue}{t_1} \arrow[r, "\colorformula{red}{\eta_2} \cdot \colorformula{blue}{t_1}"] &
\colorformula{red}{t_2} \cdot \colorformula{blue}{t_1} \arrow[dd, "\gamma" description] &
\colorformula{red}{t_2} \cdot X \arrow[l, "\colorformula{red}{t_2} \cdot \colorformula{blue}{\eta_1}"'] \\
\colorformula{blue}{t_1} \arrow[u, "\lambda^{-1}"] \arrow[d, "\rho^{-1}"'] &&
\colorformula{red}{t_2} \arrow[u, "\rho^{-1}"'] \arrow[d, "\lambda^{-1}"] \\
\colorformula{blue}{t_1} \cdot X \arrow[r, "\colorformula{blue}{t_1} \cdot \colorformula{red}{\eta_2}"'] &
\colorformula{blue}{t_1} \cdot \colorformula{red}{t_2} &
X \cdot \colorformula{red}{t_2} \arrow[l, "\colorformula{blue}{\eta_1} \cdot \colorformula{red}{t_2}"]
\end{tikzcd}\end{center}
\begin{center}\begin{tikzcd}
\colorformula{red}{t_2} \cdot (\colorformula{blue}{t_1} \cdot \colorformula{blue}{t_1}) \arrow[r, "\colorformula{red}{t_2} \cdot \colorformula{blue}{\mu_1}"] \arrow[d, "\alpha^{-1}"'] &
\colorformula{red}{t_2} \cdot \colorformula{blue}{t_1} \arrow[ddddd, "\gamma" description] &
(\colorformula{red}{t_2} \cdot \colorformula{red}{t_2}) \cdot \colorformula{blue}{t_1} \arrow[l, "\colorformula{red}{\mu_2} \cdot \colorformula{blue}{t_1}"'] \arrow[d, "\alpha"] \\
(\colorformula{red}{t_2} \cdot \colorformula{blue}{t_1}) \cdot \colorformula{blue}{t_1} \arrow[d, "\gamma \cdot \colorformula{blue}{t_1}"'] &&
\colorformula{red}{t_2} \cdot (\colorformula{red}{t_2} \cdot \colorformula{blue}{t_1}) \arrow[d, "\colorformula{red}{t_2} \cdot \gamma"] \\
(\colorformula{blue}{t_1} \cdot \colorformula{red}{t_2}) \cdot \colorformula{blue}{t_1} \arrow[d, "\alpha"'] &&
\colorformula{red}{t_2} \cdot (\colorformula{blue}{t_1} \cdot \colorformula{red}{t_2}) \arrow[d, "\alpha^{-1}"] \\
\colorformula{blue}{t_1} \cdot (\colorformula{red}{t_2} \cdot \colorformula{blue}{t_1}) \arrow[d, "\colorformula{blue}{t_1} \cdot \gamma"'] &&
(\colorformula{red}{t_2} \cdot \colorformula{blue}{t_1}) \cdot \colorformula{red}{t_2} \arrow[d, "\gamma \cdot \colorformula{red}{t_2}"] \\
\colorformula{blue}{t_1} \cdot (\colorformula{blue}{t_1} \cdot \colorformula{red}{t_2}) \arrow[d, "\alpha^{-1}"'] &&
(\colorformula{blue}{t_1} \cdot \colorformula{red}{t_2}) \cdot \colorformula{red}{t_2} \arrow[d, "\alpha"] \\
(\colorformula{blue}{t_1} \cdot \colorformula{blue}{t_1}) \cdot \colorformula{red}{t_2} \arrow[r, "\colorformula{blue}{\mu_1} \cdot \colorformula{red}{t_2}"'] &
\colorformula{blue}{t_1} \cdot \colorformula{red}{t_2} &
\colorformula{blue}{t_1} \cdot (\colorformula{red}{t_2} \cdot \colorformula{red}{t_2}) \arrow[l, "\colorformula{blue}{t_1} \cdot \colorformula{red}{\mu_2}"]
\end{tikzcd}\end{center}
that is (the first two equalities encode the first diagram, the second two the second)
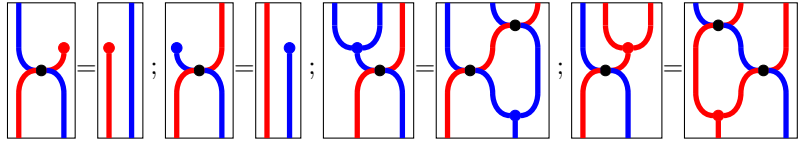
\begin{center}\begin{tikzpicture}[scale=0.3]
    \filldraw[draw=black, fill=white] (-0.5, -3) rectangle (2.5, 3);

    \node (uuul1) at (0, 3) [] {};
    \node (ul1) at (0, 1) [] {};
    \node (ur1) at (2, 1) [bullet, red] {};
    \node (c1) at (1, 0) [bullet] {};
    \node (dl1) at (0, -1) [] {};
    \node (dr1) at (2, -1) [] {};
    \node (dddl1) at (0, -3) [] {};
    \node (dddr1) at (2, -3) [] {};

    \path[draw, line width=2pt, blue] (uuul1.center) -- (ul1.center);
    \path[draw, line width=2pt, blue] (ul1.center) to [out=-90, in=180] (c1);
    \path[draw, line width=2pt, red] (ur1) to [out=-90, in=0] (c1);
    \path[draw, line width=2pt, red] (c1) to [out=180, in=90] (dl1.center);
    \path[draw, line width=2pt, blue] (c1) to [out=0, in=90] (dr1.center);
    \path[draw, line width=2pt, red] (dl1.center) -- (dddl1.center);
    \path[draw, line width=2pt, blue] (dr1.center) -- (dddr1.center);

    \node at (3, 0) [] {=};

    \filldraw[draw=black, fill=white] (3.5, -3) rectangle (5.5, 3);

    \node (uuur2) at (5, 3) [] {};
    \node (ul2) at (4, 1) [bullet, red] {};
    \node (dddl2) at (4, -3) [] {};
    \node (dddr2) at (5, -3) [] {};

    \path[draw, line width=2pt, blue] (uuur2.center) -- (dddr2.center);
    \path[draw, line width=2pt, red] (ul2) -- (dddl2.center);

    \node at (6, 0) [] {;};

    \filldraw[draw=black, fill=white] (6.5, -3) rectangle (9.5, 3);

    \node (uuur3) at (9, 3) [] {};
    \node (ul3) at (7, 1) [bullet, blue] {};
    \node (ur3) at (9, 1) [] {};
    \node (c3) at (8, 0) [bullet] {};
    \node (dl3) at (7, -1) [] {};
    \node (dr3) at (9, -1) [] {};
    \node (dddl3) at (7, -3) [] {};
    \node (dddr3) at (9, -3) [] {};

    \path[draw, line width=2pt, red] (uuur3.center) -- (ur3.center);
    \path[draw, line width=2pt, blue] (ul3) to [out=-90, in=180] (c3);
    \path[draw, line width=2pt, red] (ur3.center) to [out=-90, in=0] (c3);
    \path[draw, line width=2pt, red] (c3) to [out=180, in=90] (dl3.center);
    \path[draw, line width=2pt, blue] (c3) to [out=0, in=90] (dr3.center);
    \path[draw, line width=2pt, red] (dl3.center) -- (dddl3.center);
    \path[draw, line width=2pt, blue] (dr3.center) -- (dddr3.center);

    \node at (10, 0) [] {=};
    
    \filldraw[draw=black, fill=white] (10.5, -3) rectangle (12.5, 3);

    \node (uuul4) at (11, 3) [] {};
    \node (ur4) at (12, 1) [bullet, blue] {};
    \node (dddl4) at (11, -3) [] {};
    \node (dddr4) at (12, -3) [] {};

    \path[draw, line width=2pt, red] (uuul4.center) -- (dddl4.center);
    \path[draw, line width=2pt, blue] (ur4) -- (dddr4.center);

    \node at (13, 0) [] {;};
    
    \filldraw[draw=black, fill=white] (13.5, -3) rectangle (17.5, 3);

    \node (uuull5) at (14, 3) [] {};
    \node (uuur5) at (16, 3) [] {};
    \node (uuurr5) at (17, 3) [] {};
    \node (uull5) at (14, 2) [] {};
    \node (uur5) at (16, 2) [] {};
    \node (ul5) at (15, 1) [bullet, blue] {};
    \node (urr5) at (17, 1) [] {};
    \node (cr5) at (16, 0) [bullet] {};
    \node (dl5) at (15, -1) [] {};
    \node (drr5) at (17, -1) [] {};
    \node (dddl5) at (15, -3) [] {};
    \node (dddrr5) at (17, -3) [] {};

    \path[draw, line width=2pt, blue] (uuull5.center) -- (uull5.center);
    \path[draw, line width=2pt, blue] (uuur5.center) -- (uur5.center);
    \path[draw, line width=2pt, red] (uuurr5.center) -- (urr5.center);
    \path[draw, line width=2pt, blue] (uull5.center) to [out=-90, in=180] (ul5);
    \path[draw, line width=2pt, blue] (uur5.center) to [out=-90, in=0] (ul5);
    \path[draw, line width=2pt, blue] (ul5) to [out=-90, in=180] (cr5);
    \path[draw, line width=2pt, red] (urr5.center) to [out=-90, in=0] (cr5);
    \path[draw, line width=2pt, red] (cr5) to [out=180, in=90] (dl5.center);
    \path[draw, line width=2pt, blue] (cr5) to [out=0, in=90] (drr5.center);
    \path[draw, line width=2pt, red] (dl5.center) -- (dddl5.center);
    \path[draw, line width=2pt, blue] (drr5.center) -- (dddrr5.center);

    \node at (18, 0) [] {=};
    
    \filldraw[draw=black, fill=white] (18.5, -3) rectangle (23.5, 3);

    \node (uuull6) at (19, 3) [] {};
    \node (uuuc6) at (21, 3) [] {};
    \node (uuurr6) at (23, 3) [] {};
    \node (uur6) at (22, 2) [bullet] {};
    \node (ull6) at (19, 1) [] {};
    \node (uc6) at (21, 1) [] {};
    \node (urr6) at (23, 1) [] {};
    \node (cl6) at (20, 0) [bullet] {};
    \node (dll6) at (19, -1) [] {};
    \node (dc6) at (21, -1) [] {};
    \node (drr6) at (23, -1) [] {};
    \node (ddr6) at (22, -2) [bullet, blue] {};
    \node (dddll6) at (19, -3) [] {};
    \node (dddr6) at (22, -3) [] {};

    \path[draw, line width=2pt, blue] (uuull6.center) -- (ull6.center);
    \path[draw, line width=2pt, blue] (uuuc6.center) to [out=-90, in=180] (uur6);
    \path[draw, line width=2pt, red] (uuurr6.center) to [out=-90, in=0] (uur6);
    \path[draw, line width=2pt, red] (uur6) to [out=180, in=90] (uc6.center);
    \path[draw, line width=2pt, blue] (uur6) to [out=0, in=90] (urr6.center);
    \path[draw, line width=2pt, blue] (ull6.center) to [out=-90, in=180] (cl6);
    \path[draw, line width=2pt, red] (uc6.center) to [out=-90, in=0] (cl6);
    \path[draw, line width=2pt, blue] (urr6.center) -- (drr6.center);
    \path[draw, line width=2pt, red] (cl6) to [out=180, in=90] (dll6.center);
    \path[draw, line width=2pt, blue] (cl6) to [out=0, in=90] (dc6.center);
    \path[draw, line width=2pt, red] (dll6.center) -- (dddll6.center);
    \path[draw, line width=2pt, blue] (dc6.center) to [out=-90, in=180] (ddr6);
    \path[draw, line width=2pt, blue] (drr6.center) to [out=-90, in=0] (ddr6);
    \path[draw, line width=2pt, blue] (ddr6) -- (dddr6.center);

    \node at (24, 0) [] {;};
    
    \filldraw[draw=black, fill=white] (24.5, -3) rectangle (28.5, 3);

    \node (uuull7) at (25, 3) [] {};
    \node (uuul7) at (26, 3) [] {};
    \node (uuurr7) at (28, 3) [] {};
    \node (uul7) at (26, 2) [] {};
    \node (uurr7) at (28, 2) [] {};
    \node (ull7) at (25, 1) [] {};
    \node (ur7) at (27, 1) [bullet, red] {};
    \node (cl7) at (26, 0) [bullet] {};
    \node (dll7) at (25, -1) [] {};
    \node (dr7) at (27, -1) [] {};
    \node (dddll7) at (25, -3) [] {};
    \node (dddr7) at (27, -3) [] {};

    \path[draw, line width=2pt, blue] (uuull7.center) -- (ull7.center);
    \path[draw, line width=2pt, red] (uuul7.center) -- (uul7.center);
    \path[draw, line width=2pt, red] (uuurr7.center) -- (uurr7.center);
    \path[draw, line width=2pt, red] (uul7.center) to [out=-90, in=180] (ur7);
    \path[draw, line width=2pt, red] (uurr7.center) to [out=-90, in=0] (ur7);
    \path[draw, line width=2pt, blue] (ull7.center) to [out=-90, in=180] (cl7);
    \path[draw, line width=2pt, red] (ur7) to [out=-90, in=0] (cl7);
    \path[draw, line width=2pt, red] (cl7) to [out=180, in=90] (dll7.center);
    \path[draw, line width=2pt, blue] (cl7) to [out=0, in=90] (dr7.center);
    \path[draw, line width=2pt, red] (dll7.center) -- (dddll7.center);
    \path[draw, line width=2pt, blue] (dr7.center) -- (dddr7.center);

    \node at (29, 0) [] {=};
    
    \filldraw[draw=black, fill=white] (29.5, -3) rectangle (34.5, 3);

    \node (uuull8) at (30, 3) [] {};
    \node (uuuc8) at (32, 3) [] {};
    \node (uuurr8) at (34, 3) [] {};
    \node (uul8) at (31, 2) [bullet] {};
    \node (ull8) at (30, 1) [] {};
    \node (uc8) at (32, 1) [] {};
    \node (urr8) at (34, 1) [] {};
    \node (cr8) at (33, 0) [bullet] {};
    \node (dll8) at (30, -1) [] {};
    \node (dc8) at (32, -1) [] {};
    \node (drr8) at (34, -1) [] {};
    \node (ddl8) at (31, -2) [bullet, red] {};
    \node (dddl8) at (31, -3) [] {};
    \node (dddrr8) at (34, -3) [] {};

    \path[draw, line width=2pt, blue] (uuull8.center) to [out=-90, in=180] (uul8);
    \path[draw, line width=2pt, red] (uuuc8.center) to [out=-90, in=0] (uul8);
    \path[draw, line width=2pt, red] (uuurr8.center) -- (urr8.center);
    \path[draw, line width=2pt, red] (uul8) to [out=180, in=90] (ull8.center);
    \path[draw, line width=2pt, blue] (uul8) to [out=0, in=90] (uc8.center);
    \path[draw, line width=2pt, red] (ull8.center) -- (dll8.center);
    \path[draw, line width=2pt, blue] (uc8.center) to [out=-90, in=180] (cr8);
    \path[draw, line width=2pt, red] (urr8.center) to [out=-90, in=0] (cr8);
    \path[draw, line width=2pt, red] (cr8) to [out=180, in=90] (dc8.center);
    \path[draw, line width=2pt, blue] (cr8) to [out=0, in=90] (drr8.center);
    \path[draw, line width=2pt, red] (dll8.center) to [out=-90, in=180] (ddl8);
    \path[draw, line width=2pt, red] (dc8.center) to [out=-90, in=0] (ddl8);
    \path[draw, line width=2pt, blue] (drr8.center) -- (dddrr8.center);
    \path[draw, line width=2pt, red] (ddl8) -- (dddl8.center);

    \node at (35, 0) [] {;};
\end{tikzpicture}\end{center}
\end{defin}

In the case of monads over categories (that is, monads in the 2-category of categories), it's well-known that distributive laws correspond 1-to-1 to lifts of one monad to the Eilenberg-Moore category for the other, allowing one to compose the two monads (for the precise statement, see the first proposition in \cite{distributive-laws}). This fails in general 2-categories\footnote{Though under the appropriate assumptions, that is existence of Eilenberg-Moore objects, the same theorem still holds: the proof is almost identical to that in \cite{distributive-laws}}, but one of the two implications still holds in full generality

\begin{prop}\thmheader{Distributive laws allow for monad composition}

In a 2-category $\mathcal{B}$, fix a distributive law between monads
\[\gamma : (X, \colorformula{blue}{t_1}, \colorformula{blue}{\eta_1}, \colorformula{blue}{\mu_1}) \rightsquigarrow (X, \colorformula{red}{t_2}, \colorformula{red}{\eta_2}, \colorformula{red}{\mu_2})\]
We can define a new monad $\colorformula{purple}{\mathcal{T}} = (X, \colorformula{purple}{t}, \colorformula{purple}{\eta}, \colorformula{purple}{\mu})$ with $\colorformula{purple}{t} = \colorformula{blue}{t_1} \cdot \colorformula{red}{t_2}$.

\begin{proof}

First, we define the monad $\colorformula{purple}{\mathcal{T}}$
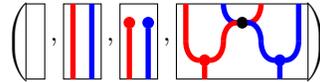
\begin{center}\begin{tikzpicture}[scale=0.25]
    \node at (0, 0) [] {\Bigg{(}};

    \filldraw[draw=black, fill=white] (0.5, -2) rectangle (1.5, 2);

    \node at (2, 0) [] {,};

    \filldraw[draw=black, fill=white] (2.5, -2) rectangle (4.5, 2);

    \node (uul1) at (3, 2) [] {};
    \node (uur1) at (4, 2) [] {};
    \node (ddl1) at (3, -2) [] {};
    \node (ddr1) at (4, -2) [] {};

    \path[draw, line width=2pt, red] (uul1.center) -- (ddl1.center);
    \path[draw, line width=2pt, blue] (uur1.center) -- (ddr1.center);

    \node at (5, 0) [] {,};

    \filldraw[draw=black, fill=white] (5.5, -2) rectangle (7.5, 2);

    \node (ul2) at (6, 1) [bullet, red] {};
    \node (ur2) at (7, 1) [bullet, blue] {};
    \node (ddl2) at (6, -2) [] {};
    \node (ddr2) at (7, -2) [] {};

    \path[draw, line width=2pt, red] (ul2) -- (ddl2.center);
    \path[draw, line width=2pt, blue] (ur2) -- (ddr2.center);

    \node at (8, 0) [] {,};

    \filldraw[draw=black, fill=white] (8.5, -2) rectangle (15.5, 2);

    \node (uulll3) at (9, 2) [] {};
    \node (uul3) at (11, 2) [] {};
    \node (uur3) at (13, 2) [] {};
    \node (uurrr3) at (15, 2) [] {};
    \node (uc3) at (12, 1) [bullet] {};
    \node (clll3) at (9, 0) [] {};
    \node (cl3) at (11, 0) [] {};
    \node (cr3) at (13, 0) [] {};
    \node (crrr3) at (15, 0) [] {};
    \node (dll3) at (10, -1) [bullet, red] {};
    \node (drr3) at (14, -1) [bullet, blue] {};
    \node (ddll3) at (10, -2) [] {};
    \node (ddrr3) at (14, -2) [] {};

    \path[draw, line width=2pt, red] (uulll3.center) -- (clll3.center);
    \path[draw, line width=2pt, blue] (uul3.center) to [out=-90, in=180] (uc3);
    \path[draw, line width=2pt, red] (uur3.center) to [out=-90, in=0] (uc3);
    \path[draw, line width=2pt, blue] (uurrr3.center) -- (crrr3.center);
    \path[draw, line width=2pt, red] (uc3) to [out=180, in=90] (cl3.center);
    \path[draw, line width=2pt, blue] (uc3) to [out=0, in=90] (cr3.center);
    \path[draw, line width=2pt, red] (clll3.center) to [out=-90, in=180] (dll3);
    \path[draw, line width=2pt, red] (cl3.center) to [out=-90, in=0] (dll3);
    \path[draw, line width=2pt, blue] (cr3.center) to [out=-90, in=180] (drr3);
    \path[draw, line width=2pt, blue] (crrr3.center) to [out=-90, in=0] (drr3);
    \path[draw, line width=2pt, red] (dll3) -- (ddll3.center);
    \path[draw, line width=2pt, blue] (drr3) -- (ddrr3.center);

    \node at (16, 0) [] {\Bigg{)}};
\end{tikzpicture}\end{center}

Showing $\colorformula{violet}{\mathcal{T}}$ is a monad amounts to proving unitality and associativity: using string diagrams this is straightforward. The following string diagram calculation proves unitality
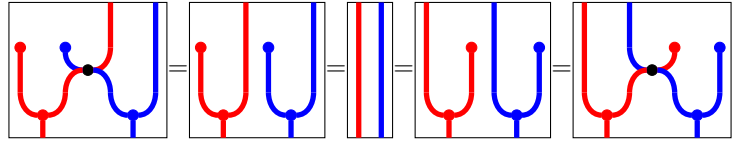
\begin{center}\begin{tikzpicture}[scale=0.3]
    \filldraw[draw=black, fill=white] (-0.5, -3) rectangle (6.5, 3);

    \node (uuur1) at (4, 3) [] {};
    \node (uuurrr1) at (6, 3) [] {};
    \node (ulll1) at (0, 1) [bullet, red] {};
    \node (ul1) at (2, 1) [bullet, blue] {};
    \node (ur1) at (4, 1) [] {};
    \node (cc1) at (3, 0) [bullet] {};
    \node (dlll1) at (0, -1) [] {};
    \node (dl1) at (2, -1) [] {};
    \node (dr1) at (4, -1) [] {};
    \node (drrr1) at (6, -1) [] {};
    \node (ddll1) at (1, -2) [bullet, red] {};
    \node (ddrr1) at (5, -2) [bullet, blue] {};
    \node (dddll1) at (1, -3) [] {};
    \node (dddrr1) at (5, -3) [] {};

    \path[draw, line width=2pt, red] (uuur1.center) -- (ur1.center);
    \path[draw, line width=2pt, blue] (uuurrr1.center) -- (drrr1.center);
    \path[draw, line width=2pt, red] (ulll1) -- (dlll1.center);
    \path[draw, line width=2pt, blue] (ul1) to [out=-90, in=180] (cc1);
    \path[draw, line width=2pt, red] (ur1.center) to [out=-90, in=0] (cc1);
    \path[draw, line width=2pt, red] (cc1) to [out=180, in=90] (dl1.center);
    \path[draw, line width=2pt, blue] (cc1) to [out=0, in=90] (dr1.center);
    \path[draw, line width=2pt, red] (dlll1.center) to [out=-90, in=180] (ddll1);
    \path[draw, line width=2pt, red] (dl1.center) to [out=-90, in=0] (ddll1);
    \path[draw, line width=2pt, blue] (dr1.center) to [out=-90, in=180] (ddrr1);
    \path[draw, line width=2pt, blue] (drrr1.center) to [out=-90, in=0] (ddrr1);
    \path[draw, line width=2pt, red] (ddll1) -- (dddll1.center);
    \path[draw, line width=2pt, blue] (ddrr1) -- (dddrr1.center);

    \node at (7, 0) [] {=};

    \filldraw[draw=black, fill=white] (7.5, -3) rectangle (13.5, 3);

    \node (uuul2) at (10, 3) [] {};
    \node (uuurrr2) at (13, 3) [] {};
    \node (ulll2) at (8, 1) [bullet, red] {};
    \node (ur2) at (11, 1) [bullet, blue] {};
    \node (dlll2) at (8, -1) [] {};
    \node (dl2) at (10, -1) [] {};
    \node (dr2) at (11, -1) [] {};
    \node (drrr2) at (13, -1) [] {};
    \node (ddll2) at (9, -2) [bullet, red] {};
    \node (ddrr2) at (12, -2) [bullet, blue] {};
    \node (dddll2) at (9, -3) [] {};
    \node (dddrr2) at (12, -3) [] {};

    \path[draw, line width=2pt, red] (uuul2.center) -- (dl2.center);
    \path[draw, line width=2pt, blue] (uuurrr2.center) -- (drrr2.center);
    \path[draw, line width=2pt, red] (ulll2) -- (dlll2.center);
    \path[draw, line width=2pt, blue] (ur2) -- (dr2.center);
    \path[draw, line width=2pt, red] (dlll2.center) to [out=-90, in=180] (ddll2);
    \path[draw, line width=2pt, red] (dl2.center) to [out=-90, in=0] (ddll2);
    \path[draw, line width=2pt, blue] (dr2.center) to [out=-90, in=180] (ddrr2);
    \path[draw, line width=2pt, blue] (drrr2.center) to [out=-90, in=0] (ddrr2);
    \path[draw, line width=2pt, red] (ddll2) -- (dddll2.center);
    \path[draw, line width=2pt, blue] (ddrr2) -- (dddrr2.center);

    \node at (14, 0) [] {=};

    \filldraw[draw=black, fill=white] (14.5, -3) rectangle (16.5, 3);

    \node (uuul3) at (15, 3) [] {};
    \node (uuur3) at (16, 3) [] {};
    \node (dddl3) at (15, -3) [] {};
    \node (dddr3) at (16, -3) [] {};

    \path[draw, line width=2pt, red] (uuul3.center) -- (dddl3.center);
    \path[draw, line width=2pt, blue] (uuur3.center) -- (dddr3.center);

    \node at (17, 0) [] {=};

    \filldraw[draw=black, fill=white] (17.5, -3) rectangle (23.5, 3);

    \node (uuulll4) at (18, 3) [] {};
    \node (uuur4) at (21, 3) [] {};
    \node (ul4) at (20, 1) [bullet, red] {};
    \node (urrr4) at (23, 1) [bullet, blue] {};
    \node (dlll4) at (18, -1) [] {};
    \node (dl4) at (20, -1) [] {};
    \node (dr4) at (21, -1) [] {};
    \node (drrr4) at (23, -1) [] {};
    \node (ddll4) at (19, -2) [bullet, red] {};
    \node (ddrr4) at (22, -2) [bullet, blue] {};
    \node (dddll4) at (19, -3) [] {};
    \node (dddrr4) at (22, -3) [] {};

    \path[draw, line width=2pt, red] (uuulll4.center) -- (dlll4.center);
    \path[draw, line width=2pt, blue] (uuur4.center) -- (dr4.center);
    \path[draw, line width=2pt, red] (ul4) -- (dl4.center);
    \path[draw, line width=2pt, blue] (urrr4) -- (drrr4.center);
    \path[draw, line width=2pt, red] (dlll4.center) to [out=-90, in=180] (ddll4);
    \path[draw, line width=2pt, red] (dl4.center) to [out=-90, in=0] (ddll4);
    \path[draw, line width=2pt, blue] (dr4.center) to [out=-90, in=180] (ddrr4);
    \path[draw, line width=2pt, blue] (drrr4.center) to [out=-90, in=0] (ddrr4);
    \path[draw, line width=2pt, red] (ddll4) -- (dddll4.center);
    \path[draw, line width=2pt, blue] (ddrr4) -- (dddrr4.center);

    \node at (24, 0) [] {=};

    \filldraw[draw=black, fill=white] (24.5, -3) rectangle (31.5, 3);

    \node (uuulll5) at (25, 3) [] {};
    \node (uuul5) at (27, 3) [] {};
    \node (ul5) at (27, 1) [] {};
    \node (ur5) at (29, 1) [bullet, red] {};
    \node (urrr5) at (31, 1) [bullet, blue] {};
    \node (cc5) at (28, 0) [bullet] {};
    \node (dlll5) at (25, -1) [] {};
    \node (dl5) at (27, -1) [] {};
    \node (dr5) at (29, -1) [] {};
    \node (drrr5) at (31, -1) [] {};
    \node (ddll5) at (26, -2) [bullet, red] {};
    \node (ddrr5) at (30, -2) [bullet, blue] {};
    \node (dddll5) at (26, -3) [] {};
    \node (dddrr5) at (30, -3) [] {};

    \path[draw, line width=2pt, red] (uuulll5.center) -- (dlll5.center);
    \path[draw, line width=2pt, blue] (uuul5.center) -- (ul5.center);
    \path[draw, line width=2pt, blue] (ul5.center) to [out=-90, in=180] (cc5);
    \path[draw, line width=2pt, red] (ur5) to [out=-90, in=0] (cc5);
    \path[draw, line width=2pt, blue] (urrr5) -- (drrr5.center);
    \path[draw, line width=2pt, red] (cc5) to [out=180, in=90] (dl5.center);
    \path[draw, line width=2pt, blue] (cc5) to [out=0, in=90] (dr5.center);
    \path[draw, line width=2pt, red] (dlll5.center) to [out=-90, in=180] (ddll5);
    \path[draw, line width=2pt, red] (dl5.center) to [out=-90, in=0] (ddll5);
    \path[draw, line width=2pt, blue] (dr5.center) to [out=-90, in=180] (ddrr5);
    \path[draw, line width=2pt, blue] (drrr5.center) to [out=-90, in=0] (ddrr5);
    \path[draw, line width=2pt, red] (ddll5) -- (dddll5.center);
    \path[draw, line width=2pt, blue] (ddrr5) -- (dddrr5.center);
\end{tikzpicture}\end{center}
The first and last equalities follow from the axioms of a distributive law, the middle ones are unitality for the two monads. The next string diagram calculation proves associativity
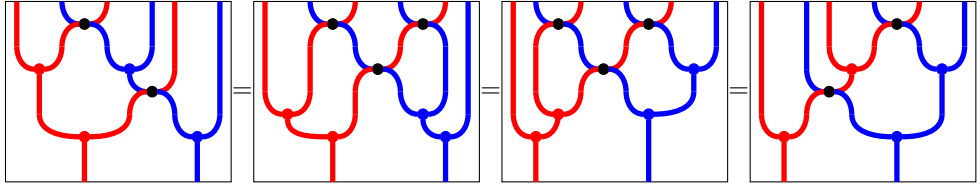
\begin{center}\begin{tikzpicture}[scale=0.3]
    \filldraw[draw=black, fill=white] (-0.5, -4) rectangle (9.5, 4);

    \node (uuuulllll1) at (0, 4) [] {};
    \node (uuuulll1) at (2, 4) [] {};
    \node (uuuul1) at (4, 4) [] {};
    \node (uuuurr1) at (6, 4) [] {};
    \node (uuuurrr1) at (7, 4) [] {};
    \node (uuuurrrrr1) at (9, 4) [] {};
    \node (uuull1) at (3, 3) [bullet] {};
    \node (uulllll1) at (0, 2) [] {};
    \node (uulll1) at (2, 2) [] {};
    \node (uul1) at (4, 2) [] {};
    \node (uurr1) at (6, 2) [] {};
    \node (ullll1) at (1, 1) [bullet, red] {};
    \node (ur1) at (5, 1) [bullet, blue] {};
    \node (urrr1) at (7, 1) [] {};
    \node (crr1) at (6, 0) [bullet] {};
    \node (dllll1) at (1, -1) [] {};
    \node (dr1) at (5, -1) [] {};
    \node (drrr1) at (7, -1) [] {};
    \node (drrrrr1) at (9, -1) [] {};
    \node (ddll1) at (3, -2) [bullet, red] {};
    \node (ddrrrr1) at (8, -2) [bullet, blue] {};
    \node (ddddll1) at (3, -4) [] {};
    \node (ddddrrrr1) at (8, -4) [] {};

    \path[draw, line width=2pt, red] (uuuulllll1.center) -- (uulllll1.center);
    \path[draw, line width=2pt, blue] (uuuulll1.center) to [out=-90, in=180] (uuull1);
    \path[draw, line width=2pt, red] (uuuul1.center) to [out=-90, in=0] (uuull1);
    \path[draw, line width=2pt, blue] (uuuurr1.center) -- (uurr1.center);
    \path[draw, line width=2pt, red] (uuuurrr1.center) -- (urrr1.center);
    \path[draw, line width=2pt, blue] (uuuurrrrr1.center) -- (drrrrr1.center);
    \path[draw, line width=2pt, red] (uuull1) to [out=180, in=90] (uulll1.center);
    \path[draw, line width=2pt, blue] (uuull1) to [out=0, in=90] (uul1.center);
    \path[draw, line width=2pt, red] (uulllll1.center) to [out=-90, in=180] (ullll1);
    \path[draw, line width=2pt, red] (uulll1.center) to [out=-90, in=0] (ullll1);
    \path[draw, line width=2pt, blue] (uul1.center) to [out=-90, in=180] (ur1);
    \path[draw, line width=2pt, blue] (uurr1.center) to [out=-90, in=0] (ur1);
    \path[draw, line width=2pt, red] (ullll1) -- (dllll1.center);
    \path[draw, line width=2pt, blue] (ur1) to [out=-90, in=180] (crr1);
    \path[draw, line width=2pt, red] (urrr1.center) to [out=-90, in=0] (crr1);
    \path[draw, line width=2pt, red] (crr1) to [out=180, in=90] (dr1.center);
    \path[draw, line width=2pt, blue] (crr1) to [out=0, in=90] (drrr1.center);
    \path[draw, line width=2pt, red] (dllll1.center) to [out=-90, in=180] (ddll1);
    \path[draw, line width=2pt, red] (dr1.center) to [out=-90, in=0] (ddll1);
    \path[draw, line width=2pt, blue] (drrr1.center) to [out=-90, in=180] (ddrrrr1);
    \path[draw, line width=2pt, blue] (drrrrr1.center) to [out=-90, in=0] (ddrrrr1);
    \path[draw, line width=2pt, red] (ddll1) -- (ddddll1.center);
    \path[draw, line width=2pt, blue] (ddrrrr1) -- (ddddrrrr1.center);

    \node at (10, 0) [] {=};

    \filldraw[draw=black, fill=white] (10.5, -4) rectangle (20.5, 4);

    \node (uuuulllll2) at (11, 4) [] {};
    \node (uuuulll2) at (13, 4) [] {};
    \node (uuuul2) at (15, 4) [] {};
    \node (uuuurr2) at (17, 4) [] {};
    \node (uuuurrrr2) at (19, 4) [] {};
    \node (uuuurrrrr2) at (20, 4) [] {};
    \node (uuull2) at (14, 3) [bullet] {};
    \node (uuurrr2) at (18, 3) [bullet] {};
    \node (uulll2) at (13, 2) [] {};
    \node (uul2) at (15, 2) [] {};
    \node (uurr2) at (17, 2) [] {};
    \node (uurrrr2) at (19, 2) [] {};
    \node (ur2) at (16, 1) [bullet] {};
    \node (clllll2) at (11, 0) [] {};
    \node (clll2) at (13, 0) [] {};
    \node (cl2) at (15, 0) [] {};
    \node (crr2) at (17, 0) [] {};
    \node (crrrr2) at (19, 0) [] {};
    \node (dllll2) at (12, -1) [bullet, red] {};
    \node (dl2) at (15, -1) [] {};
    \node (drrr2) at (18, -1) [bullet, blue] {};
    \node (drrrrr2) at (20, -1) [] {};
    \node (ddll2) at (14, -2) [bullet, red] {};
    \node (ddrrrr2) at (19, -2) [bullet, blue] {};
    \node (ddddll2) at (14, -4) [] {};
    \node (ddddrrrr2) at (19, -4) [] {};

    \path[draw, line width=2pt, red] (uuuulllll2.center) -- (clllll2.center);
    \path[draw, line width=2pt, blue] (uuuulll2.center) to [out=-90, in=180] (uuull2);
    \path[draw, line width=2pt, red] (uuuul2.center) to [out=-90, in=0] (uuull2);
    \path[draw, line width=2pt, blue] (uuuurr2.center) to [out=-90, in=180] (uuurrr2);
    \path[draw, line width=2pt, red] (uuuurrrr2.center) to [out=-90, in=0] (uuurrr2);
    \path[draw, line width=2pt, blue] (uuuurrrrr2.center) -- (drrrrr2.center);
    \path[draw, line width=2pt, red] (uuull2) to [out=180, in=90] (uulll2.center);
    \path[draw, line width=2pt, blue] (uuull2) to [out=0, in=90] (uul2.center);
    \path[draw, line width=2pt, red] (uuurrr2) to [out=180, in=90] (uurr2.center);
    \path[draw, line width=2pt, blue] (uuurrr2) to [out=0, in=90] (uurrrr2.center);
    \path[draw, line width=2pt, red] (uulll2.center) -- (clll2.center);
    \path[draw, line width=2pt, blue] (uul2.center) to [out=-90, in=180] (ur2);
    \path[draw, line width=2pt, red] (uurr2.center) to [out=-90, in=0] (ur2);
    \path[draw, line width=2pt, blue] (uurrrr2.center) -- (crrrr2.center);
    \path[draw, line width=2pt, red] (ur2) to [out=180, in=90] (cl2.center);
    \path[draw, line width=2pt, blue] (ur2) to [out=0, in=90] (crr2.center);
    \path[draw, line width=2pt, red] (clllll2.center) to [out=-90, in=180] (dllll2);
    \path[draw, line width=2pt, red] (clll2.center) to [out=-90, in=0] (dllll2);
    \path[draw, line width=2pt, red] (cl2.center) -- (dl2.center);
    \path[draw, line width=2pt, blue] (crr2.center) to [out=-90, in=180] (drrr2);
    \path[draw, line width=2pt, blue] (crrrr2.center) to [out=-90, in=0] (drrr2);
    \path[draw, line width=2pt, red] (dllll2) to [out=-90, in=180] (ddll2);
    \path[draw, line width=2pt, red] (dl2.center) to [out=-90, in=0] (ddll2);
    \path[draw, line width=2pt, blue] (drrr2) to [out=-90, in=180] (ddrrrr2);
    \path[draw, line width=2pt, blue] (drrrrr2.center) to [out=-90, in=0] (ddrrrr2);
    \path[draw, line width=2pt, red] (ddll2) -- (ddddll2.center);
    \path[draw, line width=2pt, blue] (ddrrrr2) -- (ddddrrrr2.center);

    \node at (21, 0) [] {=};

    \filldraw[draw=black, fill=white] (21.5, -4) rectangle (31.5, 4);

    \node (uuuulllll3) at (22, 4) [] {};
    \node (uuuullll3) at (23, 4) [] {};
    \node (uuuull3) at (25, 4) [] {};
    \node (uuuur3) at (27, 4) [] {};
    \node (uuuurrr3) at (29, 4) [] {};
    \node (uuuurrrrr3) at (31, 4) [] {};
    \node (uuulll3) at (24, 3) [bullet] {};
    \node (uuurr3) at (28, 3) [bullet] {};
    \node (uullll3) at (23, 2) [] {};
    \node (uull3) at (25, 2) [] {};
    \node (uur3) at (27, 2) [] {};
    \node (uurrr3) at (29, 2) [] {};
    \node (uurrrrr3) at (31, 2) [] {};
    \node (ul3) at (26, 1) [bullet] {};
    \node (urrrr3) at (30, 1) [bullet, blue] {};
    \node (cllll3) at (23, 0) [] {};
    \node (cll3) at (25, 0) [] {};
    \node (cr3) at (27, 0) [] {};
    \node (crrrr3) at (30, 0) [] {};
    \node (dlllll3) at (22, -1) [] {};
    \node (dlll3) at (24, -1) [bullet, red] {};
    \node (drr3) at (28, -1) [bullet, blue] {};
    \node (ddllll3) at (23, -2) [bullet, red] {};
    \node (ddddllll3) at (23, -4) [] {};
    \node (ddddrr3) at (28, -4) [] {};

    \path[draw, line width=2pt, red] (uuuulllll3.center) -- (dlllll3.center);
    \path[draw, line width=2pt, blue] (uuuullll3.center) to [out=-90, in=180] (uuulll3);
    \path[draw, line width=2pt, red] (uuuull3.center) to [out=-90, in=0] (uuulll3);
    \path[draw, line width=2pt, blue] (uuuur3.center) to [out=-90, in=180] (uuurr3);
    \path[draw, line width=2pt, red] (uuuurrr3.center) to [out=-90, in=0] (uuurr3);
    \path[draw, line width=2pt, blue] (uuuurrrrr3.center) -- (uurrrrr3.center);
    \path[draw, line width=2pt, red] (uuulll3) to [out=180, in=90] (uullll3.center);
    \path[draw, line width=2pt, blue] (uuulll3) to [out=0, in=90] (uull3.center);
    \path[draw, line width=2pt, red] (uuurr3) to [out=180, in=90] (uur3.center);
    \path[draw, line width=2pt, blue] (uuurr3) to [out=0, in=90] (uurrr3.center);
    \path[draw, line width=2pt, red] (uullll3.center) -- (cllll3.center);
    \path[draw, line width=2pt, blue] (uull3.center) to [out=-90, in=180] (ul3);
    \path[draw, line width=2pt, red] (uur3.center) to [out=-90, in=0] (ul3);
    \path[draw, line width=2pt, blue] (uurrr3.center) to [out=-90, in=180] (urrrr3);
    \path[draw, line width=2pt, blue] (uurrrrr3.center) to [out=-90, in=0] (urrrr3);
    \path[draw, line width=2pt, red] (ul3) to [out=180, in=90] (cll3.center);
    \path[draw, line width=2pt, blue] (ul3) to [out=0, in=90] (cr3.center);
    \path[draw, line width=2pt, blue] (urrrr3) -- (crrrr3.center);
    \path[draw, line width=2pt, red] (cllll3.center) to [out=-90, in=180] (dlll3);
    \path[draw, line width=2pt, red] (cll3.center) to [out=-90, in=0] (dlll3);
    \path[draw, line width=2pt, blue] (cr3.center) to [out=-90, in=180] (drr3);
    \path[draw, line width=2pt, blue] (crrrr3.center) to [out=-90, in=0] (drr3);
    \path[draw, line width=2pt, red] (dlllll3.center) to [out=-90, in=180] (ddllll3);
    \path[draw, line width=2pt, red] (dlll3.center) to [out=-90, in=0] (ddllll3);
    \path[draw, line width=2pt, blue] (drr3) -- (ddddrr3.center);
    \path[draw, line width=2pt, red] (ddllll3) -- (ddddllll3.center);

    \node at (32, 0) [] {=};

    \filldraw[draw=black, fill=white] (32.5, -4) rectangle (42.5, 4);

    \node (uuuulllll4) at (33, 4) [] {};
    \node (uuuulll4) at (35, 4) [] {};
    \node (uuuull4) at (36, 4) [] {};
    \node (uuuur4) at (38, 4) [] {};
    \node (uuuurrr4) at (40, 4) [] {};
    \node (uuuurrrrr4) at (42, 4) [] {};
    \node (uuurr4) at (39, 3) [bullet] {};
    \node (uull4) at (36, 2) [] {};
    \node (uur4) at (38, 2) [] {};
    \node (uurrr4) at (40, 2) [] {};
    \node (uurrrrr4) at (42, 2) [] {};
    \node (ulll4) at (35, 1) [] {};
    \node (ul4) at (37, 1) [bullet, red] {};
    \node (urrrr4) at (41, 1) [bullet, blue] {};
    \node (cll4) at (36, 0) [bullet] {};
    \node (dlllll4) at (33, -1) [] {};
    \node (dlll4) at (35, -1) [] {};
    \node (dl4) at (37, -1) [] {};
    \node (drrrr4) at (41, -1) [] {};
    \node (ddllll4) at (34, -2) [bullet, red] {};
    \node (ddrr4) at (39, -2) [bullet, blue] {};
    \node (ddddllll4) at (34, -4) [] {};
    \node (ddddrr4) at (39, -4) [] {};

    \path[draw, line width=2pt, red] (uuuulllll4.center) -- (dlllll4.center);
    \path[draw, line width=2pt, blue] (uuuulll4.center) -- (ulll4.center);
    \path[draw, line width=2pt, red] (uuuull4.center) -- (uull4.center);
    \path[draw, line width=2pt, blue] (uuuur4.center) to [out=-90, in=180] (uuurr4);
    \path[draw, line width=2pt, red] (uuuurrr4.center) to [out=-90, in=0] (uuurr4);
    \path[draw, line width=2pt, blue] (uuuurrrrr4.center) -- (uurrrrr4.center);
    \path[draw, line width=2pt, red] (uuurr4) to [out=180, in=90] (uur4.center);
    \path[draw, line width=2pt, blue] (uuurr4) to [out=0, in=90] (uurrr4.center);
    \path[draw, line width=2pt, red] (uull4.center) to [out=-90, in=180] (ul4);
    \path[draw, line width=2pt, red] (uur4.center) to [out=-90, in=0] (ul4);
    \path[draw, line width=2pt, blue] (uurrr4.center) to [out=-90, in=180] (urrrr4);
    \path[draw, line width=2pt, blue] (uurrrrr4.center) to [out=-90, in=0] (urrrr4);
    \path[draw, line width=2pt, blue] (ulll4.center) to [out=-90, in=180] (cll4);
    \path[draw, line width=2pt, red] (ul4) to [out=-90, in=0] (cll4);
    \path[draw, line width=2pt, blue] (urrrr4) -- (drrrr4.center);
    \path[draw, line width=2pt, red] (cll4) to [out=180, in=90] (dlll4.center);
    \path[draw, line width=2pt, blue] (cll4) to [out=0, in=90] (dl4.center);
    \path[draw, line width=2pt, red] (dlllll4.center) to [out=-90, in=180] (ddllll4);
    \path[draw, line width=2pt, red] (dlll4.center) to [out=-90, in=0] (ddllll4);
    \path[draw, line width=2pt, blue] (dl4.center) to [out=-90, in=180] (ddrr4);
    \path[draw, line width=2pt, blue] (drrrr4.center) to [out=-90, in=0] (ddrr4);
    \path[draw, line width=2pt, red] (ddllll4) -- (ddddllll4.center);
    \path[draw, line width=2pt, blue] (ddrr4) -- (ddddrr4.center);    
\end{tikzpicture}\end{center}
The first and last equalities follow from the axioms of a distributive law, the middle one is associativity for the two monads.

\end{proof}
\end{prop}

As it turns out, distributive laws can be encoded as monads too

\begin{prop}\thmheader{Distributive laws as monads}

The data of a distributive law $\gamma : (X, t, \eta^t, \mu^t) \rightsquigarrow (X, s, \eta^S, \mu^s)$ between two monads $\mathcal{T} = (X, t, \eta^t, \mu^t), \mathcal{S} = (X, s, \eta^s, \mu^s)$ over the same 0-cell $X : \mathcal{B}$ is equivalently encoded as a monad
\[\Gamma = ((X, t, \eta^t, \mu^t), (s, \gamma), \eta^s, \mu^s) : \mathbf{Mnd}(\mathbf{Mnd}(\mathcal{B}))\]
in the 2-category $\mathbf{Mnd}(\mathcal{B}) \simeq \mathbf{Lax}[\mathbb{1}, \mathcal{B}]$

\begin{proof}
The proof consists in just writing down the axioms for both and matching them:

\begin{enumerate}
\item the fact that $(X, s, \eta^S, \mu^S)$ is a monad follows from the fact that $(s, \gamma)$ is a monad on $(X, t, \eta^t, \mu^t)$ with unit $\eta^s$ and multiplication $\mu^s$;
\item the fact that $\gamma$ is a distributive law follows from the assumptions that $(s, \gamma)$ is a 1-cell and $\eta^s$, $\mu^s$ are 2-cells in $\mathbf{Mnd}(\mathcal{B})$.
\end{enumerate}

\todo{Write down explicitly? Maybe not necessary...}
\end{proof}
\end{prop}

This suggests we \textit{define} the 2-category of distributive laws as ``monads in monads''

\begin{defin}\thmheader{2-category of distributive laws}

We define the 2-category $\mathbf{Dist}(\mathcal{B}) := \mathbf{Mnd}(\mathbf{Mnd}(\mathcal{B}))$.
\end{defin}

As we did for the 2-category $\mathbf{Mnd}(\mathcal{B})$, it is worth giving a concrete description of $\mathbf{Dist}(\mathcal{B})$

\begin{rmk}\thmheader{2-category of distributive laws, explicitly}

Fix a (strict, for the sake of convenience) 2-category $\mathcal{B}$. The 2-category $\mathbf{Dist}(\mathcal{B})$ can be presented\footnote{Meaning $\mathbf{Dist}(\mathcal{B})$ it is strictly equivalent to the 2-category being defined: there's strict 2-functors back-and-forth (defined in the obvious way), which compose to the identity 2-functor on the nose.} in the following way:
\begin{enumerate}
\item its 0-cells are 8-tuples $(X, t, s, \eta^t, \eta^s, \mu^t, \mu^s, \sigma)$, with
\begin{itemize}
\item $X : \mathcal{B}$
\item $t, s : X \to X : \mathcal{B}$
\item $\eta^t : \mathbf{Id}_X \to t : X \to X : \mathcal{B}$
\item $\eta^s : \mathbf{Id}_X \to s : X \to X : \mathcal{B}$
\item $\mu^t : t \cdot t \to t : X \to X : \mathcal{B}$
\item $\mu^s : s \cdot s \to s : X \to X : \mathcal{B}$
\item $\sigma : s \cdot t \to t \cdot s : X \to X : \mathcal{B}$
\end{itemize}
satisfying the following axioms:
\begin{enumerate}
\item $(X, t, \eta^t, \mu^t) : \mathbf{Mnd}(\mathcal{B})$ is a monad:
\begin{center}\begin{tikzcd}
X \cdot t \arrow[r, equals] \arrow[d, "\eta^t \cdot t"] & t \arrow[d, equals] & t \cdot X \arrow[l, equals] \arrow[d, "t \cdot \eta^t"] & t \cdot t \cdot t \arrow[r, "t \cdot \mu^t"] \arrow[d, "\mu^t \cdot t"] & t \cdot t \arrow[d, "\mu^t"] \\
t \cdot t \arrow[r, "\mu^t"] & t & t \cdot t \arrow[l, "\mu^t"'] & t \cdot t \arrow[r, "\mu^t"] & t
\end{tikzcd}\end{center}
\item $(X, s, \eta^s, \mu^s) : \mathbf{Mnd}(\mathcal{B})$ is a monad:
\begin{center}\begin{tikzcd}
X \cdot s \arrow[r, equals] \arrow[d, "\eta^s \cdot s"] & s \arrow[d, equals] & s \cdot X \arrow[l, equals] \arrow[d, "s \cdot \eta^s"] & s \cdot s \cdot s \arrow[r, "s \cdot \mu^s"] \arrow[d, "\mu^s \cdot s"] & s \cdot s \arrow[d, "\mu^s"] \\
s \cdot s \arrow[r, "\mu^s"] & s & s \cdot s \arrow[l, "\mu^s"'] & s \cdot s \arrow[r, "\mu^s"] & s
\end{tikzcd}\end{center}
\item $\sigma : (X, t, \eta^t, \mu^t) \rightsquigarrow (X, s, \eta^s, \mu^s)$ is a distributive law:
\begin{center}\begin{tikzcd}
X \cdot t \arrow[d, equals] \arrow[r, "\eta^s \cdot t"] & s \cdot t \arrow[dd, "\sigma"] & s \cdot X \arrow[l, "s \cdot \eta^t"'] \arrow[d, equals] & s \cdot s \cdot t \arrow[r, "\mu^s \cdot t"] \arrow[d, "s \cdot \sigma"] & s \cdot t \arrow[dd, "\sigma"] & s \cdot t \cdot t \arrow[l, "s \cdot \mu^t"'] \arrow[d, "\sigma \cdot t"] \\
t \arrow[d, equals] && s \arrow[d, equals] & s \cdot t \cdot s \arrow[d, "\sigma \cdot s"] && t \cdot s \cdot t \arrow[d, "t \cdot \sigma"] \\
t \cdot X \arrow[r, "t \cdot \eta^s"] & t \cdot s & X \cdot s \arrow[l, "\eta^t \cdot s"'] & t \cdot s \cdot s \arrow[r, "t \cdot \mu^s"] & t \cdot s & t \cdot t \cdot s \arrow[l, "\mu^t \cdot s"']
\end{tikzcd}\end{center}
\end{enumerate}
\item its 1-cells are triples $(f, \phi, \psi) : (X, t, s, \eta^t, \eta^s, \mu^t, \mu^s, \sigma) \to (Y, r, u, \eta^r, \eta^u, \mu^r, \mu^u, \rho)$ where
\begin{itemize}
\item $f : X \to Y : \mathcal{B}$
\item $\phi : f \cdot t \to r \cdot f$
\item $\psi : f \cdot s \to u \cdot f$
\end{itemize}
satisfying the following axioms:
\begin{enumerate}
\item $(f, \phi) : (X, t, \eta^t, \mu^t) \to (Y, r, \eta^r, \mu^r) : \mathbf{Mnd}(\mathcal{B})$
\begin{center}\begin{tikzcd}
f \cdot X \arrow[d, equals] \arrow[r, "f \cdot \eta^t"] & f \cdot t \arrow[dd, "\phi"] & f \cdot t \cdot t \arrow[l, "f \cdot \mu^t"'] \arrow[d, "\phi \cdot t"] \\
f \arrow[d, equals] && r \cdot f \cdot t \arrow[d, "r \cdot \phi"] \\
Y \cdot f \arrow[r, "\eta^r \cdot f"] & r \cdot f & r \cdot r \cdot f \arrow[l, "\mu^r \cdot f"']
\end{tikzcd}\end{center}
\item $(f, \psi) : (X, s, \eta^s, \mu^s) \to (Y, u, \eta^u, \mu^u) : \mathbf{Mnd}(\mathcal{B})$
\begin{center}\begin{tikzcd}
f \cdot X \arrow[d, equals] \arrow[r, "f \cdot \eta^s"] & f \cdot s \arrow[dd, "\psi"] & f \cdot s \cdot s \arrow[l, "f \cdot \mu^s"'] \arrow[d, "\psi \cdot s"] \\
f \arrow[d, equals] && u \cdot f \cdot s \arrow[d, "u \cdot \psi"] \\
Y \cdot f \arrow[r, "\eta^u \cdot f"] & u \cdot f & u \cdot u \cdot f \arrow[l, "\mu^u \cdot f"']
\end{tikzcd}\end{center}
\item $\psi : (f, \phi) \cdot (s, \sigma) \to (u, \rho) \cdot(f, \phi)$ is a 2-cell in $\mathbf{Mnd}(\mathcal{C})$
\begin{center}\begin{tikzcd}
f \cdot s \cdot t \arrow[r, "f \cdot \sigma"] \arrow[d, "\psi \cdot t"] & f \cdot t \cdot s \arrow[d, "\phi \cdot s"] \\
u \cdot f \cdot t \arrow[d, "u \cdot \phi"] & r \cdot f \cdot s \arrow[d, "r \cdot \psi"] \\
u \cdot r \cdot f \arrow[r, "\rho \cdot f"] & r \cdot u \cdot f
\end{tikzcd}\end{center}
\end{enumerate}
\item given composable 1-cells $(f_1, \phi_1, \psi_1)$, $(f_2, \phi_2, \psi_2)$ their composition is given by
\[(f_2, \phi_2, \psi_2) \cdot (f_1, \phi_1, \psi_1) := (f_2 \cdot f_1, (\phi_2 \cdot f_1) \circ (f_2 \cdot \phi_1), (\psi_2 \cdot f_1) \circ (f_2 \cdot \psi_1))\]
\item its 2-cells $\gamma : (f, \phi_1, \psi_1) \to (g, \phi_2, \psi_2)$ are 2-cells $\gamma : f \to g : X \to Y : \mathcal{B}$ satisfying the following axioms:
\begin{enumerate}
\item $\gamma : (f, \phi_1) \to (g, \phi_2) : (X, t, \eta^t, \mu^t) \to (Y, r, \eta^r, \mu^r) : \mathbf{Mnd}(\mathcal{B})$
\begin{center}\begin{tikzcd}
f \cdot t \arrow[r, "\phi_1"] \arrow[d, "\gamma \cdot t"] & r \cdot f \arrow[d, "r \cdot \gamma"] \\
g \cdot t \arrow[r, "\phi_2"] & r \cdot g
\end{tikzcd}\end{center}
\item $\gamma : (f, \psi_1) \to (g, \psi_2) : (X, s, \eta^s, \mu^s) \to (Y, u, \eta^u, \mu^u) : \mathbf{Mnd}(\mathcal{B})$
\begin{center}\begin{tikzcd}
f \cdot s \arrow[r, "\psi_1"] \arrow[d, "\gamma \cdot s"] & u \cdot f \arrow[d, "u \cdot \gamma"] \\
g \cdot s \arrow[r, "\psi_2"] & u \cdot g
\end{tikzcd}\end{center}
\end{enumerate}
\item given vertically composable 2-cells $\gamma_1$, $\gamma_2$ their vertical composite is given by $\gamma_2 \circ \gamma_1$
\item given horizontally composable 2-cells $\gamma_1$, $\gamma_2$ their horizontal composite is given by $\gamma_2 \cdot \gamma_1$
\end{enumerate}
\begin{proof}
Immediate. \todo{prove?}
\end{proof}
\end{rmk}

We can therefore define 2-functors $\mathbf{Dist}(\mathcal{B}) \to \mathbf{Mnd}(\mathcal{B})$, encoding the recovery of the underlying monads (as well as the composite we described in \ref{Distributive laws allow for monad composition}):

\begin{rmk}\thmheader{Monads from distributive laws}

Given a (strict, for the sake of conciseness) 2-category $\mathcal{B}$, we have strict 2-functors
\[\mathcal{U}_1, \mathcal{U}_2, \mathcal{C} : \mathbf{Dist}(\mathcal{B}) \to \mathbf{Mnd}(\mathcal{B})\]
defined as follows:
\[\mathcal{U}_1((a, t, \eta^t, \mu^t), (s, \sigma), \eta^s, \mu^s) = (a, t, \eta^t, \mu^t)\]
\[\mathcal{U}_2((a, t, \eta^t, \mu^t), (s, \sigma), \eta^s, \mu^s) = (a, s, \eta^s, \mu^s)\]
\[\mathcal{C}((a, t, \eta^t, \mu^t), (s, \sigma), \eta^s, \mu^s) = (a, s \cdot t, \eta^s \cdot \eta^t, (s \cdot \sigma \cdot t) \circ (\mu^s \cdot \mu^t))\]
\begin{proof}

We need to complete the definition of $\mathcal{U}_1, \mathcal{U}_2, \mathcal{C}$ to 2-functors: give the actions on 1- and 2-cells, and prove that makes them into 2-functors. Fix a 1-cell in $\mathbf{Dist}(\mathcal{B})$
\[((f, \phi), \psi) : ((a_1, t_1, \eta^{t}_{1}, \mu^{t}_{1}), (s_1, \sigma_1), \eta^{s}_{1}, \mu^{s}_{1}) \to ((a_2, t_2, \eta^{t}_{2}, \mu^{t}_{2}), (s_2, \sigma_2), \eta^{s}_{2}, \mu^{s}_{2})\]
We then define $\mathcal{U}_1((f, \phi), \psi) := (f, \phi)$, $\mathcal{U}_2((f, \phi), \psi) = (f, \psi)$ and $\mathcal{C}((f, \phi), \psi) = (f, (s_2 \cdot \phi) \circ (\psi \cdot t_1))$. Similarly, fix a 2-cell $\gamma : ((f_1, \phi_1), \psi_1) \to ((f_2, \phi_2), \psi_2)$: we (unsurprisingly) define $\mathcal{U}_1(\gamma) = \mathcal{U}_2(\gamma) = \mathcal{C}(\gamma) = \gamma$. Moreover, the images of 1-cells (resp. 2-cells) are 1-cells (2-cells) in $\mathbf{Mnd}(\mathcal{B})$: the argument for $\mathcal{U}_i$ follows immediately from our previous remark \ref{2-category of distributive laws, explicitly} (2.a and 4.a for $\mathcal{U}_1$, 2.b and 4.b for $\mathcal{U}_2$ and 2.c), while that for $\mathcal{C}$ follows from 2.a - 2.c and 4.a - 4.b.

That $\mathcal{U}_1$ and $\mathcal{U}_2$ are strict 2-functors is an immediate consequence of our previous remark, \ref{2-category of distributive laws, explicitly}. We will sketch the proof for $\mathcal{U}_1$, as the one for $\mathcal{U}_2$ is identical. Fix composable 1-cells
\[T_1 \xrightarrow{F_1 = ((f_1, \phi_1), \psi_1)} T_2 \xrightarrow{F_2 = ((f_2, \phi_2), \psi_2)} T_3\]
We want to show $\mathcal{U}_1(F_2 \cdot F_1) = \mathcal{U}_1(F_2) \cdot \mathcal{U}_1(F_2)$. Recall that
\[F_2 \cdot F_1 = ((f_2 \cdot f_1, (\phi_2 \cdot f_1) \circ (f_2 \cdot \phi_1)), (\psi_2 \cdot f_1) \circ (f_2 \cdot \psi_1))\]
so that $\mathcal{U}_1(F_2 \cdot F_1) = (f_2 \cdot f_1, (\phi_2 \cdot f_1) \circ (f_2 \cdot \phi_1))$. On the other hand, we have $\mathcal{U}_1(F_1) = (f_1, \phi_1)$ and $\mathcal{U}_1(F_2) = (f_2, \phi_2)$, implying
\[\mathcal{U}_1(F_2) \cdot \mathcal{U}_1(F_1) = (f_2, \phi_2) \cdot (f_1, \phi_1) = (f_2 \cdot f_1, (\phi_2 \cdot f_1) \circ (f_2 \cdot \phi_1))\]
Similarly, we can see that it sends identity 1-cells to identity 1-cells. Likewise, it is trivial to check that it preserves identity and (vertical, horizontal) composites of 2-cells: $\mathcal{U}_1$ is a strict 2-functor.

The case for $\mathcal{C}$ is slightly less immediate: consider composable 1-cells $F_1, F_2$ as before. We can compute
\begin{align*}
\mathcal{C}(F_2) \cdot \mathcal{C}(F_1) &= (f_2, (s_3 \cdot \phi_2) \circ (\psi_2 \cdot t_2)) \cdot (f_1, (s_2 \cdot \phi_1) \circ (\psi_1 \cdot t_1)) \\
&= (f_2 \cdot f_1, (s_2 \cdot \phi_2 \cdot f_1) \circ (\psi_2 \cdot t_2 \cdot f_1) \circ (f_2 \cdot s_2 \cdot \phi_1) \circ (f_2 \cdot \psi_1 \cdot t_1))
\end{align*}
\begin{align*}
\mathcal{C}(F_2 \cdot F_1) &= \mathcal{C}((f_2 \cdot f_1, (\phi_2 \cdot f_1) \circ (f_2 \cdot \phi_1)), (\psi_2 \cdot f_1) \circ (f_2 \cdot \psi_1)) \\
&= (f_2 \cdot f_1, (s_2 \cdot \phi_2 \cdot f_1) \circ (s_2 \cdot f_2 \cdot \phi_1) \circ (\psi_2 \cdot f_1 \cdot t_1) \circ (f_2 \cdot \psi_1 \cdot t_1)
\end{align*}
The first components are clearly equal; the second components are equal because of interchange. Showing that identity 1-cells are also preserved on the nose is immediate, as is showing identity and (vertical, horizontal) composites of 2-cells are preserved. Hence $\mathcal{C}$ is a strict 2-functor.
\end{proof}
\end{rmk}

This allows us to define what a parametric distributive law is:

\begin{defin}\thmheader{Parametric Distributive Law}

Given 2-categories $\mathcal{B}$ and $\mathcal{C}$, the category of $\mathcal{C}$-parametric distributive laws in $\mathcal{B}$ is the 2-functor 2-category
\[\mathbf{PDist}(\mathcal{C}, \mathcal{B}) := \mathbf{St}_{\mathbf{St}}[\mathcal{C}, \mathbf{Dist}(\mathcal{B})]\]
\end{defin}

An immediate consequence of \ref{Distributive laws as monads} is then
\[\mathbf{PDist}(\mathcal{C}, \mathcal{B}) = \mathbf{St}_{\mathbf{St}}[\mathcal{C}, \mathbf{Lax}[\mathbb{1}, \mathbf{Lax}[\mathbb{1}, \mathcal{B}]]]\]

The absence of a well-behaved left 2-adjoint to the $\mathbf{Lax}$ 2-functor makes this somewhat annoying to deal with; the main purpose of the next section is to alleviate this issue.

%% ## SECTION 2
\section{The lax functor classifier}

In this section we will present a construction that is known to experts in the field\footnote{E.g. as described in T. Johnson-Freyd's post \cite{delaxing-object}}, but for which we are not aware of any reference: the ``lax functor classifier''. Given a 2-category $\mathcal{B}$, its lax functor classifier $\widehat{\mathcal{B}}$ is a 2-category such that, for any other 2-category $\mathcal{C}$, we have an equivalence
\[\mathbf{Lax}[\mathcal{B}, \mathcal{C}] \simeq \mathbf{Ps}_{\mathbf{Lax}}[\widehat{\mathcal{B}}, \mathcal{C}]\]

It is reminiscent of a related construction, sometimes called the  ``lax morphism classifier'', introduced by Blackwell, Kelly and Power in \cite{2-dimensional-monad-theory}. In order to give a definition that is not too combinatorial in nature, we shall first recall the definition of the (augmented) simplex and interval categories, together with the duality that relates them.

\begin{rmk}\thmheader{Simplex and interval categories}

Recall that the category of preorders $\mathbf{PreOrd}$ is a full subcategory of the category of categories $\mathbf{Cat}$:
\[\mathbf{PreOrd} \hookrightarrow \mathbf{Cat}\]

The augmented simplex category, denoted by $\Delta_+ \hookrightarrow \mathbf{PreOrd}$ is the full subcategory of finite, totally ordered sets and the simplex category, denoted by $\Delta \hookrightarrow \Delta_+$ is the full subcategory of \textit{non-empty} such. We will employ the following, widespread conventions:
\begin{enumerate}
\item $\forall n \in \mathbb{N}$, write $[n-1] := \{0 \le \dots \le n-1\} : \Delta_+$ for the\footnote{Unique up to isomorphism in $\mathbf{PreOrd}$} poset with $n$ elements;
\item $\forall n \in \mathbb{N}$ and $0 \le k \le n$, write $\delta^{n}_{k} : [n-1] \to [n]$ for the unique injective, monotonic function such that $\delta^{n}_{k}(k) = k+1$;
\item $\forall n \in \mathbb{N}$ and $0 \le k \le n$, write $\sigma^{n}_{k} : [n+1] \to [n]$ for the unique surjective, monotonic function such that $\sigma^{n}_{k}(k) = k = \sigma^{n}_{k}(k+1)$.
\end{enumerate}

Furthermore the augmented interval category, denoted by $\mathbf{I}_+ \hookrightarrow \mathbf{PreOrd}$ is the subcategory of finite intervals (that is, totally ordered sets with \textit{not necessarily distinct} top and bottom element) and interval maps (monotonic functions preserving the top and bottom elements), and the interval category, denoted by $\mathbf{I} \hookrightarrow \mathbf{I}_+$ is for the full subcategory on objects with \textit{distinct} top and bottom elements.

Finally, recall that $\Delta^{op}_{+} \simeq \mathbf{I}_+$, the equivalence being induced by the hom-functors
\[\Delta_+[-, [1]] : \Delta^{op}_{+} \to \mathbf{I}_+\]
\[\mathbf{I}_+[-, [1]] : \mathbf{I}^{op}_{+} \to \Delta_+\]

Given $k \in \{-1, 0, \dots\}$, the poset structure on $\Delta_+[[k], [1]] = \{f : [k] \to [1] : \Delta_+\}$ is the usual, pointwise preorder (which happens to be a total order because the codomain is $[1]$), with top and bottom elements being the constant functions on $1$ and $0$ respectively. Analogously, given $k \in \mathbb{N}$, the poset structure on $\mathbf{I}_+[[k], [1]] = \{f : [k] \to [1] : \mathbf{I}_+\}$ is pointwise. Finally, recall that this equivalence restricts to an equivalence $\Delta^{op} \simeq \mathbf{I}$.
\end{rmk}

In what follows we will keep writing $[k]$ for both the (posetal) category and its corresponding locally discrete 2-category.

\begin{defin}\thmheader{The lax functor classifier}

The lax functor classifier $\widehat{\mathcal{B}}$ of a 2-category $\mathcal{B}$ is defined as follows:
\begin{enumerate}
\item Its 0-cells are the same as those of $\mathcal{B}$;
\item Its 1-cells $(k, \mathcal{F}) : a \to b : \widehat{\mathcal{B}}$ are pairs where $k \in \mathbb{N}$, $\mathcal{F} : [k] \to \mathcal{B}$ is a semistrict pseudofunctor (meaning it preserves identity 1-cells \textit{on the nose}), and $\mathcal{F}(0) = a, \; \mathcal{F}(k) = b$;
\item Identity 1-cells $\mathbb{1}_a := (0, !_a) : a \to a : \widehat{\mathcal{B}}$ are given by the strict 2-functors $!_a : [0] \to \mathcal{B}$ with $!_a(0) = a$ (and hence $!_a(\mathbb{1}_0) = \mathbb{1}_a$);
\item Its 2-cells $(\phi, \alpha) : (m, \mathcal{F}) \to (n, \mathcal{G}) : a \to b : \widehat{\mathcal{B}}$ are pairs where $\alpha : [n] \to [m] : \mathbf{I}_+$, $\mathcal{F} \circ \alpha$ and $\mathcal{G}$ agree on objects and $\phi : \mathcal{F} \circ \alpha \to \mathcal{G}$ is an ICON;
\item Identity 2-cells $\mathbb{1}_{(m, \mathcal{F})} : (m, \mathcal{F}) \to (m, \mathcal{F})$ are the 2-cells $(\mathbb{1}_{\mathcal{F}}, \mathbb{1}_{[m]})$
\item Given $(\phi, \alpha) : (m, \mathcal{G}) \to (n, \mathcal{H})$ and $(\psi, \beta) : (l, \mathcal{F}) \to (m, \mathcal{G})$, their vertical composite is defined by
\[(\phi, \alpha) \circ (\psi, \beta) := (\phi \circ (\psi \alpha), \beta \circ \alpha) : (l, \mathcal{F}) \to (n, \mathcal{H})\]
where $\psi \alpha : \mathcal{F} \circ \beta \circ \alpha \to \mathcal{G} \circ \alpha$ denotes the ICON obtained by whiskering $\psi$ on the right by the (strict) 2-functor $\alpha$.
\item Given $(n, \mathcal{G}) : b \to c$ and $(m, \mathcal{F}) : a \to b$, their (horizontal) composite is defined by
\[(n, \mathcal{G}) \cdot (m, \mathcal{F}) := (n + m, \mathcal{G} \star \mathcal{F})\]
where $\mathcal{G} \star \mathcal{F}$ is given by 2-categorical Day convolution with respect to (the dual of) the monoidal structure $(-\vee-)$ on $\mathbf{I}_+$ obtained by transferring ordinal sum over the equivalence $\Delta^{op}_{+} \simeq \mathbf{I}_+$. Explicitly, define on objects
\begin{align*}
(\mathcal{G} \star \mathcal{F})(k) := \begin{cases}
\mathcal{F}(k) & \text{if } k \le m \\
\mathcal{G}(k - m) & \text{if } m \le k
\end{cases}
\end{align*}
(which is well-defined because $\mathcal{F}(k) = b = \mathcal{G}(0)$) and on morphisms $f : i \to j : [m+n]$
\begin{align*}
(\mathcal{G} \star \mathcal{F}) := \begin{cases}
\mathcal{F}(f) & \text{if } j < m \\
\mathcal{G}(f) & \text{if } m < i \\
\mathcal{F}(f_1) \cdot \mathcal{G}(f_2) & \text{if } i \le m \le j \text{, hence } f : i \xrightarrow{f_2} m \xrightarrow{f_1} j
\end{cases}
\end{align*}
\item Given
\[(\psi, \beta) : (n_1, \mathcal{G}_1) \to (n_2, \mathcal{G}_2) : b \to c\]
\[(\phi, \alpha) : (m_1, \mathcal{F}_1) \to (m_2, \mathcal{F}_2) : a \to b\]
define their horizontal composite by $(\psi, \beta) \cdot (\phi, \alpha) := (\psi \star \phi, \beta \vee \alpha)$ where $(- \vee -)$ is the aforementioned monoidal structure on $\mathbf{I}_+$, and $\psi \star \phi$ is Day convolution of ICONs. Explicitly, define $\beta \vee \alpha$ as:
\begin{align*}
(\beta \vee \alpha)(k) := \begin{cases}
\alpha(k) & \text{if } k \le m_1 \\
\beta(k-m_1)+n_1 & \text{if } m_1 \le k
\end{cases}
\end{align*}
(which is well-defined since $\alpha(m_1) = n_1$ and $\beta(0) = 0$) and $\psi \star \phi$, for any $f : i \to j : [m_1 + m_2]$
\begin{align*}
(\psi \star \phi)_f := \begin{cases}
\phi_f & \text{if } j < m_1 \\
\psi_f & \text{if } m_1 < i \\
\phi_{f_1} \cdot\psi_{f_2} & \text{if } i \le m_1 \le j \text{, hence } f : i \xrightarrow{f_2} m_1 \xrightarrow{f_1} j
\end{cases}
\end{align*}
\item The left and right unitors are defined as follows
\[\lambda_{(k, \mathcal{F})} := (\lambda^{\mathcal{F}}, \mathbb{1}_{[k]}) : \mathbb{1}_b \cdot (k, \mathcal{F}) \to (k, \mathcal{F})\]
\[\rho_{(k, \mathcal{F})} := (\rho^{\mathcal{F}}, \mathbb{1}_{[k]}) : (k, \mathcal{F}) \cdot \mathbb{1}_a \to (k, \mathcal{F})\]
with $\lambda^{\mathcal{F}}, \rho^{\mathcal{F}}$ being the ICONs with, for all $f : i \to j : [k]$:
\begin{align*}
\lambda^{\mathcal{F}}_{f} := \begin{cases}
\mathbb{1}_{\mathcal{F}(f)} &\text{if } j < k \\
\lambda_{\mathcal{F}(f)} &\text{if } j = k
\end{cases} \;\;\;;\;\;\;
\rho^{\mathcal{F}}_{f} := \begin{cases}
\rho_{\mathcal{F}(f)} & \text{if } i = 0 \\
\mathbb{1}_{\mathcal{F}(f)} & \text{if } i > 0
\end{cases}
\end{align*}
\item The associator
\[\alpha_{(l, \mathcal{F}), (m, \mathcal{G}), (n, \mathcal{H})} : ((l, \mathcal{F}) \cdot (m, \mathcal{G})) \cdot (n, \mathcal{H}) \to (l, \mathcal{F}) \cdot ((m, \mathcal{G}) \cdot (n, \mathcal{H}))\]
is defined as
\[\alpha_{(l, \mathcal{F}), (m, \mathcal{G}), (n, \mathcal{H})} := (\alpha^{\mathcal{F}, \mathcal{G}, \mathcal{H}}, \mathbb{1}_{[n+m+l]})\]
where $\alpha^{\mathcal{F}, \mathcal{G}, \mathcal{H}}$ is the ICON with, for any $f : i \to j : [n + m + l]$,
\[\alpha^{\mathcal{F}, \mathcal{G}, \mathcal{H}}_{f} := \alpha_{\mathcal{F}(f_1), \mathcal{G}(f_2), \mathcal{H}(f_3)} \text{ whenever } i \le n \le n + m \le j\]
(after factoring $f$ appropriately), and the appropriate identity 2-cell in all other cases.
\end{enumerate}
\begin{proof} It is clear that, for any $a, b : \mathcal{B}$, the 1-cells and 2-cells defined above make up a category: checking that (vertical) composition is associative and unital is a simple exercise.

It is equally straightforward (though slightly less immediate) to see that horizontal composition of 1-cells satisfies the usual triangle and pentagon identities. We are only left to check that horizontal composition is functorial, that is, interchange: given the following
\[(m_1, \mathcal{F}_1) \xrightarrow{(\phi_1, \alpha_1)} (m_2, \mathcal{F}_2) \xrightarrow{(\phi_2, \alpha_2)} (m_3, \mathcal{F}_3) : b \to c\]
\[(n_1, \mathcal{G}_1) \xrightarrow{(\psi_1, \beta_1)} (n_2, \mathcal{G}_2) \xrightarrow{(\psi_2, \beta_2)} (n_3, \mathcal{G}_3) : a \to b\]
we need to check that
\[((\phi_2, \alpha_2) \circ (\phi_1, \alpha_1)) \cdot ((\psi_2, \beta_2) \circ (\psi_1, \beta_1)) = ((\phi_2, \alpha_2) \cdot (\psi_2, \beta_2)) \circ ((\phi_1, \alpha_1) \cdot (\psi_1, \beta_1))\]
It is easy to check that the left-hand side equals
\[((\phi_2 \circ (\phi_1 \alpha_2)) \star (\psi_2 \circ (\psi_1 \beta_2)), (\alpha_1 \circ \alpha_2) \vee (\beta_1 \circ \beta_2))\]
and the right-hand side equals
\[((\phi_2 \star \psi_2) \circ ((\phi_1 \star \psi_1) (\alpha_2 \vee \beta_2)), (\alpha_1 \vee \beta_1) \circ (\alpha_2 \vee \beta_2))\]

One can explicitly check that $(\alpha_1 \circ \alpha_2) \vee (\beta_1 \circ \beta_2) = (\alpha_1 \vee \beta_1) \circ (\alpha_2 \vee \beta_2)$, so the second components agree. For the first component, notice how
\begin{enumerate}
\item $(\phi_1 \star \psi_1)(\alpha_2 \vee \beta_2) = (\phi_1 \alpha_2) \star (\psi_1 \beta_2)$
\item $(\phi_2 \circ \theta_1) \star (\psi_2 \circ \theta_2) = (\phi_2 \star \psi_2) \circ (\theta_1 \star \theta_2)$
\end{enumerate}
These two facts imply that the first components agree, concluding the proof.
\end{proof}
\end{defin}

We can now prove the main result in this section:

\begin{thm}\thmheader{The lax functor classifier classifies lax functors}

Given 2-categories $\mathcal{B}, \mathcal{C}$, there is a biequivalence
\[\mathbf{Lax}[\mathcal{B}, \mathcal{C}] \simeq \mathbf{Ps}_{\mathbf{Lax}}[\widehat{\mathcal{B}}, \mathcal{C}]\]
\begin{proof}
We start by defining 2-functors
\[(\widehat{-}) : \mathbf{Lax}[\mathcal{B}, \mathcal{C}] \to \mathbf{Ps}_{\mathbf{Lax}}[\widehat{\mathcal{B}}, \mathcal{C}]\]
\[(\widetilde{-}) : \mathbf{Ps}_{\mathbf{Lax}}[\widehat{\mathcal{B}}, \mathcal{C}] \to \mathbf{Lax}[\mathcal{B}, \mathcal{C}]\]
then show they are each other's inverses.

Start with $\widehat{(-)}$: we define
\begin{enumerate}
    \item On objects (that is, lax functors $\mathcal{F} : \mathcal{B} \to \mathcal{C}$) define $\widehat{\mathcal{F}}$ as
    \begin{itemize}
        \item On objects $a : \widehat{\mathcal{B}}$, define  $\widehat{\mathcal{F}}(a) := \mathcal{F}(a)$,
        \item On morphisms $f = (k_f, \mathcal{H}_f) : a \to b : \widehat{\mathcal{B}}$, define $\widehat{\mathcal{F}}(f) := \mathcal{H}_{\widehat{\mathcal{F}}(f)}(0 \le k_f)$ where
        \[\mathcal{H}_{\widehat{\mathcal{F}}(f)}(i) := \mathcal{F}(\mathcal{H}(i)) \;\; \forall i : [k_f]\]
        \begin{align*}
        \mathcal{H}_{\widehat{\mathcal{F}}(f)}(i \le j) := \begin{cases}
        \mathbb{1} &\text{  if  } j = i \\
        \mathcal{F}(\mathcal{H}(i \le i+1)) &\text{  if  } j = i+1 \\
        \mathcal{H}_{\widehat{\mathcal{F}}(f)}(i+1 \le j) \cdot \mathcal{H}_{\widehat{\mathcal{F}}(f)}(i \le i+1) &\text{  if  } j > i+1
        \end{cases}
        \end{align*}
        \item On 2-cells $\Gamma = (\phi_\Gamma, \alpha_\Gamma) : f \to g : a \to b : \widehat{\mathcal{B}}$ we further distinguish two cases:
        \begin{itemize}
            \item If $\phi_\Gamma$ is an identity, we decompose\footnote{Recall how $I_+ \simeq \Delta_+$, hence its morphisms are generated under composition and ordinal sum by the comonoid $([1], \epsilon, \delta)$} $\alpha$. We define $\widehat{\mathcal{F}}(\mathbb{1}, \epsilon) := \mathcal{F}_u$ and $\widehat{\mathcal{F}}(\mathbb{1}, \delta) := \mathcal{F}_c$.
            \item If $\alpha_\Gamma$ is an identity, we define $\widehat{\mathcal{F}}(\Gamma) := (\phi_{\widehat{\mathcal{F}}(\Gamma)})_{0 \le k_f}$ where
            \begin{align*}
            \widehat{(\phi_\Gamma, \mathbb{1})} := \begin{cases}
            \mathbb{1} &\text{  if  } j = i \\
            \mathcal{F}((\phi_\Gamma)_{i \le i+1}) & \text{  if  } j = i+1 \\
            (\phi_{\widehat{\mathcal{F}}(\Gamma)})_{i+1 \le j} \cdot (\phi_{\widehat{\mathcal{F}}(\Gamma)})_{i \le i+1} &\text{  if  } j > i+1
            \end{cases}
            \end{align*}
        \end{itemize}
        \item The unit constraint $\widehat{\mathcal{F}}_u$ is the identity,
        \item The composition constraint $\widehat{\mathcal{F}}_c$ is the appropriate composition of associators.
    \end{itemize}
    Showing that $\widehat{F}$ is a pseudofunctor is a simple exercise.

    \item On morphisms (that is, lax 2-transformations $\Phi : \mathcal{F} \to \mathcal{G}$) define $\widehat{\Phi}$ as
    \begin{itemize}
        \item on objects $\widehat{\Phi}_a := \Phi_a$
        \item on morphisms $f = (k_f, \mathcal{H}_f) : a \to b$,
        \[\widehat{\Phi}_f := (\Phi_{f_{k_f}} \cdot \mathcal{H}_{\widehat{\mathcal{F}}(f)}(0 \le k_f - 1)) \cdot ... \cdot (\mathcal{H}_{\widehat{\mathcal{F}}(f)}(1 \le k_f) \cdot \Phi_{f_{1}})\]
        (where the horizontal compositions are to be associated to the left)
    \end{itemize}
    As for $\widehat{\mathcal{F}}$, showing that $\widehat{\Phi}$ is a lax 2-transformation is a simple exercise.

    \item On 2-cells (that is, modifications : $\Gamma : \Phi \to \Psi : \mathcal{F} \to \mathcal{G}$) we define $\widehat{\Gamma}_a := \Gamma_a$. Showing it is a modification is, once more, a simple exercise.
\end{enumerate}
Showing that $\widehat{(-)} : \mathbf{Lax}[\mathcal{B}, \mathcal{C}] \to \mathbf{Ps}_{\mathbf{Lax}}[\widehat{\mathcal{B}}, \mathcal{C}]$ is also a straightforward calculation.

Similarly, we define $\widetilde{(-)}$:
\begin{enumerate}
    \item On objects (that is, pseudofunctors $\mathcal{F} : \widehat{\mathcal{B}} \to \mathcal{C}$), we define $\widetilde{\mathcal{F}} : \mathcal{B} \to \mathcal{C}$ by
    \begin{itemize}
        \item On objects $a : \mathcal{C}$, as $\widetilde{\mathcal{F}}(a) = \mathcal{F}(a)$,
        \item On morphisms $f : a \to b$, as $\widetilde{\mathcal{F}}(f) := \mathcal{F}(!_f)$, where $!_f : [1] \to \mathcal{C}$ is the only semistrict pseudofunctor such that $(!_f)(0 \le 1) = f$
        \item On 2-cells $\phi : f \to g : a \to b$, as $\widetilde{\mathcal{F}}(\phi) := \mathcal{F}(!_\phi, \mathbb{1})$ where $!_\phi : !_f \to !_g$ is the only ICON with $(!_\phi)_{0 \le 1} = \phi$.
        \item The unit constraint $\widetilde{\mathcal{F}}_u$ is given by $\mathcal{F}(\mathbb{1}, \epsilon)$
        \item The composition constraint $\widetilde{\mathcal{F}}$ is given by $\mathcal{F}(\mathbb{1}, \delta)$.
    \end{itemize}
    Once again, showing this defines a lax functor is a simple exercise.
    
    \item On morphisms (that is, lax 2-transformations $\Phi : \mathcal{F} \to \mathcal{G}$), we define $\widetilde{\Phi}$ by
    \begin{itemize}
        \item On objects, $\widetilde{\Phi}_a := \Phi_a$,
        \item On morphisms, $\widetilde{\Phi}_f := \Phi_{!_f}$
    \end{itemize}
    Which can be easily seen as a lax 2-transformation.

    \item On 2-cells (that is, modifications $\Gamma : \Phi \to \Psi : \mathcal{F} \to \mathcal{G}$), we just stipulate $\widetilde{\Gamma}_a := \Gamma_a$. This clearly defines a modification.
\end{enumerate}

What's left to check here is that the functors $(\widetilde{-})$ and $(\widehat{-})$ define a biequivalence (in fact, this turns out to be an equivalence of 2-categories). This is easily verified to be the case.
\end{proof}
\end{thm}

This construction is closely related to (and heavily inspired by) the notion of a \textit{lax morphism classifier}, introduced by Blackwell, Kelly and Power \cite{2-dimensional-monad-theory}. The following lemma is immediate:

\begin{lm}\thmheader{Lax functor classifier as lax morphism classifier}

Given a monoidal category $\mathbb{C} := (\mathcal{C}, \otimes, I, \dots)$, the delooping of the lax morphism classifier is biequivalent to the lax functor classifier of the delooping:
\[\mathbb{B}\widehat{\mathbb{C}} \simeq \widehat{\mathbb{B}\mathbb{C}}\]
\end{lm}
As an immediate corollary, we get a description of the lax functor classifier of the terminal 2-category:
\[\widehat{*} = \widehat{\mathbb{B} *} \simeq \mathbb{B} \widehat{*} \simeq \mathbb{B}\Delta_+\]
From this it's easy to see that
\[\mathbf{Mnd}(\mathcal{B}) \simeq \mathbf{Lax}[\mathbb{1}, \mathcal{B}] \simeq \mathbf{St}_{\mathbf{Lax}}[\mathbb{B} \Delta_\alpha, \mathcal{B}]\]
\[\mathbf{Dist}(\mathcal{B}) := \mathbf{Mnd}(\mathbf{Mnd}(\mathcal{B})) \simeq \mathbf{St}_{\mathbf{Lax}}[\mathbb{B} \Delta_\alpha, \mathbf{St}_{\mathbf{Lax}}[\mathbb{B}\Delta_\alpha, \mathcal{B}]]\]

It would be convenient to ``curry'' this presentation of distributive laws, so that we could characterize $\mathbf{Mnd}^n(\mathcal{B})$ as a 2-functor category with codomain $\mathcal{B}$ (rather than an \textit{iterated} 2-functor category). We can do something very close by employing (the lax version of) the Gray tensor product.

%% ## SECTION 3
\section{The Gray tensor product}

Consider the \textit{category} $\mathbf{2Cat}$ of strict 2-categories and strict 2-functors. As it turns out (cft. \cite{gray-product} theorems I.4.9 and I.4.14), this category can be endowed with the structure of a monoidal closed category in such a way that $\mathbf{Ps}[\mathcal{C}, \mathcal{D}]$ is its internal hom: its left adjoint is then called the ``Gray tensor product'', usually written as just $\otimes$. There's also variants\footnote{The lax variant is the one that was originally discussed in Gray's paper, thm I.4.9 in \cite{gray-product}} that make $\mathbf{Lax}$, or $\mathbf{CoLax}$, the internal hom; we'll refer to them as $\otimes_l$ and $\otimes_c$.

Spelling out the adjunction in the lax case, we get
\[\mathbf{2Cat}[\mathcal{B}, \mathbf{Lax}[\mathcal{C}, \mathcal{D}]] \simeq \mathbf{2Cat}[\mathcal{B} \otimes_l \mathcal{C}, \mathcal{D}]\]

It is important to remember that what we get is \textit{not} an equivalence of 2-categories: $\mathbf{2Cat}$ is a category, and the above construction is known to not be amenable to be extended to a 2-functorial (or 3-functorial) one\footnote{See this post by T. Champion for a counterexample: \cite{Gray-not-functorial}}. What follows is a sketch of the explicit construction for the lax Gray tensor product:

\begin{thm}\thmheader{The lax Gray tensor product}

Given 2-categories $\mathcal{B}$, $\mathcal{C}$, we define $\mathcal{B} \otimes_l \mathcal{C}$ to be the 2-category whose
\begin{enumerate}
\item Objects $(b, c) : \mathcal{B} \otimes_l \mathcal{C}$ are pairs of objects $b : \mathcal{B}$, $c : \mathcal{C}$;
\item 1-cells $h : (b , c) \to (b\prime , c\prime) : \mathcal{B} \otimes_l \mathcal{C}$ are composable sequences of pairs of 1-cells $h = (f_1 , g_2) \cdot ... \cdot (f_n, g_n)$ where for each $i \in \{1, \dots, n\}$ either $f = \mathbf{Id}$ or $g = \mathbf{Id}$, subject to the relations 
\[(f_1 , \mathbf{Id}) \circ (f_2 , \mathbf{Id}) = (f_1 \circ f_2 , \mathbf{Id})\]
\[(\mathbf{Id} , g_1) \circ (\mathbf{Id} , g_2) = (\mathbf{Id} , g_1 \circ g_2)\]
\item 2-cells are generated (under horizontal and vertical composition) by the following three families:
\[(\phi , \mathbf{Id}) : (f_1 , \mathbf{Id}) \to (f_2 , \mathbf{Id})\]
\[(\mathbf{Id}, \psi) : (\mathbf{Id} , g_1) \to (\mathbf{Id}, g_2)\]
\[\gamma_{f , g} : (f , \mathbf{Id}) \circ (\mathbf{Id}, g) \to (\mathbf{Id}, g) \circ (f , \mathbf{Id})\]
subject to the following relations, where we write $\circ$ for vertical and $\cdot$ for horizontal composition:
\begin{enumerate}
\item for any $f_3 \xrightarrow{\phi_2} f_2 \xrightarrow{\phi_1} f_1 : b_1 \to b_2 : \mathcal{B}$ and $c : \mathcal{C}$,
\[(\phi_1 , \mathbf{Id}) \circ (\phi_2 , \mathbf{Id}) = (\phi_1 \circ \phi_2 , \mathbf{Id}) : (f_3 , \mathbf{Id}) \to (f_1 , \mathbf{Id})\]
\item for any $g_3 \xrightarrow{\psi_2} g_2 \xrightarrow{\psi_1} g_1 : c_1 \to c_2 : \mathcal{C}$ and $b : \mathcal{B}$,
\[(\mathbf{Id}, \psi_1) \circ (\mathbf{Id}, \psi_2) = (\mathbf{Id}, \psi_1 \circ \psi_2) : (\mathbf{Id}, g_3) \to (\mathbf{Id}, g_1) :\]
\item for any $f_1 \xrightarrow{\phi_1} f_2 : b_2 \to b_3 : \mathcal{B}$, $f_3 \xrightarrow{\phi_2} f_4 : b_1 \to b_2 : \mathcal{B}$ and $c : \mathcal{C}$,
\[(\phi_1 , \mathbf{Id}) \cdot (\phi_2 , \mathbf{Id}) = (\phi_1 \cdot \phi_2 , \mathbf{Id}) : (f_1 \cdot f_3 , \mathbf{Id}) \to (f_2 \cdot f_4 , \mathbf{Id})\]
\item for any $g_1 \xrightarrow{\psi_1} g_2 : c_2 \to c_3 : \mathcal{C}$, $g_3 \xrightarrow{\psi_2} g_4 : c_1 \to c_2 : \mathcal{C}$ and $b : \mathcal{B}$,
\[(\mathbf{Id}, \psi_1) \cdot (\mathbf{Id}, \psi_2) = (\mathbf{Id}, \psi_1 \cdot \psi_2) : (\mathbf{Id}, g_1 \cdot g_3) \to (\mathbf{Id}, g_2 \cdot g_4)\]
\item for all $f : b_1 \to b_2 : \mathcal{B}$,
\[\gamma_{f , \mathbf{Id}} = \mathbf{Id} : (f , \mathbf{Id}) \to (f , \mathbf{Id})\]
\item for all $g : c_1 \to c_2 : \mathcal{C}$,
\[\gamma_{\mathbf{Id} , g} =  \mathbf{Id} : (\mathbf{Id}, g) \to (\mathbf{Id}, g) \]
\item for any $b_3 \xrightarrow{f_2} b_2 \xrightarrow{f_1} b_1 : \mathcal{B}$ and $g : c_1 \to c_2 : \mathcal{C}$,
\[\gamma_{(f_1 \cdot f_2), g} = (\gamma_{f_1 , g} \cdot (f_2 , \mathbf{Id})) \circ ((f_1 , \mathbf{Id}) \cdot \gamma_{f_2 , g}) : (f_1 \cdot f_2 , \mathbf{Id}) \cdot (\mathbf{Id}, g) \to (\mathbf{Id}, g) \cdot (f_1 \cdot f_2, \mathbf{Id})\]
\item for any $c_3 \xrightarrow{g_2} c_2 \xrightarrow{g_1} c_1 : \mathcal{C}$ and $f : b_1 \to b_2 : \mathcal{B}$,
\[\gamma_{f , (g_1 \cdot g_2)} = ((\mathbf{Id} , g_1) \cdot \gamma_{f , g_2}) \circ (\gamma_{f , g_1} \cdot (\mathbf{Id}, g_2)) : (f , \mathbf{Id}) \cdot (\mathbf{Id}, g_1 \cdot g_2) \to (\mathbf{Id}, g_1 \cdot g_2) \cdot (f, \mathbf{Id})\]
\item for any $f_1 \xrightarrow{\phi} f_2 : b_1 \to b_2 : \mathcal{B}$ and $g_1 \xrightarrow{\psi} g_2 : c_1 \to c_2 : \mathcal{C}$,
\[\gamma_{f_2 , g_2} \circ ((\phi , \mathbf{Id}) \cdot (\mathbf{Id}, \psi)) = ((\mathbf{Id} , \psi) \cdot (\phi , \mathbf{Id})) \circ \gamma_{f_1 , g_1} : (f_1 , \mathbf{Id}) \cdot (\mathbf{Id}, g_1) \to (\mathbf{Id}, g_2) \cdot (f_2 , \mathbf{Id})\]
\end{enumerate}
\end{enumerate}
\end{thm}

With this tool, we can get a (1-)categorical description of parametric monads and parametric distributive laws. If we replace the 2-category of parametric distributive laws
\[\mathbf{PDist}(\mathcal{B}, \mathcal{C}) := \mathbf{St}_{\mathbf{St}}[\mathcal{B}, \mathbf{Dist}(\mathcal{C})]\]
with a less structured \textit{category} of parametric distributive laws
\[PDist(\mathcal{C}, \mathcal{D}) := \mathbf{2Cat}[\mathcal{C}, \mathbf{Dist}(\mathcal{D})]\]
then we can apply the tensor-hom adjunction for the lax Gray tensor product, obtaining the following characterization:

\[PDist(\mathcal{B}, \mathcal{C}) \simeq \mathbf{2Cat}[(\mathcal{B} \otimes_l \mathbb{B} \Delta_\alpha) \otimes_l \mathbb{B} \Delta_\alpha, \mathcal{C}]\]

Something completely analogous can be done for parametric monads $PMnd$ (replacing the previously defined $\mathbf{PMnd}$) 
\[PMnd(\mathcal{B}, \mathcal{C}) \simeq \mathbf{2Cat}[\mathcal{B} \otimes_l \mathbb{B} \Delta_\alpha, \mathcal{C}]\]
and (as we will discuss later) we can do the same for iterated distributive laws. This suggests defining operations on 2-categories
\[\mathbf{mnd}(\mathcal{C}) := \mathcal{C} \otimes_l \mathbb{B}\Delta_\alpha\]
\[\mathbf{dist}(\mathcal{C}) := \mathbf{mnd}^2(\mathcal{C}) = (\mathcal{C} \otimes_l \mathbb{B}\Delta_\alpha) \otimes_l \mathbb{B}\Delta_\alpha\]
so that
\[PMnd(\mathcal{C}, \mathcal{B}) \simeq \mathbf{2Cat}[\mathbf{mnd}(\mathcal{C}), \mathcal{B}]\]
\[PDist(\mathcal{C}, \mathcal{B}) \simeq \mathbf{2Cat}[\mathbf{dist}(\mathcal{C}), \mathcal{B}]\]

Before moving on, we want to state the following explicit characterization of iterated Gray tensor products (which will be useful later)

\begin{thm}\thmheader{Iterated lax Gray tensor products}

Given $n$ 2-categories $\mathcal{C}_1, \dots, \mathcal{C}_n$, we define $\bigotimes_{i = 1}^{n} \mathcal{C}_i$ by $\bigotimes_{i = 1}^{1} := \mathcal{C}_1$ and $\bigotimes_{i = 1}^{k+1}\mathcal{C}_i := (\bigotimes_{i = 1}^{k} \mathcal{C}_i) \otimes_l \mathcal{C}_{k+1}$.

\begin{enumerate}
\item the 0-cells in $\bigotimes_{i = 1}^{n} \mathcal{C}_i$ are $n$-tuples of 0-cells $(a_1, \dots, a_n)$ with $a_i : \mathcal{C}_i$;
\item the 1-cells $(a_1, \dots, a_n) \to (b_1, \dots, b_n)$ are generated (under composition) by $\{f\}_i : (a_1, \dots, a_i, \dots, a_n) \to (a_1, \dots, a_{i}^{\prime}, \dots, a_n)$ for $f : a_i \to a_{i}^{\prime} : \mathcal{C}_i$, subject to the relation $\{f\}_i \cdot \{g\}_i = \{f \cdot g\}_i$ for all choices of $f, g, i$.
\item the 2-cells are generated, under horizontal and vertical composition, by the following two families:
\begin{itemize}
\item $\{\phi\}_i : \{f\}_i \to \{g\}_i$ for $\phi : f \to g$ in $\mathcal{C}_i$
\item $\gamma^{i, j}_{f, g} : \{f\}_i \cdot \{g\}_j \to \{g\}_j \cdot \{f\}_i$ for any $i < j$
\end{itemize}
subject to the following relations (whenever well-typed):
\begin{enumerate}
\item compatibility with vertical and horizontal composition: for each $i$
\[\{\phi\}_i \circ \{\psi\}_i = \{\phi \circ \psi\}_i \;\; , \;\; \{\phi\}_i \cdot \{\psi\}_i = \{\phi \cdot \psi\}_i\]
\item compatibility with identities: for each $i < j$
\[\gamma^{i, j}_{f, \mathbf{Id}} = \{f\}_i \;\; , \;\; \gamma^{i, j}_{\mathbf{Id}, f} = \{f\}_j\]
\item compatibility with composition of 1-cells: for each $i < j$
\[\gamma^{i, j}_{(f_1 \cdot f_1), g} = (\gamma^{i, j}_{f_1, g} \cdot \{f_2\}_i) \circ (\{f_1\}_i \cdot \gamma^{i, j}_{f_2, g})\]
\[\gamma^{i, j}_{f, (g_1 \cdot g_2)} = (\{g_1\}_j \cdot \gamma^{i, j}_{f, g_2}) \circ (\gamma^{i, j}_{f, g_1} \cdot \{g_2\}_j)\]
\item naturality: for each $i < j$
\[(\{\psi\}_j \cdot \{\phi\}_i) \circ \gamma^{i, j}_{f_1, g_1} = \gamma^{i, j}_{f_2, g_2} \circ (\{\phi\}_i \cdot \{\psi\}_j)\]
\item the Yang-Baxter equation(s): for each $i < j < k$
\[(\{h\}_k \cdot \gamma^{i, j}_{f, g}) \circ (\gamma^{i, k}_{f, h} \cdot \{g\}_j) \circ (\{f\}_i \cdot \gamma^{j, k}_{g, h}) = (\gamma^{j, k}_{g, h} \cdot \{f\}_i) \circ (\{g\}_j \cdot \gamma^{i, k}_{f, h}) \circ (\gamma^{i, j}_{f, g} \cdot \{h\}_k)\]
\end{enumerate}
\end{enumerate}
\begin{proof}
By induction on $n$; the base case is trivial. Indeed, for $n = 2$ conditions (a - d) are exactly the ones spelled out in \ref{The lax Gray tensor product}. Notice the shift in notation: we write $\{f\}_1 := (f , \mathbf{Id})$ and $\{f\}_2 := (\mathbf{Id}, f)$ for 1- and 2-cells, and $\gamma^{1, 2}_{f, g} := \gamma_{f, g}$. Condition (e) doesn't apply, since it requires at least $n = 3$.

Consider now the inductive case $n +1$ (for $n \geq 2$):
\[\bigotimes_{i = 1}^{n+1} \mathcal{C}_i= (\bigotimes_{i=1}^{n} \mathcal{C}_i) \otimes_l \mathcal{C}_{n+1}\]
By inductive hypothesis, the characterization in question applies to $\bigotimes_{i = 1}^{n} \mathcal{C}_i$. For $i < j \leq n$ we write
$\{f\}_i := (\{f\}_i , \mathbf{id})$ and $\{f\}_{n+1} := (\mathbf{Id}, f)$ for 1- and 2-cells, and $\gamma^{i, j}_{f, g} := (\gamma^{i, j}_{f, g} , \mathbf{Id})$. Moreover, we write $\gamma^{i, n+1}_{f, g} := \gamma_{\{f\}_i, g}$.

The characterization for 1-cells holds trivially; we will only discuss the relations imposed on 2-cells. Clearly, if all indices are strictly less than $n+1$ the relations hold by inductive hypothesis;
\begin{enumerate}
\item[(a)] holds by definition for $i = n+1$;
\item[(b)] for $i \leq n$, we want to prove
\[\gamma^{i, n+1}_{f, \mathbf{Id}} := \gamma_{\{f\}_1, \mathbf{Id}} = (\{f\}_i, \mathbf{Id}) =: \{f\}_i\]
\[\gamma^{i, n+1}_{\mathbf{Id}, f} := \gamma_{\mathbf{Id}, f} = (\mathbf{Id}, f) =: \{f\}_{n+1}\]
both of which hold by the definition of Gray tensor product \ref{The lax Gray tensor product}, conditions (e) and (f) respectively.
\item[(c)] for $i \leq n$, we want to prove
\[\gamma^{i, n+1}_{(f_1 \cdot f_2), g} = (\gamma^{i, n+1}_{f_1, g} \cdot \{f_2\}_i) \circ (\{f_1\}_i \cdot \gamma^{i, n+1}_{f_2, g})\]
\[\gamma^{i, n+1}_{f, (g_1 \cdot g_2)} = (\{g_1\}_{n+1} \cdot \gamma^{i, n+1}_{f, g_2}) \circ (\gamma^{i, n+1}_{f, g_1} \cdot \{g_2\}_{n+1})\]
Start with the first equation: the left hand side is defined as
\[\gamma^{i, n+1}_{(f_1 \cdot f_2), g} := \gamma_{\{f_1 \cdot f_2\}_i, g} = \gamma_{(\{f_1\}_i \cdot \{f_2\}_i), g}\]
while the right hand side is
\[(\gamma^{i, n+1}_{f_1, g} \cdot \{f_2\}_i) \circ (\{f_1\}_i \cdot \gamma^{i, n+1}_{f_2, g}) := (\gamma_{\{f_1\}_i , g} \cdot (\{f_2\}_i , \mathbf{Id})) \circ ((\{f_1\}_i , \mathbf{Id}) \cdot \gamma_{\{f_2\}_i , g})\]
So our claim reduces to
\[\gamma_{(\{f_1\}_i \cdot \{f_2\}_i), g} = (\gamma_{\{f_1\}_i , g} \cdot (\{f_2\}_i , \mathbf{Id})) \circ ((\{f_1\}_i , \mathbf{Id}) \cdot \gamma_{\{f_2\}_i , g})\]
which is exactly point (g) in \ref{The lax Gray tensor product}. The proof for the other equation is completely analogous.
\item[(d)] for $i \leq n$ we want to prove
\[(\{\psi\}_{n+1} \cdot \{\phi\}_i) \circ \gamma^{i, n+1}_{f_1 , g_1} = \gamma^{i, n+1}_{f_2 , g_2} \circ (\{\phi\}_i \cdot \{\psi\}_{n+1})\]
Again, we unroll our notation obtaining
\[((\mathbf{Id}, \psi) \cdot (\{\phi\}_i , \mathbf{Id})) \circ \gamma_{\{f_1\}_i , g_1} = \gamma_{\{f_2\}_i , g_2} \circ ((\{\phi\}_i , \mathbf{Id}) \cdot (\mathbf{Id}, \psi))\]
which is exactly condition (i) in \ref{The lax Gray tensor product}.
\item[(e)] for $i < j \leq n$, we want to prove
\[(\{h\}_{n+1} \cdot \gamma^{i, j}_{f, g}) \circ (\gamma^{i, n+1}_{f, h} \cdot \{g\}_j) \circ (\{f\}_i \cdot \gamma^{j, n+1}_{g, h}) = (\gamma^{j, n+1}_{g, h} \cdot \{f\}_i) \circ (\{g\}_j \cdot \gamma^{i, n+1}_{f, h}) \circ (\gamma^{i, j}_{f, g} \cdot \{h\}_{n+1})\]
Unrolling our notation one last time, we get
\[((\mathbf{Id} , h) \cdot (\gamma^{i,j}_{f, g}, \mathbf{Id})) \circ (\gamma_{\{f\}_i, h} \cdot (\{g\}_i , \mathbf{Id})) \circ ((\{f\}_i , \mathbf{Id}) \cdot \gamma_{\{g\}_j, h})\]
for the left hand side, and
\[(\gamma_{\{g\}_{j} , h} \cdot (\{f\}_i , \mathbf{Id})) \circ ((\{g\}_j , \mathbf{Id}) \cdot \gamma_{\{f\}_i, h}) \circ ((\gamma^{i, j}_{f, g}, \mathbf{Id}) \cdot (\mathbf{Id}, h))\]
for the right hand side. As is easy to see, we can apply point (g) in \ref{The lax Gray tensor product} on both sides and get
\[((\mathbf{Id} , h) \cdot (\gamma^{i,j}_{f, g}, \mathbf{Id})) \circ \gamma_{(\{f\}_i \cdot \{g\}_j), h}\]
for the left hand side, and
\[\gamma_{(\{g\}_j \cdot \{f\}_i), h} \circ ((\gamma^{i, j}_{f, g}, \mathbf{Id}) \cdot (\mathbf{Id}, h))\]
for the right hand side. We can now invoke point (i) in \ref{The lax Gray tensor product} to conclude that they must equal each other.
\end{enumerate}
\end{proof}
\end{thm}

As a corollary of this, whenever we perform two or more applications of the Gray tensor product, we are going to introduce some Yang-Baxter equations.

%% ## SECTION 4
\section{Parametric distributive laws, explicitly}

In this section, we give a concrete description of the data that goes into a parametric distributive law. In order to do so, we will specifically look at $\mathcal{C}$-parametric distributive laws in the case where the 2-category of parameters is the terminal 2-category $\mathcal{C} = \mathbb{1}$.

As is easy to compute $\mathbb{1} \otimes_l \mathcal{C} = \mathcal{C}$, hence $\mathbb{1} \otimes_l \mathbb{B}\Delta_\alpha = \mathbb{B}\Delta_\alpha$. This means that the category of $\mathbb{1}$-parametric distributive laws\footnote{Notice how this is different from the 2-category of distributive laws $\mathbf{Dist}(\mathcal{B}) \simeq \mathbf{St}_{\mathbf{Lax}}[\mathbb{B} \Delta_\alpha, \mathbf{St}_{\mathbf{lax}}[\mathbb{B} \Delta_\alpha, \mathcal{B}]]$: while 0-cells are the same, the 1-cells differ already! Moreover, the former is ``just'' a category while the latter is a fully fledged 2-category.}
\[PDist(\mathbb{1}, \mathcal{B}) \simeq \mathbf{2Cat}[\mathbb{B} \Delta_\alpha \otimes_l \mathbb{B} \Delta_\alpha, \mathcal{B}]\]

It is worth spelling out the data that goes into such a parametric distributive law more explicitly, stressing what pieces of data correspond to the two monads and what to the actual distributive law.

\begin{rmk}\thmheader{Distributive laws as strict functors}

A $\mathbb{1}$-parametric distributive law $\mathcal{D} : PDist(\mathbb{1}, \mathcal{B})$ is given by the following data
\begin{enumerate}
\item A $0$-cell $a : \mathcal{B}$,
\item $1$-cells $s, t : a \to a : \mathcal{B}$
\item $2$-cells $\eta^s : \mathbf{Id}_a \to s$, $\mu^s : s \cdot s \to s$, $\eta^t : \mathbf{Id}_a \to t$, $\mu^t : t \cdot t \to t$ and $\sigma : s \cdot t \to t \cdot s$
\end{enumerate}
subject to the appropriate axioms: $(s, \eta^s, \mu^s)$ and $(t, \eta^t, \mu^t)$ are required to be monads, and $\sigma$ is required to be a distributive law in the sense of J. Beck (\cite{distributive-laws}).
\begin{proof}

Recall how $\Delta_\alpha$ is the free monoidal category with a monoid object\footnote{more precisely: it is the initial object in the 2-category whose 0-cells are monoidal categories with a specified monoid object, whose 1-cells are ``pointed'' monoidal functors, and whose 2-cells are monoidal natural transformations}. This implies that 
\begin{enumerate}
\item Its objects are generated (under its tensor product, sometimes called ordinal sum $\oplus$) by the singleton ordinal $[0]$: the empty ordinal $[-1]$ would then be its 0-fold ordinal sum; and in general we have that $[n] = [0]^{\oplus n+1}$.
\item Its morphisms are generated (under both composition and ordinal sum) by the unit and multiplication of the monoid structure on $[0]$: $\eta : [-1] \to [0]$ and $\mu : [1] \to [0]$.
\end{enumerate}
Hence its delooping admits an almost identical characterization for its 0-, 1- and 2-cells. From the construction of the Gray tensor product, we can then notice that in order to define a strict 2-functor $\mathbb{B}\Delta_\alpha \otimes_l \mathbb{B}\Delta_\alpha \to \mathcal{B}$ we only need to specify where the two copies of the data we just listed end up, together with the ``swap'' 2-cell $\gamma_{[0], [0]}$. Explicitly, we only need to care about
\begin{enumerate}
\item The 0-cell $(*, *)$, whose image $\mathcal{D}(*, *)$ we can call $a$;
\item the 1-cells $([0], \mathbf{Id})$ and $(\mathbf{Id}, [0])$, whose image we can call $s$ and $t$ respectively;
\item the 2-cells $(\eta, \mathbf{Id})$, $(\mu, \mathbf{Id})$, $(\mathbf{Id}, \eta)$, $(\mathbf{Id}, \mu)$ and $\gamma_{[0], [0]}$. We can call their images $\eta^s, \mu^s, \eta^t, \mu^t, \sigma$ respectively.
\end{enumerate}
\end{proof}
\end{rmk}

This should not be surprising: we set up our notion of parametric distributive laws exactly for this result to follow. What might turn out to be (at least slightly) surprising is that we can now produce a similarly explicit description of parametric iterated distributive laws too, which we define just as

\[PDist^n(\mathcal{C}, \mathcal{B}) := \mathbf{2Cat}[\mathcal{C}, \mathbf{Mnd}^n(\mathcal{B})] \simeq \mathbf{2Cat}[\mathbf{mnd}^n(\mathcal{C}), \mathcal{B}]\]

Notice how, as an artifact of our counting convention, a 1-fold parametric distributive law is a parametric monad; and a 2-fold such is a parametric distributive law. The definition produces something new for $n > 2$.

\begin{rmk}\thmheader{3-fold distributive laws as strict functors}

A $\mathbb{1}$-parametric 3-fold distributive law $\mathcal{D} : PDist^3(\mathbb{1}, \mathcal{B})$ is given by the data of
\begin{enumerate}
\item A 0-cell $a : \mathcal{B}$,
\item 1-cells $t_1, t_2, t_3 : a \to a$,
\item 2-cells $\eta^i : \mathbf{Id}_a \to t_i$ and $\mu^i : t_i \cdot t_i \to t_i$ for $i \in \{1, 2, 3\}$, as well as $\sigma_{i, j}$ as before for $i < j \in \{1, 2, 3\}$
\end{enumerate}
subject to the appropriate axioms: $(t_i, \eta^i, \mu^i)$ are required to be monads, $\gamma_{i, j}$ to be distributive laws, and (as spelled out by E. Cheng in \cite{iterated-distributive-laws}) such distributive laws are required to satisfy the Yang-Baxter equation.
\begin{proof}
This is a corollary of the characterization of $\mathbb{1}$-parametric distributive laws in \ref{Distributive laws as strict functors} and our remarks about the universal property of $\Delta_\alpha$ therein, together with \ref{Iterated lax Gray tensor products}. 
\end{proof}
\end{rmk}

The previous result immediately generalizes to $\mathbb{1}$-parametric n-fold distributive laws, without any significant change to the proof:

\begin{rmk}\thmheader{n-fold distributive laws as strict functors}

A $\mathbb{1}$-parametric n-fold distributive law $\mathcal{D} : PDist^n(\mathbb{1}, \mathcal{B})$ is given by the following data:
\begin{enumerate}
\item A 0-cell $a : \mathcal{B}$,
\item 1-cells $t_1, \dots, t_n : a \to a$,
\item 2-cells $\eta^i$ and $\mu^i$ for $1 \le i \le n$ making $t_i$ a monad on $a$
\item 2-cells $\sigma_{i, j}$ for $1 \le i < j \le n$ distributive laws, satisfying the Yang-Baxter equation when appropriate, i.e. for all choices of $1 \le i < j < k \le n$.
\end{enumerate}
\end{rmk}

%% ## SECTION 5
\section{Example: the Writer monad}

In this section, we will show how the ``Writer'' monad admits a parametric distributive law. We will achieve this by explicitly constructing a lift of an arbitrary (appropriately structured) monad to its Eilenberg-Moore category.

First, some definitions:

\begin{defin}\thmheader{The Writer monad}
Let $\mathcal{C}$ be a category with a terminal object $*$ and all binary products $A \times B$. Given a monoid object $\mathbf{M} = (M, u, m) : \mathbf{Mon}(\mathcal{C})$, we can define a functor
\[\mathbf{Writer}_\mathbf{M} : \mathcal{C} \to \mathcal{C}\]
\[\mathbf{Writer}_\mathbf{M}(A \xrightarrow{f} B) := (M \times A) \xrightarrow{M \times f} (M \times B)\]
As well as natural transformations
\[\eta : \mathbf{Id} \to \mathbf{Writer}_\mathbf{M}\]
\[\mu : \mathbf{Writer}_{\mathbf{M}}^{2} \to \mathbf{Writer}_\mathbf{M}\]
by setting
\[\eta_A = A \simeq * \times A \xrightarrow{u \times A} M \times A \]
\[\mu_A : M \times (M \times A) \simeq (M \times M) \times A \xrightarrow{m \times A} M \times A\]
Naturality is immediate to verify. This choice makes $(\mathbf{Writer}_\mathbf{M}, \eta, \mu)$ into a monad.

\begin{proof}
We need to show that $\eta$ and $\mu$ satisfy unitality and associativity. For unitality, consider

\begin{center}\begin{tikzcd}
& * \times (M \times A) \arrow[dl, "u \times (M \times A)"'] \arrow[d, "\alpha^{-1}"] &
M \times A \arrow[l, "\lambda"'] \arrow[dl, "\lambda \times A"] \arrow[dr, "\rho \times A"'] \arrow[r, "M \times \lambda"] \arrow[dd, equals] &
M \times (* \times A) \arrow[d, "\alpha^{-1}"'] \arrow[dr, "M \times (u \times A)"] & \\
M \times (M \times A) \arrow[dr, "\alpha^{-1}"'] &
(* \times M) \times A \arrow[d, "(u \times M) \times A"] &&
(M \times *) \times A \arrow[d, "(M \times u) \times A"'] &
M \times (M \times A) \arrow[dl, "\alpha^{-1}"] \\
& (M \times M) \times A \arrow[r, "m \times A"'] & M \times A & (M \times M) \times A \arrow[l, "m \times A"] &
\end{tikzcd}\end{center}
The triangles on the top left and right commute by the coherence theorem, the two squares on the sides by interchange. The remaining cells in the center are exactly the unitality assumption on the monoid $\mathbf{M} = (M, u, m)$ tensored by $A$.
A similar argument works for associativity: consider the following, where juxtaposition is to be read as the cartesian product $\times$:

\begin{center}\begin{tikzcd}
M (M (M A)) \arrow[rr, "M \alpha^{-1}"] \arrow[dd, "\alpha^{-1}"'] &&
M ((M M) A) \arrow[r, "M (m A)"] \arrow[d, "\alpha^{-1}"] &
M (M A) \arrow[d, "\alpha^{-1}"] \\
&& (M (M M)) A \arrow[r, "(M m) A"] &
(M M) A \arrow[dd, "m A"] \\
(M M) (M A) \arrow[r, "\alpha^{-1}"] \arrow[d, "m (M A)"'] &
((M M) M) A \arrow[ur, "\alpha A"] \arrow[d, "(m M) A"'] && \\
M (M A) \arrow[r, "\alpha^{-1}"] &
(M M) A \arrow[rr, "m A"] &&
M A
\end{tikzcd}\end{center}
The two square cells commute by interchange, the pentagon on the top left is the pentagon axiom for monoidal categories. The remaining pentagon is the associativity condition on $\mathbf{M} = (M, u, m)$ tensored by $A$.
\end{proof}
\end{defin}

There's a more compact way of describing the above definition: it is the monad sending an object $A : \mathcal{C}$ to the (carrier of the) corresponding \textit{free $\mathbf{M}$-module}

\[\mathbf{M} \times A := (M \times A, m \times A)\]

More explicitly, we have the category $\mathbf{Mod}_\mathbf{M}(\mathcal{C})$ of $\mathbf{M}$-modules\footnote{We don't need to choose any handedness for these modules, since the cartesian monoidal structure is always symmetric; for the sake of concreteness however, whenever such a choice is to be made, we'll work with left modules.} whose objects are pairs $(A, \phi)$ with $\phi : M \times A \to A$ satisfying appropriate unitality and associativity conditions; morphisms $f : (A, \phi) \to (B, \psi)$ are the morphisms in the base category $f : A \to B : \mathcal{C}$ that make the obvious square commute.

The category $\mathbf{Mod}_\mathbf{M}(\mathcal{C})$ always comes equipped with a forgetful functor $\mathcal{U} : \mathbf{Mod}_\mathbf{M}(\mathcal{C}) \to \mathcal{C}$, which turns out to have a left adjoint; the associated monad is exactly $\mathbf{Writer}_\mathbf{M}$.

In the following, we will assume $\mathcal{C}$ to be a Cartesian Closed Category, implying that it is enriched over itself (though somewhat milder assumptions would suffice). We will also assume fixed a monoid object $\mathbf{M} = (M, u, m) : \mathbf{Mon}(\mathcal{C})$.

\begin{lm}\thmheader{Strong and enriched functors}

Given an endofunctor $\mathcal{F} : \mathcal{C} \to \mathcal{C}$, the following are equivalent
\begin{itemize}
\item A tensorial strength $\tau_{A, B} : A \times \mathcal{F}(B) \to \mathcal{F}(A \times B)$ interacting appropriately with the unitors and associators of $\mathcal{C}$;
\item A $\mathcal{C}$-enrichment for $\mathcal{F}$, meaning a natural family of morphisms interacting appropriately with identities and compositions
\[\mathcal{F}_{A, B} : [A, B]_{\mathcal{C}} \to [\mathcal{F}(A), \mathcal{F}(B)]_{\mathcal{C}}\]
\end{itemize}
\begin{proof}
The result is standard, so we only sketch the construction back-and-forth (for a full proof, see sections 3.2 and 3.3 of \cite{strength-enrichment}).

First, assume we have a tensorial strength $\tau_{A, B} : A \times \mathcal{F}(B) \to \mathcal{F}(A \times B)$. We can define
\[[A, B]_{\mathcal{C}} \times \mathcal{F}(A) \xrightarrow{\tau} \mathcal{F}([A, B]_{\mathcal{C}} \times A) \xrightarrow{Ev} \mathcal{F}(B)\]
where $Ev$ is the evaluation map, the counit of the tensor-hom adjunction. The morphism we just defined is the mate (under the tensor-hom adjunction in $\mathcal{C}$) of the desired enrichment
\[\mathcal{F}_{A, B} : [A, B]_{\mathcal{C}} \to [\mathcal{F}(A), \mathcal{F}(B)]_{\mathcal{C}}\]

Second, assume we have an enrichment
\[\mathcal{F}_{A, B} : [A, B]_{\mathcal{C}} \to [\mathcal{F}(A), \mathcal{F}(B)]_{\mathcal{C}}\]
We can define
\[A \times \mathcal{F}(B) \xrightarrow{} [B, (A \times B)]_{\mathcal{C}} \times \mathcal{F}(B) \xrightarrow{} [\mathcal{F}(B), \mathcal{F}(A \times B)]_{\mathcal{C}} \times \mathcal{F}(B) \xrightarrow{} \mathcal{F}(A \times B)\]
\end{proof}
where the first map is the unit of the tensor-hom adjunction (tensored by $\mathcal{F}(B)$), the second is the enrichment, followed by the evaluation map.
\end{lm}

An analogous statement holds for monads: the notions of enriched monad and strong monad coincide (see proposition 5.8 in \cite{strong-enriched-monads}).

We can now state our main result:

\begin{thm}\thmheader{Enriched monads lift to modules}

Given an enriched monad $(\mathcal{T}, \eta, \mu)$ on $\mathcal{C}$, there's a corresponding monad $(\hat{\mathcal{T}}, \hat{\eta}, \hat{\mu})$ on $\mathbf{Mod}_{\mathbf{M}}(\mathcal{C})$ such that
\[\mathcal{U} \circ \hat{\mathcal{T}} = \mathcal{T} \circ \mathcal{U}\]
\begin{proof}

We start by defining the monad $\hat{\mathcal{T}} : \mathbf{Mod}_{\mathbf{M}}(\mathcal{C}) \to \mathbf{Mod}_{\mathbf{M}}(\mathcal{C})$

On objects, we can just define $\hat{\mathcal{T}}(A, \phi) := (\mathcal{T}(A), \hat{\phi})$ with
\[M \times \mathcal{T}(A) \xrightarrow{\sigma} \mathcal{T}(M \times A) \xrightarrow{\mathcal{T}(\phi)} \mathcal{T}(A)\]
where $\sigma$ is the tensorial strength of $\mathcal{T}$.

To see that $\hat{\mathcal{T}}(A, \phi) : \mathbf{Mod}_{\mathbf{M}}(\mathcal{C})$, consider the following:
\begin{center}\begin{tikzcd}
* \times \mathcal{T}(A) \arrow[r, "u \times \mathcal{T}(\phi)"] \arrow[d, "\sigma"] &
M \times \mathcal{T}(A) \arrow[d, "\sigma"] &
M \times \mathcal{T}(M \times A) \arrow[l, "M \times \mathcal{T}(\phi)"'] \arrow[d, "\sigma"'] &
M \times (M \times \mathcal{T}(A)) \arrow[l, "M \times \sigma"'] \arrow[d, "\alpha^{-1}"'] \\
\mathcal{T}(* \times A) \arrow[r, "\mathcal{T}(u \times A)"] \arrow[d, "\mathcal{T}(\sim)"] &
\mathcal{T}(M \times A) \arrow[d, "\mathcal{T}(\phi)"] &
\mathcal{T}(M \times (M \times A)) \arrow[l, "\mathcal{T}(M \times \phi)"'] \arrow[dr, "\mathcal{T}(\alpha^{-1})" description] &
(M \times M) \times \mathcal{T}(A) \arrow[d, "\sigma"'] \\
\mathcal{T}(A) \arrow[r, equals] &
\mathcal{T}(A) &
\mathcal{T}(M \times A) \arrow[l, "\mathcal{T}(\phi)"'] &
\mathcal{T}((M \times M) \times A) \arrow[l, "\mathcal{T}(m \times A)"']
\end{tikzcd}\end{center}
The two squares on the top (left and center) commute by naturality of $\sigma$, the pentagon on the right involving the associator twice by the fact that $(\mathcal{T}, \sigma)$ is a strong functor. The two remaining cells commute by the assumption that $(A, \phi)$ is an $\mathbf{M}$-module.

Let now $f : (A, \phi) \to (B, \psi) : \mathbf{Mod}_{\mathbf{M}}(\mathcal{C})$. We define $\hat{\mathcal{T}}(f) := \mathcal{T}(f)$; to see that it is a morphism of $\mathbf{M}$-modules, consider

\begin{center}\begin{tikzcd}
M \times \mathcal{T}(A) \arrow[r, "M \times \mathcal{T}(f)"] \arrow[d, "\sigma"] &
M \times \mathcal{T}(B) \arrow[d, "\sigma"] \\
\mathcal{T}(M \times A) \arrow[r, "\mathcal{T}(M \times f)"] \arrow[d, "\mathcal{T}(\phi)"] &
\mathcal{T}(M \times B) \arrow[d, "\mathcal{T}(\psi)"] \\
\mathcal{T}(A) \arrow[r, "\mathcal{T}(f)"] &
\mathcal{T}(B)
\end{tikzcd}\end{center}
The top square commutes by naturality of $\sigma$, the bottom one commute by the assumption that $f$ is a morphism of $\mathbf{M}$-modules.

All that remains is defining $\hat{\eta}$ and $\hat{\mu}$, in a way such that $(\hat{\mathcal{T}}, \hat{\eta}, \hat{\mu})$ is a monad on $\mathbf{Mod}_{\mathbf{M}}(\mathcal{C})$; we define $\hat{\eta}_{(A, \phi)} := \eta_A$ and $\hat{\mu}_{(A, \phi)} := \mu_A$.

We need to show that the components $\hat{\eta}_{(A, \phi)}$, $\hat{\mu}_{(A, \phi)}$ are morphisms of $\mathbf{M}$-modules. For $\hat{\eta}$, consider:
\begin{center}\begin{tikzcd}
M \times A \arrow[r, "M \times \eta_A"] \arrow[dd, "\phi"] \arrow[dr, "\eta_{M \times A}"] &
M \times \mathcal{T}(A) \arrow[d, "\sigma"] \\
&
\mathcal{T}(M \times A) \arrow[d, "\mathcal{T}(\phi)"] \\
A \arrow[r, "\eta_A"] &
\mathcal{T}(A)
\end{tikzcd}\end{center}
The triangle commutes by the assumption that $\sigma$ makes $(\mathcal{T}, \eta, \mu)$ into a strong monad, the square by naturality of $\eta$. As for $\hat{\mu}$, consider
\begin{center}\begin{tikzcd}
M \times \mathcal{T}(\mathcal{T}(A)) \arrow[d, "\sigma"] \arrow[r, "M \times \mu_A"] &
M \times \mathcal{T}(A) \arrow[dd, "\sigma"] \\
\mathcal{T}(M \times \mathcal{T}(A)) \arrow[d, "\mathcal{T}(\sigma)"] & \\
\mathcal{T}(\mathcal{T}(M \times A)) \arrow[d, "\mathcal{T}(\mathcal{T}(\phi))"] \arrow[r, "\mu_{M \times A}"] &
\mathcal{T}(M \times A) \arrow[d, "\mathcal{T}(\phi)"] \\
\mathcal{T}(\mathcal{T}(A)) \arrow[r, "\mu_A"] &
\mathcal{T}(A)
\end{tikzcd}\end{center}
The pentagon commutes by $\sigma$ being the tensorial strength of the monad $(\mathcal{T}, \eta, \mu)$, the square by naturality of $\mu$.

It only remains to prove that $(\hat{\mathcal{T}}, \hat{\eta}, \hat{\mu})$ is a monad; the proof follows immediately by the assumption that $(\mathcal{T}, \eta, \mu)$ is.
\end{proof}
\end{thm}

This allows us to construct a parametric distributive law with parameters given by the 2-category of enriched monads over $\mathcal{C}$:

\begin{thm}\thmheader{Parametric distributive law for the Writer monad}

Define $\mathcal{B} := \mathbb{B} \mathbf{Cat}_{\mathcal{C}}[\mathcal{C}, \mathcal{C}]$ to be the delooping of the (monoidal) category of $\mathcal{C}$-enriched endofunctors of $\mathcal{C}$, and $\mathcal{D} := \mathbf{Mnd}(\mathcal{C})$.

There is an object in $PDist(\mathcal{D}, \mathcal{C})$, that is a strict 2-functor

\[(\mathbf{Mnd}(\mathbb{B} \mathbf{Cat}_{\mathcal{C}}[\mathcal{C}, \mathcal{C}]) \otimes_l \mathbb{B} \Delta_\alpha) \otimes_l \Delta_\alpha \to \mathbb{B} \mathbf{Cat}_{\mathcal{C}}[\mathcal{C}, \mathcal{C}]\]
such that the underlying monads, encoded as strict 2-functors
\[\mathbf{Mnd}(\mathbb{B} \mathbf{Cat}_{\mathcal{C}}[\mathcal{C}, \mathcal{C}]) \to \mathbf{Mnd}(\mathbb{B} \mathbf{Cat}_{\mathcal{C}}[\mathcal{C}, \mathcal{C}])\]
are, respectively, constantly the writer monad and the identity.
\begin{proof}
Such a 2-functor would correspond to a 2-functor
\[\mathbf{DWriter}_{\mathbf{M}} : \mathbf{Mnd}(\mathbb{B} \mathbf{Cat}_{\mathcal{C}}[\mathcal{C}, \mathcal{C}]) \to \mathbf{Dist}(\mathbb{B} \mathbf{Cat}_{\mathcal{C}}[\mathcal{C}, \mathcal{C}])\]
which we will describe in detail.

On objects $(*, T, \eta^T, \mu^T)$, the functor acts as follows (under the identification between distributive laws and monads in monads):
\[\mathbf{DWriter}_{\mathbf{M}}(*, T, \eta^T, \mu^T) := ((*, \mathbf{Writer}_{\mathbf{M}}, \eta^{\mathbf{M}}, \mu^{\mathbf{M}}), (T, \sigma^T), \eta^T, \mu^T)\]

On 1-cells of monads $(\Phi, \phi) : (*, T, \eta^T, \mu^T) \to (*, S, \eta^S, \mu^S)$, the functor acts as
\[\mathbf{DWriter}_{\mathbf{M}}(\Phi, \phi) := ((\Phi, \sigma^\Phi), \phi)\]
where $\sigma^\Phi$ is the tensorial strength the functor $\Phi$ is canonically endowed with.

Finally, on 2-cells $\gamma : (\Phi, \phi) \to (\Psi, \psi)$ the functor  acts as
\[\mathbf{DWriter}_{\mathbf{M}}(\gamma) = \gamma\]
Checking that this is indeed a 2-functor is a long exercise.
\end{proof}
\end{thm}

%% ## SECTION 6
\section{Example: the Exception monad}

In this section we will show that the monad computer scientists call ``Either'' admits a parametric distributive law. Our strategy will be similar to what we did for the Writer monad in the section above: we will directly show that any other monad lifts to its Eilenberg-Moore category; everything else will trivially follow.

We start by defining the writer monad and proving an equivalent definition that will prove extremely useful:

\begin{defin}\thmheader{The Either monad}

Let $\mathcal{C}$ be a category with an initial object $\varnothing$ and all binary coproducts $A + B$. Then for any object $X : \mathcal{C}$ we can define a functor
\[\mathbf{Either}_X : \mathcal{C} \to \mathcal{C}\]
\[\mathbf{Either}_X(A \xrightarrow{f} B) := (X + A) \xrightarrow{X + f} (x + B)\]
We also define natural transformations
\[\eta : \mathbf{Id} \to \mathbf{Either}_X\]
\[\mu : \mathbf{Either}^{2}_{X} \to \mathbf{Either}_{X}\]
by setting
\[\eta_A = r : A \to X + A\]
\[\mu_A = \langle \ell , id \rangle : X + (X + A) \to X + A\]
Naturality for $\eta$ and $\mu$ is trivial to verify. This makes the triple $(\mathbf{Either}_x, \eta, \mu)$ into a monad; this is the either monad.

\begin{proof}
We need to prove that $\eta$ and $\mu$ indeed satisfy unitality and associativity. We start with unitality:
\begin{center}\begin{tikzcd}
& X + A \arrow[dl, "{\langle \ell , r \circ r \rangle}"'] \arrow[dr, "r"] \arrow[dd, equals] & \\
X + (X + A) \arrow[dr, "{\langle \ell , id \rangle}"'] && X + (X + A) \arrow[dl, "{\langle \ell , id \rangle}"] \\
& X + A & \\
\end{tikzcd}\end{center}
We want to prove that
\[\langle \ell , id \rangle \circ \langle \ell , r \circ r \rangle = id = \langle \ell , id \rangle \circ r\]
Recall that, by the universal property of coproducts, for any $\phi : B + C \to D : \mathcal{C}$, we get
\[\phi = \langle (\phi \circ \ell) , (\phi \circ r) \rangle\]
So we can just compute:
\[\langle \ell , id \rangle \circ r = id\]
showing right unitality, and
\[\langle \ell , id \rangle \circ \langle \ell , (r \circ r) \rangle \circ \ell = \langle \ell , r \rangle \circ \ell = \ell\]
\[\langle \ell , id \rangle \circ \langle \ell , (r \circ r) \rangle \circ r = \langle \ell , id \rangle \circ r \circ r = id \circ r = r\]
showing left unitality (since $id = \langle \ell , r \rangle$).

The argument for associativity is similar:
\begin{center}\begin{tikzcd}
X + (X + (X + A)) \arrow[rr, "{\langle \ell , r \circ \langle \ell , id \rangle \rangle}"] \arrow[dd, "{\langle \ell , id \rangle}"'] && X + (X + A) \arrow[dd, "{\langle \ell , id \rangle}"] \\ \\
X + (X + A) \arrow[rr, "{\langle \ell , id \rangle}"'] && X + A
\end{tikzcd} \end{center}
We want to show the two morphisms are equal, and by the aforementioned universal property it's enough to show they are equal when composed with the coprojections. So we can compute:
\[\langle \ell , id \rangle \circ \langle \ell , id \rangle \circ \ell = \langle \ell , id \rangle \circ \ell = \ell\]
\[\langle \ell , id \rangle \circ \langle \ell , r \circ \langle \ell , id \rangle \rangle \circ \ell = \langle \ell , id \rangle \circ \ell = \ell\]
and similarly
\[\langle \ell , id \rangle \circ \langle \ell , id \rangle \circ r = \langle \ell , id \rangle \circ id = \langle \ell , id \rangle\]
\[\langle \ell , id \rangle \circ \langle \ell , r \circ \langle \ell , id \rangle \rangle \circ r = \langle \ell , id \rangle \circ r \circ \langle \ell , id \rangle = id \circ \langle \ell , id \rangle = \langle \ell , id \rangle\]
\end{proof}
\end{defin}

There is a more compact way of defining this monad, as the monad whose algebras are modules over monoids in the cocartesian monoidal structure $(\mathcal{C}, +, \varnothing, \dots)$ (which is the monoidal structure induced by the initial object $\varnothing$ and coproducts $+$). We'll go through the details below, as this characterization will turn out very useful.

\begin{defin}\thmheader{Cocartesian monoidal categories}

Let $\mathcal{C}$ admit an initial object $\varnothing$ and binary coproducts $+$. Then we can endow $\mathcal{C}$ with a monoidal structure $(\mathcal{C}, +, \varnothing, \lambda, \rho, \alpha)$, called the cocartesian monoidal structure.

\begin{proof}
We need to define the unitors $\lambda, \rho$ and the associator $\alpha$, as well as to prove the triangle and pentagon identities.
\begin{enumerate}
\item $\lambda_A : \varnothing + A \xrightarrow{\langle \epsilon_A , id \rangle} A$, where $\epsilon_A : \varnothing \to A$ is the initial morphism;
\item $\rho_A : A + \varnothing \xrightarrow{\langle id , \epsilon_A \rangle} A$ with the same definition for $\epsilon_A$;
\item $\alpha_{A, B, C} : (A + B) + C \xrightarrow{\langle \langle \ell , r \circ \ell \rangle , r \circ r \rangle} A + (B + C)$
\end{enumerate}
each with the obvious inverse. The triangle identity is then equivalently
\[\rho + B = (A + \lambda) \circ \alpha : (A + \varnothing) + B \to A + B\]
The two morphisms are equal if and only if they are equal when composed with all the canonical injections
\[\ell \circ \ell : A \to (A + \varnothing) + B \text{ ; } \ell \circ r : \varnothing \to (A + \varnothing) + B \text{ ; } r : B \to (A + \varnothing) + B\]
Before proceeding with the proof, recall that given $f : X \to X^\prime$, $g : Y \to Y^\prime$, we have
\[(f + g) := \langle \ell \circ f , r \circ g \rangle : X + Y \to X^\prime + Y^\prime\]

We can now compute:
\begin{enumerate}
\item $\ell \circ \ell$:
\[(\rho + b) \circ \ell \circ \ell = \langle \ell \circ \rho , r \rangle \circ \ell \circ \ell = \ell \circ \rho \circ \ell = \ell \circ \langle id , \epsilon \rangle \circ \ell = \ell \circ id = \ell\]
\[(A + \lambda) \circ \alpha \circ \ell \circ \ell = (A + \lambda) \circ \langle \langle \ell , r \circ \ell \rangle , r \circ r \rangle \circ \ell \circ \ell = (A + \lambda) \circ \langle \ell , r \circ \ell \rangle \circ \ell = (A + \lambda) \circ \ell = \langle \ell , r \circ \lambda \rangle \circ \ell = \ell\]
\item $\ell \circ r$: since the domain is $\varnothing$, the claim follows from initiality
\item $r$:
\[(\rho + B) \circ r = \langle \ell \circ \rho , r \rangle \circ r = r\]
\[(A + \lambda) \circ \alpha \circ r = (A + \lambda) \circ \langle \langle \ell , r \circ \ell \rangle , r \circ r \rangle \circ r = (A + \lambda) \circ r \circ r = \langle \ell , r \circ \lambda \rangle \circ r \circ r = r \circ \lambda \circ r = r \circ \langle \epsilon , id \rangle \circ r = r \circ id = r\]
\end{enumerate}
A similar argument allows us to prove the pentagon identity; we omit the explicit details for the sake of brevity.
\end{proof}
\end{defin}

\begin{lm}\thmheader{Monoids in cocartesian monoidal categories}

Let $(\mathcal{C}, +, \varnothing)$ be a cocartesian monoidal category. There is an isomorphism of categories
\[\mathbf{Mon}(\mathcal{C}, +, \varnothing, \dots) \simeq \mathcal{C}\]

\begin{proof}
We will only define a functor $\mathcal{F} : \mathcal{C} \to \mathbf{Mon}(\mathcal{C}, +, \varnothing)$ which is inverse to the usual forgetful $\mathcal{U} : \mathbf{Mon}(\mathcal{C}, +, \varnothing) \to \mathcal{C}$.

\[\mathcal{F}(A) := (A, \epsilon, \langle id , id \rangle)\]
is clearly seen to be a monoid object: for unitality
\begin{center}\begin{tikzcd}
\varnothing + A \arrow[d, "\epsilon + A"']
& A \arrow[l, "r"'] \arrow[r, "\ell"] \arrow[d, equals]
& A + \varnothing \arrow[d, "A + \epsilon"] \\
A + A \arrow[r, "{\langle id , id \rangle}"'] &
A &
A + A \arrow[l, "{\langle id , id \rangle}"]
\end{tikzcd}\end{center}
we can compute the following:
\[\langle id , id \rangle \circ (\epsilon + A) \circ r = \langle id , id \rangle \circ \langle \ell \circ \epsilon , r \rangle \circ r = \langle id , id \rangle \circ r = id\]
\[\langle id , id \rangle \circ (A + \epsilon) \circ \ell = \langle id , id \rangle \circ \langle \ell , r \circ \epsilon \rangle \circ \ell = \langle id , id \rangle \circ \ell = id\]
The argument for proving associativity is similar (compose with the canonical injections, then compute) but quite long; we'll omit it. The same applies to the proof that $\mathcal{F}(f) := f$ is a well-defined morphism of monoid objects, for any morphism $f$ in $\mathcal{C}$.

We now need to prove that the functor $\mathcal{F}$ we just defined is the inverse to the usual forgetful functor $\mathcal{U}$. That $\mathcal{U} \circ \mathcal{F} = \mathbf{Id}$ is trivial; we need to prove $\mathcal{F} \circ \mathcal{U} = \mathbf{Id}$.

On objects, this means that for any given monoid $(A, u, m) : \mathbf{Mon}(\mathcal{C}, +, \varnothing, \dots)$ we have $u = \epsilon$ and $m = \langle id , id \rangle$: the first statement follows from initiality of $\varnothing$, while the second statement is a consequence of unitality for $(A, \epsilon, m)$. Indeed, we know that
\[m = \langle m \circ \ell , m \circ r \rangle : A + A \to A\]
and by our previous computation, $(\epsilon + A) \circ r = r$ and $(A + \epsilon) \circ \ell = \ell$; from this and unitality, it follows that
\[m \circ (\epsilon + A) \circ r = m \circ r = id = m \circ \ell = m \circ (A + \epsilon) \circ \ell\]
which implies $m = \langle id , id \rangle$.

Finally, the statement that $\mathcal{F}(\mathcal{U}(f)) = f$ is trivial.
\end{proof}
\end{lm}

It should now be clear that the $\mathbf{Either}$ monad on a given object is the one obtained by ``tensoring'' (i.e. performing the coproduct) with that object (with the monad structure given by the unique monoid structure that object is endowed with).

\begin{lm}\thmheader{Modules over monoids in cocartesian monoidal categories}

Let $\mathcal{C}$ a cocartesian monoidal category, and let $A : \mathcal{C}$ be any object (so that $(A, \epsilon, \langle id , id \rangle)$ is a monoid). Then we have an isomorphism of categories
\[\mathbf{Mod}_{(A, \epsilon, \langle id , id \rangle)}(\mathcal{C}, +, \varnothing, \dots) \simeq A \downarrow\mathcal{C}\]
between the category of modules over a monoid and the coslice category.
\begin{proof}

We define functors
\[\mathcal{F} : \mathbf{Mod}_{(A, \epsilon, \langle id , id \rangle)}(\mathcal{C}, +, \varnothing, \dots) \to A \downarrow \mathcal{C}\]
\[\mathcal{G} : A \downarrow \mathcal{C} \to \mathbf{Mod}_{(A, \epsilon, \langle id , id \rangle)}(\mathcal{C}, +, \varnothing, \dots)\]
and prove they are inverses of each other.

For $\mathcal{F}$, let first $(X, \sigma) : \mathbf{Mod}_{(A, \epsilon, \langle id, id \rangle)}(\mathcal{C}, +, \varnothing, \dots)$. This means that $X : \mathcal{C}$ and $\sigma : A + X \to X$, so there exist morphisms $\phi : A \to X$ and $\psi : X \to X$ in $\mathcal{C}$ such that $\sigma = \langle \phi , \psi \rangle$. The fact that $(X, \sigma) : \mathbf{Mod}_{(A, \epsilon, \langle id , id \rangle)}(\mathcal{C}, +, \varnothing, \dots)$ means that the following two diagrams commute in $\mathcal{C}$:
\begin{center}\begin{tikzcd}
\varnothing + X \arrow[ddr, "\lambda"'] \arrow[r, "\epsilon + X"] &
A + X \arrow[dd, "\sigma"description] &
A + (A + X) \arrow[l, "A + \sigma"'] \arrow[d, "\alpha"] \\
& & (A + A) + X \arrow[d, "{\langle id , id \rangle + X}"] \\
& X & A + X \arrow[l, "\sigma"]
\end{tikzcd}\end{center}
Commutativity of the triangle on the left means that, in particular
\[\sigma \circ (\epsilon + X) \circ r = \lambda \circ r\]
and a routine computation shows that
\[\sigma \circ (\epsilon + X) \circ r = \psi\]
\[\lambda \circ r = id\]
so that $\psi = id$. Notably, $\sigma \circ (\epsilon + X) \circ \ell = \lambda \circ \ell$ is trivial because of $\varnothing$'s initiality. Similarly, the pentagon on the right is also always commutative: $(X, \sigma)$ is an object of $\mathbf{Mod}_{(A, \epsilon, \langle id , id \rangle)}(\mathcal{C}, +, \varnothing, \dots)$ if and only if $\sigma = \langle \phi , id \rangle$. This allows us to define
\[\mathcal{F}(X, \sigma) := (X, \sigma \circ \ell) = (X, \phi)\]
which is an object of $A \downarrow \mathcal{C}$.

Given $(X, \sigma), (Y, \gamma) : \mathbf{Mod}_{(A, \epsilon, \langle id , id \rangle)}(\mathcal{C}, +, \varnothing, \dots)$ and $f : X \to Y : \mathcal{C}$, the condition for $f$ to be a morphism $f : (X, \sigma) \to (Y, \gamma)$ is
\begin{center}\begin{tikzcd}
A + X \arrow[d, "\sigma"'] \arrow[r, "A + f"] & A + Y \arrow[d, "\gamma"] \\
X \arrow[r, "f"] & Y
\end{tikzcd}\end{center}

Since we proved that $\sigma = \langle \phi , id \rangle$ (and similarly, $\gamma = \langle \psi , id \rangle$), the above condition reduces to $f \circ \phi = \psi$. This means that $f : (X, \phi) \to (Y, \psi) : A \downarrow \mathcal{C}$. This completes the definition of $\mathcal{F}$.

It's easy to see that the definition of $\mathcal{G}$ is similarly straightforward: we can define $\mathcal{G}(X, \phi) := (X, \langle \phi , id \rangle)$ and $\mathcal{G}(f) = f$. Indeed, the only non-trivial condition on $\mathcal{G}(X, \phi)$ is unitality (which is satisfied because the action is given by $\langle \phi , id \rangle$), and the condition on $\mathcal{G}(f)$ is equivalent to the condition $f$ satisfies as a morphism in the coslice $A \downarrow \mathcal{C}$.

Checking that the above two functors are (strict) inverses of each other is now completely trivial.
\end{proof}
\end{lm}

This means, since $\mathbf{Either}_A$ is the monad for modules over $A$ (or more precisely, the canonical monoid structure $A$ is endowed with), that its Eilenberg-Moore category is given by the coslice in the previous lemma.

\begin{lm}\thmheader{Lifting pointed endofunctors}

Let $A : \mathcal{C}$ be an object in a cocartesian monoidal category, and $\mathcal{U} : A \downarrow \mathcal{C} \to \mathcal{C}$ the forgetful functor. Let also $(\mathcal{F}, u)$ be a pointed endofunctor on $\mathcal{C}$, meaning
\[u : \mathbf{Id} \to \mathcal{F} : \mathcal{C} \to \mathcal{C}\]

There is a pointed endofunctor $(\hat{\mathcal{F}}, \hat{u})$ on $A \downarrow \mathcal{C}$ such that $\mathcal{U} \circ \hat{\mathcal{F}} = \mathcal{F} \circ \mathcal{U}$.
\begin{proof}

We first define $\hat{\mathcal{F}} : A \downarrow \mathcal{C} \to A \downarrow \mathcal{C}$, by defining for any $(X , \phi) : A \downarrow \mathcal{C}$ the object
\[\hat{\mathcal{F}}(X, \phi) := (\mathcal{F}(X), \mathcal{F}(\phi) \circ u_A)\]
Given a morphism $f : (X, \phi) \to (Y, \psi) : A \downarrow \mathcal{C}$, (that is, a morphism $f : X \to Y : \mathcal{C}$ such that $f \circ \psi = \phi$) we get a morphism 
\[\hat{\mathcal{F}}(f) := \mathcal{F}(f) : \hat{\mathcal{F}}(X, \phi) \to \hat{\mathcal{F}}(Y, \psi) : A \downarrow \mathcal{C}\]
since
\[\mathcal{F}(f) \circ \mathcal{F}(\phi) \circ u_A = \mathcal{F}(f \circ \phi) \circ u_A = \mathcal{F}(\psi) \circ u_A\]
The fact $\hat{\mathcal{F}}$ is a functor (i.e. preserves identities and compositions) follows trivially from the fact that $\mathcal{F}$ is a functor.

We can now define $\hat{u} : \mathbf{Id} \to \hat{\mathcal{F}}$ as $\hat{u}_{(X, \phi)} := u_X$. To check that $\hat{u}_{(X, \phi)} : (X, \phi) \to \mathcal{F}(X, \phi)$, we need to prove that the following square commutes for any $(X, \phi) : A \downarrow \mathcal{C}$:
\begin{center}\begin{tikzcd}
A \arrow[r, "u_A"] \arrow[d, "\phi"'] & \mathcal{F}(A) \arrow[d, "\mathcal{F}(\phi)"] \\
X \arrow[r, "u_X"] & \mathcal{F}(X)
\end{tikzcd}\end{center}
This follows directly from $u$ being a natural transformation. The proof that $\hat{u}$ is indeed a natural transformation is completely analogous. Hence, $(\hat{\mathcal{F}}, \hat{u})$ is a pointed endofunctor of $A \downarrow \mathcal{C}$; furthermore, unwinding the definitions of $\mathcal{U}$ and $\hat{\mathcal{F}}$ proves that
\[\mathcal{U} \circ \hat{\mathcal{F}} = \mathcal{F} \circ \mathcal{U}\]
concluding the proof.
\end{proof}
\end{lm}

\begin{lm}\thmheader{Lifting monads}

With $\mathcal{C}$, $A$, $\mathcal{U}$ as in the above lemma, let $\mathcal{T} = (T, \eta, \mu)$ be a monad on $\mathcal{C}$. There exists a monad $\hat{\mathcal{T}} = (\hat{T}, \hat{\eta}, \hat{\mu})$ on $A \downarrow \mathcal{C}$ such that $\mathcal{U} \circ \hat{T} = T \circ \mathcal{U}$.
\begin{proof}

Since $(T, \eta, \mu)$ is a monad, $(T, \eta)$ is a pointed endofunctor of $\mathcal{C}$: we can invoke \ref{Lifting pointed endofunctors} and construct a pointed endofunctor $(\hat{T}, \hat{\eta})$ of $A \downarrow \mathcal{C}$ satisfying the equation $\mathcal{U} \circ \hat{T} = T \circ \mathcal{U}$; we just need to enhance $(\hat{T}, \hat{\eta})$ to a monad $(\hat{T}, \hat{\eta}, \hat{\mu})$.

Start by defining $\hat{\mu}_{(X, \phi)} := \mu_X$; to see that this defines a morphism in $A \downarrow \mathcal{C}$, we need to check that
\[\mu_X \circ T^{\circ 2}(\phi) \circ T(\eta_A) \circ \eta_A = T(\phi) \circ \eta_A\]
that is, the following diagram commutes
\begin{center}\begin{tikzcd}
A \arrow[r, "\eta_A"] \arrow[d, "\eta_A"] & T(A) \arrow[dl, equals] \arrow[r, "T(\phi)"] & T(X) \\
T(A) \arrow[r, "T(\eta_A)"'] & T^{\circ 2}(A) \arrow[u, "\mu_A"'] \arrow[r, "T^{\circ 2}(\phi)"'] & T^{\circ 2}(X) \arrow[u, "\mu_X"']
\end{tikzcd}\end{center}
The triangle on the upper left corner commutes by definition, the other one by unitality for $(T, \eta, \mu)$. The square commutes by naturality of $\mu$.

Finally, the proof that $(\hat{T}, \hat{\eta}, \hat{\mu})$ is a monad trivially reduces to the statement that $(T, \eta, \mu)$ is a monad, which we assumed.
\end{proof}
\end{lm}

As noted in section 2, such lifts correspond 1-1 to distributive laws; we now unroll the constructions and describe such distributive laws explicitly, as well as showing how they can be made ``parametric''

\begin{thm}\thmheader{Parametric distributive law for Either}

Let $\mathcal{C}$ be a cocartesian monoidal category, and $A : \mathcal{C}$. Define $\mathcal{B} := \mathbb{B} \mathbf{Cat}_*[\mathcal{C}, \mathcal{C}]$ to be the 2-category obtained by delooping the monoidal category of pointed endofunctors of $\mathcal{C}$, and $\mathcal{D} := \mathbf{Mnd}(\mathcal{B})$ (which are just monads in $\mathcal{C}$, since the point is required to agree with the monad's unit).

There is an object in $PDist(\mathcal{D}, \mathcal{B})$, that is a strict 2-functor
\[(\mathbf{Mnd}(\mathbb{B} \mathbf{Cat}_*[\mathcal{C}, \mathcal{C}]) \otimes_l \mathbb{B}\Delta_\alpha) \otimes_l \mathbb{B} \Delta_\alpha \to \mathbb{B} \mathbf{Cat}_*[\mathcal{C}, \mathcal{C}]\]
such that the underlying monads, encoded as strict 2-functors
\[\mathbf{Mnd}(\mathbb{B} \mathbf{Cat}_*[\mathcal{B}, \mathcal{B}]) \to \mathbf{Mnd}(\mathbb{B} \mathbf{Cat}_*[\mathcal{B}, \mathcal{B}])\]
\begin{proof}

Such a 2-functor corresponds to a 2-functor
\[\mathbf{DEither}_X : \mathbf{Mnd}(\mathbb{B}\mathbf{Cat}[\mathcal{C}, \mathcal{C}]) \to \mathbf{Dist}(\mathbb{B}\mathbf{Cat}[\mathcal{C}, \mathcal{C}])\]
defined in a very similar way to \ref{Parametric distributive law for the Writer monad}: we employ lemma \ref{Lifting monads} for its action on objects and lemma \ref{Lifting pointed endofunctors} for its action on 1-cells; the action on 2-cells is trivial.
\end{proof}
\end{thm}

\printbibliography

\end{document}